\documentclass[a4paper, 11pt]{article}
\usepackage[latin1]{inputenc}
\usepackage[francais]{babel}
\usepackage{epsfig}
\usepackage{amsmath}
\usepackage{amssymb}
\usepackage{latexsym}

\newtheorem{thm}{Th{\'e}or{\`e}me}[section]
\newtheorem{prop}[thm]{Proposition}
\newtheorem{cor}[thm]{Corollaire}
\newtheorem{lem}[thm]{Lemme}
\newtheorem{defn}{D{\'e}finition}[section]
\newtheorem{conj}[thm]{Conjecture}

\newcommand{\ie}{{\it i.e.}}
\newcommand{\cf}{{\it confer} }

\date{}
\author{ Aurélien Galateau }
\begin{document}
\title{ Le problème de Bogomolov effectif sur les variétés abéliennes}

\newenvironment{dem}{\textbf{Preuve}\par}
{\begin{flushright}$\Box$\end{flushright}}

\maketitle

\bigskip

{\small {\sc Résumé.} On obtient une nouvelle minoration du minimum essentiel en petite codimension sur les variétés abéliennes, sous une conjecture concernant leurs idéaux premiers ordinaires. Cette minoration, déjà connue dans le cas torique depuis les travaux d'Amoroso et David, est optimale ``à $\epsilon$ près'' en le degré de la sous-variété. 

\bigskip

{\sc Abstract.} We give a new lower bound for the essential minimum of subvarieties of abelian varieties with small codimension, under a conjecture about ordinary primes in abelian varieties. This lower bound is already known in the toric case since the work of Amoroso and David. It is the best expected, ``up to an $\epsilon$'', in the degree of the subvariety.}  

\bigskip

\section{Introduction}

On souhaite dans ce travail obtenir une minoration du minimum essentiel sur les variétés abéliennes. Une telle minoration est une version quantitative de la conjecture de Bogomolov dont la formulation originale concerne les courbes algébriques.

\bigskip

Si $C$ est une courbe algébrique de genre $g \geq 2$ définie sur $\overline{\mathbb{Q}}$ et plongée dans sa jacobienne $J(C)$, on note $\hat{h}$ la hauteur canonique sur $J(C)$. On a alors la conjecture suivante, énoncée par Bogomolov en 1981 puis démontrée par Ullmo (voir \cite{Ullmo96}):
\begin{thm}
Il existe $\epsilon > 0$ tel que $\{ x \in C(\overline{\mathbb{Q}}), \hat{h}(x) \leq \epsilon \}$ est fini.
\end{thm} 
Puisque les points de torsion sont exactement ceux de hauteur nulle, le théorème d'Ullmo généralise le résultat suivant, connu sous le nom de {\it conjecture de Manin-Mumford} et prouvé par Raynaud (dans \cite{Raynaud83}): 
\begin{thm}
Les points de torsion de $J(C)$ qui sont dans $C(\overline{\mathbb{Q}})$ sont en nombre fini. 
\end{thm}

Plus généralement, soit $V$ une sous-variété algébrique d'une variété abélienne munie d'un fibré ample et symétrique, et $\hat{h}$ la hauteur de Néron-Tate associée à ce fibré. On commence par donner un analogue en dimension supérieure de l'hypothèse faite précédemment sur le genre:
\begin{defn}
On dit que $V$ est de torsion si $V$ est la translatée d'une sous-variété abélienne par un point de torsion.
\end{defn}
Une courbe algébrique de torsion est en particulier de genre $1$. On introduit par ailleurs le minimum essentiel, pour décrire les points de petite hauteur dans $V$:

\newpage

\begin{defn}
Le minimum essentiel de $V$ est:
\begin{displaymath}
\hat{\mu}_{\mathrm{ess}}(V)= \mathrm{inf} \{ \theta > 0,
\overline{V(\theta)}^Z= V(\overline{\mathbb{Q}}) \},
\end{displaymath}
\begin{center}
où $V(\theta)= \{ x \in V(\overline{\mathbb{Q}}), \hat{h}(x) \leq
\theta \}$,
\end{center}
et $\overline{V(\theta)}^Z$ est son adhérence de Zariski.
\end{defn}
On a alors la généralisation suivante du théorème $1$, démontrée par Zhang (\cf  \cite{Zhang98}):
\begin{thm}
Soit $V$ une sous-variété propre d'une variété abélienne $A$.
Le minimum essentiel de $V$ est nul si et seulement si $V$ est de torsion.
\end{thm}
\textbf{Remarque} Ici, ``$V$ est une sous-variété propre de $A$'' signifie que $V \subsetneq A$. \\

Le résultat analogue est vrai si on remplace $A$ par un tore (\cf \cite{Zhang92}) ou plus généralement par une variété semi-abélienne (\cf \cite{David-Philippon00}).

\bigskip

On peut chercher à obtenir une version quantitative de ce
résultat, en précisant $\epsilon$ dans le théorème 1. Ceci
revient, en dimension générale, à minorer le minimum essentiel d'une variété qui n'est pas
de torsion.  Grâce au théorème des minima successifs démontré par
Zhang (dans \cite{Zhang-Positive95}), il est équivalent de minorer la hauteur d'une telle variété.
Depuis les travaux de Bombieri et Zannier (voir \cite{Bombieri-Zannier95} pour le cas torique et \cite{Bombieri-Zannier96} pour le cas abélien), on sait qu'on peut espérer obtenir une borne ``uniforme'' pour le minimum essentiel, ne dépendant que du degré de $V$ et de la variété abélienne $A$.

\bigskip

Amoroso et David obtiennent une majoration optimale aux termes logarithmiques près en le degré
de $V$ pour les sous-variétés d'un tore (voir  \cite{Amoroso-David03}). Le degré y est remplacé
par un invariant plus fin qui apparaît naturellement avec les techniques diophantiennes, l'indice
d'obstruction. 

\begin{defn}
Soit $V$ une sous-variété algébrique propre et irréductible de $S$ une variété semi abélienne munie d'un fibré ample.
On appelle indice d'obstruction de $V$, noté $\omega( V)$:
\begin{displaymath}
\omega( V)= \mathrm{inf} \{
\mathrm{deg}(Z) \},
\end{displaymath}
où l'infimum est pris sur l'ensemble des hypersurfaces irréductibles de $S$ contenant $V$.
\end{defn}
Par le plongement standard $\mathbb{G}_m^n \hookrightarrow \mathbb{P}^n$, on obtient une hauteur projective $h$ sur les points de $\mathbb{G}_m^n$, et un minimum essentiel $\hat{\mu}_{\mathrm{ess}}$ sur les sous-variétés de $\mathbb{G}_m^n$. Amoroso et David démontrent la minoration suivante:
\begin{thm}
\label{torique}
Soit $V$ une sous-variété propre (et irréductible) de $\mathbb{G}_{m}^{n}$ de codimension $r$ qui n'est contenue dans aucun translaté d'un sous-tore propre de $\mathbb{G}_{m}^{n}$. On a alors:
\begin{center}
$\hat{\mu}_{\mathrm{ess}}(V) \geq \frac{c(n)}{\omega(V)} \times (\mathrm{log}(3 \omega(V)))^{-\lambda(r)},$
\end{center}
où $c(n)$ est un réel strictement positif et $\lambda(r)=(9(3r)^{(r+1)})^{r}$.
\end{thm}

\bigskip

Dans le cas des variétés abéliennes, on dispose déjà de résultats quantitatifs et inconditionnels, mais la
dépendance en le degré n'est pas aussi bonne. On a (\cf \cite{David-Philippon02}):

\begin{thm}
\label{theoreme da-phi}
Soit $A$ une variété abélienne de genre $g \geq 2$ définie sur $\overline{\mathbb{Q}}$, principalement polarisée par un fibré $\mathcal{M}$, et $V$ une sous-variété algébrique de $A$ qui n'est pas translatée d'une sous-variété abélienne. Alors:
\begin{displaymath}
\hat{\mu}_{\mathrm{ess}}(V) \geq \frac{\mathrm{min} \{ 1;
\mathcal{R}_{\mathrm{inj}} \}^{2(b+1)}}{2^{11g^{3}} (g-k+1)
\mathrm{deg}(V)^{2k(b+1)}},
\end{displaymath}
où $k$ désigne le nombre minimal de copies de $V-V$ dont la somme est une
sous-variété abélienne de $A$ , $b$ la dimension de cette sous-variété abélienne
et $\mathcal{R}_{\mathrm{inj}}$ la plus petite norme de Riemann d'une période
d'une conjuguée de $A$.
\end{thm}
Le terme au numérateur, appelé rayon d'injectivité, est relié au terme de hauteur $h(A)$ (hauteur projective de l'origine dans le plongement associé à $\mathcal{M}^{\otimes 16}$) par le lemme ``matriciel'' de Masser (voir le lemme 6.8 de \cite{David-Philippon02}). Cette minoration est monomiale inverse en le degré, alors que dans le cas torique, elle est linéaire inverse en l'indice d'obstruction (aux termes logarithmiques près), ce qui
 correspond à une minoration en deg$(V)^{-1/ \mathrm{codim}(V)}$.

Remarquons enfin que l'hypothèse du théorème \ref{theoreme da-phi} ($V$ n'est pas un translaté de sous-variété abélienne propre) est plus faible que son analogue torique dans le théorème \ref{torique} ($V$ n'est pas {\it incluse} dans un translaté de sous-tore propre); cette différence se ressent dès qu'on obtient des résultats comparables au théorème \ref{torique} et on peut préciser la minoration sous l'hypothèse faible, en faisant intervenir la dimension du plus petit translaté de sous-tore propre contenant $V$ (voir le corollaire 1.6 de \cite{Amoroso-David03}). 

\subsection*{Résultats}

On cherche à obtenir une minoration, pour les sous-variétés de variétés abéliennes, comparable à celle connue dans le cas torique. Soit donc $A$ une variété abélienne définie sur $K$ un corps de nombres, et $L$ un fibré ample et symétrique sur $A$. Soit de plus $\mathcal{A}$ un modèle entier de $A$ sur $\mathcal{O}_K$. Les méthodes employées dans un travail antérieur (correspondant au premier chapitre de \cite{Galateau07}) laissent espérer  une minoration en codimension $r \leq 2$ sous l'hypothèse suivante (on renvoie à {\it infra}, \ref{ordinaire} pour plus de détails sur la réduction ordinaire, et à \ref{densite}, sur la définition de la densité, qui est la densité {\it naturelle}):\\
\\
\textbf{Hypothèse H1} Il existe une densité positive d'idéaux premiers $\mathfrak{p}$ en lesquels la fibre spéciale $\mathcal{A}_{\mathfrak{p}}$ est ordinaire. 

\bigskip

Pour obtenir un résultat plus large, on a essayé de travailler avec des idéaux premiers ayant un autre type de réduction. Plusieurs obstructions sont apparues, nous forçant à faire l'hypothèse suivante sur $\mathcal{A}$ (voir {\it infra}, \ref{hassewitt} pour plus de détails):\\
\\
\textbf{Hypothèse H2} Il existe une densité positive d'idéaux premiers $\mathfrak{p}$ pour lesquels le $p$-rang de la fibre spéciale $\mathcal{A}_{\mathfrak{p}}$ est égal à $0$ ou à $g$.

\bigskip

Il a été d'abord nécessaire de supposer que le $p$-rang est égal au rang de Hasse-Witt (ou est nul) pour trouver de bonnes propriétés métriques $p$-adiques pour $A$ ($p$ étant un premier de $\mathbb{Z}$), reliées au type de réduction modulo $\mathfrak{p}$, pour $ \mathfrak{p}$ un idéal premier de $\mathcal{O}_K$ divisant $p$.

\bigskip

La preuve proprement dite, de nature diophantienne, conduit ensuite à travailler avec des idéaux premiers $\mathfrak{p}$ pour lesquels la fibre $\mathcal{A}_{\mathfrak{p}}$ a un $p$-rang égal à $0$ ou $g$. En effet, lorsque l'idéal premier $ \mathfrak{p}$ de $ \mathcal{O}_K$ n'est plus à réduction ordinaire, la propriété métrique obtenue est plus faible et on a besoin, en guise de compensation, d'un grand nombre de points de torsion se réduisant sur $0$ modulo $\mathfrak{q}$ (pour $\mathfrak{q}$ divisant $\mathfrak{p}$ dans une extension idoine). Ceci est vérifié si le $p$-rang est égal à $0$.  

\newpage

On démontre le résultat conditionnel:

\begin{thm}
\label{theoreme}
Soit $A$ une variété abélienne définie sur $K$ un corps de nombre. On suppose qu'il existe un modèle entier $\mathcal{A}$ de $A$ vérifiant \textbf{H2}. Alors on a la propriété $\textbf{P}(A)$ suivante: pour toute sous-variété algébrique $V$ propre (et irréductible) de codimension $r \leq 2$ dans $A$, si $V$ n'est pas contenue dans le translaté
d'une sous-variété abélienne propre de $A$, on a:
\begin{displaymath}
\hat{\mu}_{\mathrm{ess}}(V) \geq \frac{C(A)}{\omega(V)} \times
(\mathrm{log}(3 \mathrm{deg}(V)))^{-\lambda(r)},
\end{displaymath}
où $C(A)$ est un réel strictement positif ne dépendant que de $A$ et où $\lambda(r)=(16(2r)^{(r+1)})^{r}$.
\end{thm}
La constante $C(A)$ est explicitable, et d'autant plus dans le contexte des pentes. Cependant, la méthode employée est coûteuse en terme de la hauteur de $A$ et ne peut pas égaler les meilleures minorations obtenues par David et Philippon.  

Les produits de courbes elliptiques, tout comme les variétés abéliennes CM, vérifient la propriété \textbf{H1}, qui implique clairement \textbf{H2}.
L'hypothèse \textbf{H1} est l'objet de la conjecture ``folklorique'' suivante:
\begin{conj}
\label{folklorique}
Pour toute variété abélienne $A$ définie sur $K$ un corps de nombres, quitte à étendre $K$, il existe un modèle entier $\mathcal{A}$ de $A$ vérifiant \textbf{H1}.
\end{conj}
Sous cette conjecture, la propriété \textbf{P}(A) est vérifiée pour toute variété abélienne $A$ définie sur $K$ un corps de nombres. Remarquons qu'on peut même conjecturer, comme le fait Pink (dans \cite{Pink98}, 7), que les premiers ordinaires sont en densité $1$. 

En dimension $1$, pour une courbe elliptique $E$, le résultat est connu. Plus précisément, on sait que la densité de tels idéaux est $1$ si $E$ n'est pas CM (voir \cite{Serre68}, IV, 13), au moins $1/2$ si elle est CM. 

La validité de \textbf{H1} a été étendue aux surfaces abéliennes par les travaux de Katz et Ogus (voir \cite{Ogus83} 2.7, en remarquant par le théorème de Chebotarev que les premiers de degré $1$ ont une densité positive). Le principal résultat de cet article s'applique alors sans restriction:
\begin{cor}
Soit $C$ une courbe algébrique de genre $2$, incluse dans une surface abélienne $A$. Alors on a:
\begin{displaymath}
\hat{\mu}_{\mathrm{ess}}(C) \geq \frac{c(A)}{\mathrm{deg}(C)} \times
\big(\mathrm{log}(3 \mathrm{deg}(C))\big)^{-64},
\end{displaymath}
où $c(A)$ est un réel strictement positif ne dépendant que de $A$.
\end{cor}
\textbf{Remarque} Il faudrait être beaucoup plus précis sur la constante $c(A)$ pour donner un énoncé significatif avec une courbe de genre $2$ plongée dans sa jacobienne.

\bigskip

Des conditions suffisantes pour que \textbf{H1} soit réalisée, portant sur les groupes de monodromie $G_l$ (associés à chaque nombre premier $l$) de la variété abélienne, ont été données par Noot (voir \cite{Noot95}, 2), puis Pink (\cite{Pink98}, 7).

\bigskip

Rappelons que la minoration fine du minimum essentiel permet d'obtenir des résultats en direction des conjectures formulées par Zilber sur les variétés semi-abéliennes (dans \cite{Zilber02}), puis Pink sur les variétés de Shimura mixtes (voir \cite{Pink05}, conjectures 1.2 et 1.3). Pour $S$ un sous-ensemble de $\mathbb{G}_m^n$, on note:
\begin{displaymath}
S_{\epsilon}=\{ xy, x \in S, y \in \mathbb{G}_m^n, h(y) \leq \epsilon \}.
\end{displaymath}
En utilisant le théorème \ref{torique}, Habegger a démontré le résultat suivant (voir \cite{Habegger06}):  
\begin{thm}
\label{habegger}
Soit $C$ une courbe algébrique dans $\mathbb{G}_m^n$ qui n'est pas incluse dans le translaté d'un sous-tore propre, alors il existe $\epsilon >0$ tel que $C \cap \mathcal{H}_{\epsilon}$ est fini, où: 
\begin{displaymath}
\mathcal{H}= \bigcup_{\mathrm{codim} H=2} H,
\end{displaymath}
la réunion portant sur tous les sous-groupes algébriques de $\mathbb{G}_m^n$ ayant la codimension prescrite.
\end{thm}
Ce théorème généralise à la fois le théorème 2 de \cite{Bombieri-Masser-Zannier99} et la propriété de Bogomolov pour les courbes plongées dans les tores. Récemment, Maurin a démontré la conjecture de Zilber pour une courbe plongée dans un tore, en utilisant le théorème \ref{habegger} et une inégalité de Vojta uniforme. Plus précisément, si $C$ est une courbe et $S$ un sous-ensemble de $\mathbb{G}_m^n$, notons:
\begin{displaymath}
E(C,S)=C \cap \bigcup_{\mathrm{codim} (B)=2} S \cdot B, 
\end{displaymath} 
où la réunion porte sur les sous-tores de codimension indiquée. Maurin prouve d'abord (théorème 1.5 de \cite{Maurin07}):
\begin{thm}
Soit $C$ une courbe algébrique de $\mathbb{G}_m^n$ qui n'est pas incluse dans le translaté d'un sous-tore propre et $\Gamma$ un sous-groupe de rang fini de $\mathbb{G}_m^n$. Alors il existe un réel $\epsilon > 0$ tel que l'ensemble $E(C, \Gamma_{\epsilon})$ soit fini. 
\end{thm}
Et il en déduit:
\begin{cor}
Soit $C$ une courbe algébrique irréductible de $\mathbb{G}_m^n$ non incluse dans un sous-groupe algébrique propre (pas nécessairement irréductible). Alors $C \cap \mathcal{H}$ est fini.
\end{cor}
Ce corollaire optimise le résultat principal de \cite{Bombieri-Masser-Zannier99}, qui suppose que $C$ n'est pas incluse dans un {\it translaté de sous-tore propre}. 

Le théorème \ref{theoreme} est donc la première étape, sous la conjecture \ref{folklorique}, d'un programme qui permet d'attaquer la conjecture de Zilber-Pink sur les variétés abéliennes.

\subsection*{Plan de l'article}

On démontre le théorème \ref{theoreme} en utilisant la méthode des pentes, formalisée par Bost dans \cite{Bost95}. La deuxième partie est consacrée à la démonstration d'une propriété $p$-adique obtenue par l'étude du groupe formel d'une variété abélienne en caractéristique $p$. On commence par faire quelques rappels sur le $p$-rang d'une variété abélienne, puis sur la théorie des schémas en groupes. On obtient ensuite un résultat métrique $p$-adique précis, pour les points de $p$-torsion de $A$ se réduisant sur $0$ modulo un idéal premier $\mathfrak{q}$ divisant $p$ dans une extension convenable. Ce résultat suppose en principe le choix d'une base du tangent pour chaque idéal premier $\mathfrak{q}$ mais on démontre, à l'aide d'un argument de géométrie des nombres, que si on borne les premiers, il existe une base sur $\mathcal{O}_K$ de hauteur contrôlée dans laquelle toutes les propriétés $p$-adiques sont simultanéments vérifiées. On a cherché, dans cette partie, à obtenir les résultats les plus précis en fonction du $p$-rang. 

\bigskip

Dans la troisième partie, on rappelle les définitions et résultats généraux de la théorie des pentes. Un premier fait assez inhabituel dans notre application de cette théorie est l'importance des estimations ultramétriques. Dans cette perspective, on utilise une version du théorème des pentes assez précise sur le plan ultramétrique. Puis on introduit les fibrés hermitiens qui seront utiles par la suite et on estime leur pente. La principale difficulté de cette partie réside dans la majoration de la pente maximale du fibré des sections d'un fibré ample sur une sous-variété de $A$ (avec multiplicités), le fibré d'arrivée étant habituellement formé, dans la méthode des pentes, à partir d'un nombre fini de points. Cette majoration est obtenue en suivant une idée figurant dans la thèse de Chen (\cite{Chen06}). Les résultats de Bost et Künnemann ( \cite{Bost-Kunnemann06}, améliorés par Chen en dimension $\geq 3$ dans le chapitre 5 de sa thèse) concernant la pente maximale du produit tensoriel de deux fibrés hermitiens permettent de prendre en compte la multiplicité. 

\bigskip

La preuve du théorème commence réellement dans la quatrième partie. On prend une sous-variété propre $V$ d'une variété abélienne $A$ qui n'est pas incluse dans un translaté de sous-variété abélienne et on lui associe un fermé de Zariski $X$ en vue de la fin de la preuve; on construit deux espaces vectoriels $E$ et $F$ et un morphisme de restriction entre ces espaces (paramétrés en fonction de l'indice d'obstruction de $X$ et en fonction de $A$). Puis on fixe les paramètres (degré de l'espace de sections, multiplicités, bornes pour la norme des idéaux premiers utilisés) intervenant dans cette construction et on suppose par l'absurde que le minimum essentiel de $X$ est majoré en fonction de ces paramètres. Dans toute cette partie et les suivantes, on travaille avec un plongement étiré, devenu classique dans les travaux diophantiens sur les variétés abéliennes pour passer de la hauteur projective à la hauteur de Néron-Tate.

\bigskip

Dans la partie suivante, on calcule les rangs et les normes des morphismes susceptibles de rentrer dans l'inégalité des pentes. On écrit ensuite cette inégalité, sous l'hypothèse que le morphisme soit injectif. On parvient rapidement à une contradiction. A ce stade du travail, on a montré que le minimum essentiel de $X$ est correctement minoré {\it modulo l'injectivité du morphisme}.

\bigskip

On suppose donc par l'absurde, dans la sixième partie, que le morphisme n'est pas injectif. On commence par appliquer un lemme de zéros très général d'Amoroso et David. L'utilisation de ce lemme est suivie, comme dans les cas torique et multi-elliptique, d'une phase de dénombrement et d'un argument de descente, qui permettent d'obtenir une contradiction. Le travail sur l'injectivité du morphisme s'effectue {\it exceptionnellement} après l'utilisation de l'inégalité des pentes, d'abord car il est assez long, mais surtout parce qu'il comporte une itération (dans la phase de descente) qui nécessite d'avoir déjà écrit cette inégalité. On n'a pu adapter la phase de descente qu'en petite codimension: $r \leq 2$.

\subsection*{Constantes}

Le théorème \ref{theoreme} montre l'existence d'une constante $C(A)$ ne dépendant que de $A$ et impliquée dans la minoration du minimum essentiel. Au cours de ce travail, on introduira des constantes $c_1, \ldots, c_{22}$ ne dépendant que de $A$. Le choix des paramètres fera intervenir une constante $C_0$, dépendant uniquement de $A$ elle aussi, qui sera prise {\it grande} par rapport aux constantes $c_i$ ($1 \leq i \leq 22$). La constante $C(A)$ s'exprimera alors simplement en fonction de $C_0$.

\subsection*{Remerciements}

Je souhaite ici remercier chaleureusement Sinnou David, qui m'a patiemment initié au problème de Bogomolov au cours de mon doctorat, et Eric Gaudron pour sa minutieuse relecture et ses conseils. Certains points de ce travail se sont aussi améliorés grâce aux éclairages d'Antoine Chambert-Loir, Huayi Chen, Richard Pink, Hugues Randriam et Emmanuel Ullmo. 

\section[Lemme $p$-adique]{Un lemme clé $p$-adique sur les variétés abéliennes}
\label{partie metrique}

Soit $A$ une variété abélienne définie sur $K$ un corps de nombres et soit $\mathcal{P}_A$ un ensemble de premiers bien choisi en fonction de $A$. 
Le but de cette partie est d'établir une inégalité $p$-adique concernant les points de $p$-torsion de $A$, pour $p \in \mathcal{P}_A$. On cherche à montrer, pour les variétés abéliennes, un résultat comparable à l'inégalité suivante, vraie pour tout premier $p$, toute racine $p$-ème de l'unité $\xi$ et toute place $v/p$ d'un corps de nombres quelconque contenant $\xi$: 

\begin{displaymath}
|\xi -1|_v \leq p^{-1/p}.
\end{displaymath}

\bigskip

Cette inégalité, conséquence d'une propriété de ramification bien connue sur les corps cyclotomiques (voir \cite{Ireland-Rosen80}, Proposition 13.2.7), a un analogue satisfaisant sur  les courbes elliptiques pour les premiers de bonne réduction ordinaire, s'il n'y a pas de ramification initiale; dans le cas des premiers supersinguliers, on doit remplacer $1/p$ par $1/p²$ (voir \cite{Galateau07}, 2.4.1). On s'attend donc à ce que le résultat obtenu sur $A$ dépende de la réduction de $A$ modulo $\mathfrak{p}$, un idéal de $\mathcal{O}_K$. Deux invariants associés à une variété abélienne en caractéristique positive apparaissent naturellement au cours de la discussion, à savoir le rang de la matrice de Hasse-Witt et son rang stable (ou $p$-rang de la variété). On obtient une inégalité assez précise sous l'hypothèse que ces deux invariants sont égaux, puis on trouve une base du tangent dans laquelle les propriétés $p$-adiques associées à des premiers différents (et de norme bornée) sont vraies simultanément. 

\bigskip

Dans la suite de ce travail, on n'utilisera pas ces estimations $p$-adiques dans toute leur précision. Dans ces conditions, le choix de la base adaptée n'est pas capital; il est par ailleurs peu coûteux puisqu'il n'ajoute qu'un terme négligeable dans le calcul de pente du sous-paragraphe \ref{pentesousvariete}. 

\bigskip

L'intérêt de ce résultat $p$-adique, dans la mise en \oe uvre de la méthode des pentes, sera d'estimer certaines normes ultramétriques d'un morphisme de restriction; cette estimation est le point crucial de la preuve et correspond à la phase d'extrapolation dans une preuve classique de transcendance. 

\subsection{Le $p$-rang d'une variété abélienne}
\label{ordinaire}

On fixe pour ce paragraphe et le suivant un premier $p$ et une variété abélienne $A$ définie sur un corps fini $\mathbb{F}_q \subset k= \overline{\mathbb{F}}_p$. Commençons par donner la définition du $p$-rang, à travers la proposition suivante (\cf \cite{Mumford}, page 64): 

\begin{prop}
Le sous-groupe $A[p]$ des points de $p$-torsion sur $\overline{\mathbb{F}}_p$ est de cardinal $p^{\alpha}$, où $\alpha \in [0;g]$ est un entier. De plus, le groupe $A[p]$ est isomorphe à $(\mathbb{Z}/p\mathbb{Z})^{\alpha}$.  Cet entier $\alpha$ sera appelé $p$-rang de $A$.
\end{prop}

\noindent
\textbf{Remarques} Le $p$-rang de $A$ est aussi appelé rang stable de $A$ (en raison de son lien avec la matrice de Hasse-Witt de $A$, explicité {\it infra}, \ref{hassewitt}). S'il est égal à $g$, la variété $A$ est dite {\it ordinaire}. 

\bigskip

On peut maintenant introduire les morphismes de Frobenius $\Phi_{p^j}$, pour $j \geq 1$ un entier. On définit d'abord $\Phi_{p^j}$ sur Spec$(k)$ comme étant l'identité sur l'espace topologique réduit à un point et l'élévation à la puissance $p^j$ sur $k$. On note ensuite $A^{(p^j)}$ le schéma sur Spec$(k)$ défini comme le {\it tiré en arrière} de $A$ par l'action du Frobenius $\Phi_{p^j}$ sur Spec$(k)$. Par construction, ce schéma est une variété abélienne.
\begin{defn}
Le morphisme de Frobenius:
\begin{displaymath}
\Phi_{p^j}: A \rightarrow A^{(p^j)}.
\end{displaymath}
est défini par l'élévation à la puissance $p^j$ sur le faisceau structural. 
\end{defn}
\textbf{Remarque} Si la variété abélienne $A$ est définie sur $\mathbb{F}_q$, avec $q=p^j$, alors le Frobenius $\Phi_{p^j}$ est un morphisme de $A$ dans $A$. 

\bigskip

Le schéma $A^{(p^j)}$ est une variété abélienne et le Frobenius $\Phi_{p^j}$ une isogénie purement inséparable de degré $p^{jg}$ (voir \cite{Moonen}, page 72). On vient de montrer que le degré inséparable de $[p]$ est supérieur à $p^g$. Dans le cas des courbes elliptiques, il est facile de prouver que $[p]$ se factorise à travers $\Phi_p$ ou $\Phi_{p^2}$ car son degré inséparable est exactement le degré d'un de ces deux morphismes. La factorisation par $\Phi_p$ est vraie en général, et on l'obtient en construisant explicitement le morphisme avec lequel on compose le Frobenius, appelé ``Verschiebung" (\cf \cite{Moonen}, pages 74 et 104):

\begin{lem}
Il existe une isogénie $V: A^{(p)} \rightarrow A$ telle que $[p]=V \circ \Phi_p$. De plus, $V$ et $\Phi_p$ sont duales l'une de l'autre au sens suivant: si on note $\widehat{A}$ la duale de $A$, on a une décomposition:$[p]_{\widehat{A}}= V_{\widehat{A}} \circ \Phi_{p, \widehat{A}}$ avec: 
\begin{displaymath}
\widehat{V}= \Phi_{p,\widehat{A}} \mathrm{\ et \ } \widehat{\Phi_p}= V_{\widehat{A}}.
\end{displaymath}
\end{lem}
Comme le morphisme de Frobenius est purement inséparable, le $p$-rang n'est autre que le degré séparable de $V$. Si $\alpha=g$, l'isogénie $V$ est séparable et sa différentielle en $0$ est inversible. On veut relier plus généralement $\alpha$ à la différentielle de $V$ en $0$; ceci nous amène à étudier la structure de {\it schéma en groupe} du sous-groupe $A[p]$ de $A$. 

\subsection{Schémas en groupe}
\label{hassewitt}

On fait ici quelques rappels sur la théorie des schémas en groupes; on peut trouver ces résultats dans \cite{Mumford}, partie III. 

\begin{defn}
Un schéma en groupe $G$ sur $k$ est un schéma muni d'un morphisme de multiplication $m: G \times G \rightarrow G$, d'un morphisme d'inversion $i: G \rightarrow G$ et d'un élément neutre $e: \mathrm{Spec} (k) \rightarrow G$ vérifiant les axiomes:

\begin{center}
$m \circ (m \times \mathrm{Id}_G)  =  m \circ(\mathrm{Id}_G \times m): G \times G \times G \rightarrow G,$ 
\end{center}
\begin{center}
$m \circ (e \times \mathrm{Id}_G)  =  j_1:\mathrm{Spec} (k) \times G \rightarrow G,$
\end{center} 
\begin{center}
$m \circ (\mathrm{Id}_G \times e)  =  j_2:G \times \mathrm{Spec}(k) \rightarrow G,$
\end{center}
et:
\begin{center}
$e \circ \pi=m \circ(\mathrm{Id}_G \times i) =m \circ (i \times \mathrm{Id}_G): G \rightarrow G$,
\end{center}
où $\pi: G \rightarrow \mathrm{Spec} (k)$; $j_1: \mathrm{Spec} (k) \times G \simeq G$ et $j_2:G \times \mathrm{Spec} (k) $ sont les isomorphismes canoniques.
\end{defn}

Soit $G$ un schéma en groupe. Son algèbre de Lie est le $k$-espace vectoriel des champs de vecteurs invariants par $m$ et elle est munie de la fonction de Hasse-Witt, qui associe à une dérivation $D$ la dérivation $D^p$ ($p$-ème itérée de $D$). C'est une application $\mathbb{F}_p$-linéaire (\ie \ additive et linéaire sous la multiplication par un élément de $\mathbb{F}_p$). 

On  définit $\widehat {G}$ le dual (de Cartier) de $G$ d'un schéma en groupe affine Spec$(R)$ en prenant le dual $R^*$ de $R$ et en le munissant d'une comultiplication et d'un idéal d'augmentation par dualité. 

On suppose maintenant que $G$ est un schéma en groupe fini et commutatif. On dit que $G$ est de type $l$ (resp. de type $r$) si l'espace sous-jacent est constitué d'un seul point (resp. si $G$ est réduit). On dit que $G$ est de type $(x,y)$ si $G$ est de type $x$ et $\widehat{G}$ est de type $y$. Le schéma $G$ se décompose alors de façon unique en un produit:
\begin{displaymath}
G=G_{r,r} \times G_{r,l} \times G_{l,r} \times G_{l,l},
\end{displaymath}
où $G_{x,y}$ est de type $(x,y)$ (pour plus de détails, voir \cite{Mumford}, §14). 

\bigskip

Pour le schéma en groupe $A[p]$ qui nous intéresse ici, le type $G_{r,r}$ est trivial car $A[p]$ est de cardinal une puissance de $p$. Plus précisément, on a la proposition:
\begin{prop}
On a l'isomorphisme de schémas en groupes:
\begin{displaymath}
A[p] \simeq (\mathbb{Z}/p\mathbb{Z})^{\alpha} \times (\mu_p)^{\alpha} \times G_1^0,
\end{displaymath}
où $\alpha$ est le $p$-rang de $A$, $\mu_p= \mathrm{Spec} \big(k[X]/(X^p-1)\big)$ et $G_1^0$ est de type $(l,l)$.
\end{prop}
\begin{dem}
Soit $n \in \mathbb{N}^*$. Compte tenu des structures de groupes de $A[p^n]$ et du dual $\widehat{A[p^n]}$, qui donnent les composantes réduites de ces deux schémas en groupes, et comme le dual du noyau de l'isogénie $[p]$ est le noyau de l'isogénie duale (\cf \cite{Mumford}, page 143), on a la décomposition: 
\begin{displaymath}
A[p^n]  \simeq (\mathbb{Z}/p^n \mathbb{Z})^{\alpha} \times \widehat{(\mathbb{Z}/p^n \mathbb{Z})}^{\beta} \times G_n^0,
\end{displaymath}
pour un certain entier $\beta$ et un schéma en groupe local $G_n^0$. 

\bigskip

L'algèbre de fonctions associée à $\mathbb{Z}/p^n\mathbb{Z}$ est son algèbre de groupe et, en notant $X$ l'évaluation en $1 \in \mathbb{Z}/p^n \mathbb{Z}$, on voit que l'algèbre duale est isomorphe à: 
\begin{displaymath}
\mathrm{Spec} \big(k[X]/(X^{p^n}-1)\big).
\end{displaymath}
On en déduit que: $\widehat{\mathbb{Z}/p^n \mathbb{Z}} \simeq \mu_{p^n}$. De plus, $\alpha$ et $\beta$ sont permutés par passage de $A$ à $\widehat{A}$ et puisqu'il existe une isogénie $f: A \rightarrow \widehat{A}$, en notant $K$ le cardinal de son noyau, on a: 
\begin{displaymath}
f(A[p^n]) \subset \widehat{A}[p^n] \ \  \mathrm{donc} \ \ \ p^{n \alpha} \leq K p^{n\beta}. 
\end{displaymath}
En faisant varier $n$, on obtient: $\alpha \leq \beta$. Mais comme la duale de $\widehat{A}$ est isomorphe à $A$ (\cf \cite{Mumford}, page 132), on obtient: $\alpha=\beta$. 
\end{dem}
\textbf{Remarque} Les schémas en groupes de type local-local sont les plus difficiles à comprendre. Pour plus de détails, utilisant la théorie des vecteurs de Witt, on renvoie par exemple à \cite{Pink-online} (en particulier §16 et §22). 

\bigskip

On peut maintenant faire le lien entre la différentielle de $V$ en $0$ et le $p$-rang:
\begin{prop}
Soit $\Psi$ la différentielle de $V$ en $0$. Alors le $p$-rang de $A$ est le rang de $\Psi^g$.  
\end{prop}
\begin{dem}
On passe aux algèbres de Lie dans la proposition précédente et on observe que l'application linéaire $[p]^*$ est la multiplication par $p$, donc est nulle sur Lie$(A)$. En se limitant à la partie locale en $0$, on a:
\begin{displaymath}
\mathrm{Lie}(A)=\mathrm{Lie}(A[p])= \mathrm{Lie}(\mu_p)^{\alpha} \oplus \mathrm{Lie} (G_1^0).
\end{displaymath}
En prenant comme base de $\mathrm{Lie}(\mu_p)$ la dérivation $X \partial/\partial X$, on observe que la fonction de Hasse-Witt est l'identité sur $\mathrm{Lie}(\mu_p)^{\alpha}$, alors qu'elle est nilpotente sur la partie locale-locale. De plus, par l'isomorphisme canonique: 
\begin{center}
Lie$\widehat{A} \simeq H^1(A, \mathcal{O}_A)$,
\end{center}
la fonction de Hasse-Witt correspond à l'application induite par le Frobenius sur $\mathcal{O}_A$ (\cf \cite{Mumford}, page 148), qui correspond par dualité à la différentielle de $V$ sur le tangent en $0$. Il existe donc une décomposition du tangent en $0$: $t_A=V_s+V_n$ laissée stable par $\Psi$ telle que $\Psi_{|V_s}$ soit un isomorphisme et $\Psi_{|V_n}$ soit nilpotente; de plus, l'espace vectoriel $V_s$ est de dimension $\alpha$. En itérant $g$ fois l'application $\Psi$, la partie nilpotente s'annule et on en déduit que le $p$-rang est le rang de $\Psi^g$.
\end{dem}
\textbf{Remarque} On appelle composante semi-simple de $\Psi$ l'espace vectoriel $V_s$. Cette composante semi-simple est l'image de $\Psi^g$, donc est définie sur $\mathbb{F}_q$. 

\bigskip

Le $p$-rang d'une variété abélienne n'est pas toujours égal au rang de la matrice de Hasse-Witt. Si on fixe $\alpha \leq g-1$, on peut même montrer (voir \cite{Koblitz75}, theorem 7, page 163) que sur l'espace de module des variétés abéliennes principalement polarisées de dimension $g$ avec structure de niveau fixée sur $k$, les variétés abéliennes ayant une matrice de Hasse-Witt de rang $g-1$ sont Zariski-denses dans le fermé des variétés abéliennes dont le $p$-rang est plus petit que $\alpha$. Si le $p$-rang est égal à $g-1$ ou $g$, il est automatiquement égal au rang de la matrice de Hasse-Witt.

\bigskip

En général, le $p$-rang et le rang de la matrice de Hasse-Witt sont distincts, et la partie nilpotente fait obstruction pour contrôler efficacement la norme $p$-adique de tous les paramètres. On devra donc par la suite travailler avec des idéaux premiers $\mathfrak{p}$ de bonne réduction tels que la variété abélienne modulo $\mathfrak{p}$ ait ces deux invariants égaux.\\
\\
\textbf{Hypothèse} On suppose maintenant que le $p$-rang de $A$ est égal à $0$ ou au rang de $\Psi$. \\
\\
On choisit un système de paramètres $(x_1, \ldots, x_g)$ associés à une base de différentielles invariantes, tels que $(x_1, \ldots, x_{\alpha})$ soit une base de $\mathrm{Im} \Psi^g$ (qui est égale à $\mathrm{Im} \Psi$ sauf éventuellement si le rang est nul). Notons $\widehat{\mathcal{O}}_{0,A}$ le groupe formel associé à $A$ en $0$ sur $\mathbb{F}_q$. On a (\cite{Hindry-Silverman00}, page 268):
\begin{displaymath}
\widehat{\mathcal{O}}_{0,A} \simeq \mathbb{F}_q[[x_1, \ldots, x_g]].
\end{displaymath}
On note $\textbf{V}=(V_1, \ldots, V_g)$ le $g$-uplet de séries formelles image de l'isogénie $V$ dans le groupe formel. On note aussi $\Phi_{\alpha}$ le morphisme de $\mathbb{F}_q[[x_1, \ldots, x_g]]$ qui agit sur les paramètres par: 
\begin{eqnarray*}
x_i \rightarrow x_i & \mathrm{si} & i \leq \alpha \\
x_i \rightarrow x_i^p & \mathrm{si} & i > \alpha. \\
\end{eqnarray*}

\begin{cor}
\label{factorisation}
Il existe un $g$-uplet de séries formelles $\textbf{U}=(U_1, \ldots, U_g)$ tel que $\textbf{V}$ se factorise: $\textbf{V}=\textbf{U} \circ \Phi_{\alpha}$ et tel que $d\textbf{U}$ soit inversible. 
\end{cor}
\begin{dem}
Si le $p$-rang est nul, on note $\mathbb{K}$ la clôture séparable de $V^* k(A)$ dans $k(A)$. Comme $k(A)^{p^2} \subset k(A)$ est purement inséparable, et par définition de la clôture séparable, on a l'inclusion: $k(A)^{p^2} \subset \mathbb{K}$ puis égalité en comparant les degrés de ces extensions. La factorisation sur les corps de fonctions induit une factorisation: $\textbf{V}=\textbf{U} \circ \Phi_{p^2}$, où $\textbf{U}$ est séparable.
Sa différentielle est donc inversible (voir \cite{Lang02}, VIII, proposition 5.5). 

\bigskip

Supposons maintenant $\alpha > 0$. Pour tout entier $ i \in [1,g]$, comme $V$ est une isogénie, la forme différentielle $dx_i^* V$ est encore une différentielle invariante, donc se décompose: 
\begin{displaymath}
dx_i^* V= \sum_{j=1}^g \alpha_{i,j} dx_j,
\end{displaymath}
où les $\alpha_{i,j}$ sont constants et donnés par la $i$-ème colonne de la matrice de $\Psi$ dans la base associée à $(x_1, \ldots, x_g)$. 
On en déduit par intégration que les seuls termes non-nuls dans les $V_i$ sont les termes linéaires ou des monômes en $(x_1^p, \ldots, x_g^p)$. Par choix de la base, les paramètres $(x_{\alpha+1}, \ldots, x_g)$ sont absents de la partie linéaire. On a donc bien la décomposition voulue. L'application $\Phi_{\alpha}$ est purement inséparable et son degré est le rang de: 
\begin{displaymath}
k[[x_1, \ldots, x_g]]/(x_{\alpha+1}^p, \ldots, x_g^p),
\end{displaymath}
donc:
\begin{displaymath}
\mathrm{deg} \Phi_{\alpha}= \mathrm{deg}_i \Phi_{\alpha}= p^{g - \alpha}.
\end{displaymath}
En comparant les degrés séparables et inséparables, on voit que $\textbf{U}$ est séparable, et que sa différentielle est inversible.
\end{dem}

\subsection{Retour en caractéristique nulle}
\label{densite}

Le but de ce paragraphe est de traduire la proposition précédente en un résultat $p$-adique pour les points de torsion d'une variété abélienne définie sur un corps de nombres. Soit donc $A$ une variété abélienne de dimension $g$ définie sur un corps de nombres $K$, munie d'un fibré $L$ ample et symétrique, et soit $\mathcal{A}$ un modèle entier de $A$ sur $\mathcal{O}_K$. On peut supposer (voir \cite{Hindry-Silverman00}, page 105), quitte à considérer $L^{\otimes 3}$, que le fibré $L$ est très ample (et projectivement normal). Rappelons que pour tout nombre premier $p$, il y a $p^{2g}$ points de $p$-torsion dans $A(\overline{K})$. Pour un idéal premier $\mathfrak{p}$ de $\mathcal{O}_K$ de bonne réduction divisant $p$, la fibre spéciale $\mathcal{A}_{\mathfrak{p}}$ ne contient plus que $p^{\alpha}$ points de $p$-torsion, où $\alpha$ est le $p$-rang de la fibre spéciale. 

\bigskip

Par la suite, on fera implicitement un certain nombre de choix sur $A$ (une base de sections globales pour $L$, une base de dérivations algébriques, une base d'ouverts affines) et par abus de langage, on dira qu'une constante ne dépend que de $A$ si elle dépend de $A$ et de ces choix. 

\bigskip

Précisons d'abord que si $\mathfrak{q}$ est un idéal premier de $K'$ une extension finie de $K$, dont la projection sur $\mathbb{Z}$ est $p$, on choisit la normalisation suivante pour la valuation $\mathfrak{q}$-adique:
\begin{center} 
$|p|_{\mathfrak{q}}=p^{-n_{\mathfrak{q}}}$, où $n_{\mathfrak{q}}$ est le degré local: $[K'_{\mathfrak{q}}: \mathbb{Q}_{p}]$. 
\end{center}
Cette normalisation permet d'écrire plus simplement la formule du produit et la hauteur d'un morphisme, qui intervient dans les inégalités de pentes. 

\bigskip

Les premiers de réduction ordinaire (\ie \ les premiers de bonne réduction pour lesquels le $p$-rang est égal à $g$) sont ceux pour lesquels les propriétés métriques sont les meilleures. Dans le cas d'une courbe elliptique $E$, on sait qu'ils sont de densité $1$ si $E$ n'est pas à multiplication complexe; et qu'ils sont de densité au moins $1/2$ si $E$ est à multiplication complexe. En dimension supérieure, on ne connait aucun résultat de densité comparable. On va donc travailler avec un ensemble de premiers ayant même type de réduction (pas nécessairement ordinaire) et de densité positive. On définit la densité naturelle pour les idéaux premiers de la façon suivante:
\begin{defn}
Soit $\mathcal{Q}$ un sous-ensemble de l'ensemble $\mathcal{P}$ des idéaux premiers de $\mathcal{O}_K$. On dit que $\mathcal{Q}$ a une densité naturelle $d$ si le quotient:
\begin{displaymath}
\frac{|\{\mathfrak{q} \in \mathcal{Q}, \textbf{N}(\mathfrak{q}) \leq x\}|}{|\{\mathfrak{p} \in \mathcal{P}, \textbf{N}(\mathfrak{p}) \leq x\}|} 
\end{displaymath}
tend vers $d$ quand $x \rightarrow \infty$.
\end{defn}
\textbf{Remarque} La fonction $\textbf{N}$ est la norme sur les idéaux (définie dans \cite{Samuel}: III, 5). 

\bigskip

Dans toute la discussion qui suit, on omet de préciser les ensembles indexateurs, qui sont toujours finis et dépendent de $A$.

\bigskip

La loi d'addition de $A$ est donnée sur chaque ouvert affine par des polynômes de bi-degré $(2,2)$ (\cf \cite{Lange-Ruppert85} ou \cite{David-Philippon02}, proposition 3.7) dont les coefficients sont de hauteur bornée uniquement en fonction de $A$. On a donc, si on note $(x_k)_k$ l'ensemble fini de ces coefficients: 
\begin{eqnarray}
\label{condition 1}
\forall k & : & |x_k|_{\mathfrak{p}} \leq 1, 
\end{eqnarray}
sauf pour un nombre fini d'idéaux premiers $\mathfrak{p}$ (ne dépendant que de $A$).
\bigskip

Fixons maintenant une base de dérivations algébriques $(\partial_1, \ldots, \partial_g)$ sur $A$. Quitte à prendre des idéaux premiers $\mathfrak{p}$ de $\mathcal{O}_K$ plus grands qu'une constante ne dépendant que de $A$, cette base de dérivations est encore une base de dérivations modulo $\mathfrak{p}$. De plus, on a (\cf \cite{David91}):
\begin{thm}
Pour toute fonction abélienne $f_i$ sur $A$ et pour toute dérivation $\partial_j$ :
\begin{displaymath}
\partial_j f_i= \sum_{(k,l)} y_{k,l}^{i,j} f_k f_l,
\end{displaymath}
où les $f_i$ sont une base de  fonctions abéliennes sur $A$. De plus, les $y^{i,j}_{k, l}$ sont des nombres algébriques, de hauteur bornée uniquement en fonction de $A$. 
\end{thm}
\begin{dem}
La preuve de \cite{David91}, théorème 4.1, est donnée pour une base bien spécifique, la base de dérivation de Shimura, sous l'hypothèse que la polarisation est principale. Le résultat obtenu est alors effectif. Ses arguments s'adaptent sans peine pour obtenir le résultat qualitatif dont on a besoin dans la généralité indiquée. Rappelons-en les étapes.

On commence par observer que, si on fixe deux fonctions thêta $(\theta_0, \theta_1)$ telles que: $f_i= \theta_1/ \theta_0$, on a:
\begin{displaymath}
\theta_0^2 \partial_j\left(\frac{\theta_1}{\theta_0}\right) \in \Gamma(A, L^{\otimes 2}). 
\end{displaymath}
Puis (comme $L$ est associé à un plongement projectivement normal), le morphisme: 
\begin{displaymath}
\Gamma(A, L)^{\otimes 2} \longrightarrow \Gamma (A, L^{\otimes 2}) 
\end{displaymath}
est surjectif.

On a donc l'écriture attendue, mais avec les $y_{k,l}^{i,j}$ dans $\mathbb{C}$, et comme la base de dérivations est algébrique, on en déduit que les $y_{k,l}^{i,j}$ le sont aussi. Ces coefficients étant en nombre fini, on peut trouver une borne pour leur hauteur ne dépendant que de $A$. 
\end{dem}
En appliquant ce théorème à une base de sections de $L$, il en résulte que pour tout idéal premier $\mathfrak{p}$ de $\mathcal{O}_K$ sauf un nombre fini (ne dépendant que de $A$): 
\begin{eqnarray}
\label{condition 2}
\forall (i, j, k, l) & : & |y_{k,l}^{i,j}|_\mathfrak{p}  \leq  1. 
\end{eqnarray}

\bigskip

Les premiers de mauvaise réduction pour $A$ sont en nombre fini, ainsi que les premiers de $\mathbb{Z}$ se ramifiant dans $\mathcal{O}_K$. On pose $\mathcal{P}_{A,0}$ l'ensemble des premiers $\mathfrak{p}$ de $\mathcal{O}_K$ de bonne réduction, vérifiant (\ref{condition 1}) et (\ref{condition 2}), tels que la base de dérivations algébriques soit encore une base sur $\mathcal{A}_{\mathfrak{p}}$, tels que si $(p)= \mathfrak{p} \cap \mathbb{Z}$, on a: $e_{\mathfrak{p}/p}=1$, et enfin, tels que: 
\begin{displaymath}
\sum_ {p \in \mathcal{P}_{A,0}} \frac{1}{p^2} \leq \frac{1}{3}. 
\end{displaymath}
Comme la même somme indexée par $\mathbb{N}^*$ converge, il suffit d'exclure un ensemble fini (absolu) de premiers pour que cette condition soit vérifiée.  L'ensemble $\mathcal{P}_{A,0}$ est de densité naturelle égale à $1$, et sa construction ne dépend que de $A$ et $K$.

\bigskip

On fait maintenant l'hypothèse suivante sur $\mathcal{A}$: \\
\\
\textbf{Hypothèse H3} Il y a une densité $c_0>0$ d'idéaux premiers $\mathfrak{p}$ pour lesquels:
\begin{itemize}
\item soit le $p$-rang de $\mathcal{A}_{\mathfrak{p}}$ est égal à $0$,
\item soit il est non-nul et égal au rang de la matrice de Hasse-Witt.
\end{itemize} 

\bigskip

Par le principe des tiroirs de Dirichlet, il existe un entier $k$ tel que la densité naturelle d'idéaux premiers de $\mathcal{P}_{A,0}$ vérifiant \textbf{H3} et pour lesquels la fibre spéciale a un $p$-rang égal à $k$ est supérieure ou égale à $\frac{c_0}{g+1}$. On choisit un tel entier et on le note $\alpha$. On note $\mathcal{P}_{A}$ l'ensemble des idéaux premiers $\mathfrak{p}$ de $\mathcal{P}_{A,0}$ en lesquels la variété $\mathcal{A}_{\mathfrak{p}}$ a un $p$-rang égal à $\alpha$. De plus, quitte à diviser la densité de cet ensemble par $[K:\mathbb{Q}]$, on peut supposer que deux idéaux premiers distincts de $\mathcal{P}_{A}$ ont des normes, donc des projections sur $\mathbb{Z}$, distinctes. Dans les cas non ordinaires, on obtient une inégalité métrique moins bonne. Cette perte sera compensée, en $p$-rang égal à $0$,  par le plus grand nombre de points de torsion se réduisant sur $0$.  

\bigskip

On peut maintenant traduire le corollaire précédent en propriété $p$-adique. Soit $\mathfrak{p} \in \mathcal{P}_A$ un idéal premier et $p = \mathfrak{p} \cap \mathbb{Z}$. Par choix de $\mathcal{P}_A$, la fibre spéciale $\mathcal{A}_{\mathfrak{p}}$ est lisse. Soit $\Psi_{\mathfrak{p}}$ la différentielle du Verschiebung sur $\mathcal{A}_{\mathfrak{p}}$. Par la discussion du paragraphe précédent, on peut trouver une base de paramètres algébriques en l'origine de $A$: 
\begin{displaymath}
t_{\mathfrak{p},1}, \ldots, t_{\mathfrak{p},g}
\end{displaymath}
(\ie \ dont la projection est une base de $\mathfrak{m}_0/\mathfrak{m}_0^2$, où $\mathfrak{m}_0$ est l'idéal maximal correspondant à l'origine de $A$)
telle que son image par réduction modulo $\mathfrak{p}$: 
\begin{displaymath}
\tilde{t}_{\mathfrak{p},1}, \ldots, \tilde{t}_{\mathfrak{p},g}
\end{displaymath}
soit encore une base de paramètres algébriques, avec: 
\begin{eqnarray*}
\mathrm{Im} \Psi_{\mathfrak{p}}^g & = & \mathrm{Vect} (\tilde{t}_{\mathfrak{p},1}, \ldots, \tilde{t}_{\mathfrak{p},\alpha}), \\
\mathrm{Ker}\Psi_{\mathfrak{p}}^g & = & \mathrm{Vect} (\tilde{t}_{\mathfrak{p},\alpha+1}, \ldots, \tilde{t}_{\mathfrak{p},g}). \\
\end{eqnarray*}
On note $\mathcal{O}_{\mathfrak{p}}$ l'anneau de valuation associé à $\mathfrak{p}$. Il lui correspond par tensorisation: 
\begin{displaymath}
\mathcal{A}_{\mathcal{O}_{\mathfrak{p}}}:= \mathcal{A} \times_{\mathrm{Spec} \mathcal{O}_K} \mathrm{Spec} \mathcal{O}_{\mathfrak{p}},
\end{displaymath}
et la section nulle $\epsilon_{\mathfrak{p}}$. On note $\widehat{\mathcal{A}}_{\mathcal{O}_{\mathfrak{p}}}$ le complété le long de $\epsilon_{\mathfrak{p}}$ de $\mathcal{A}_{\mathcal{O}_{\mathfrak{p}}}$.
La multiplication par $[p]$ est donnée sur le groupe formel par un $g$-uplet de séries formelles: $\textbf{F}=(F_1, \ldots, F_g)$ et par réduction modulo $\mathfrak{p}$, on obtient un $g$-uplet de séries formelles: $\tilde{\textbf{F}}=(\tilde{F}_1, \ldots, \tilde{F}_g)$. La décomposition: $[p]_{\mathfrak{p}}= V_{\mathfrak{p}} \circ \phi_{\mathfrak{p}}$ de la multiplication par $p$ sur la fibre spéciale induit une décomposition:
\begin{displaymath}
\tilde{\textbf{F}}(\tilde{\textbf{t}}_{\mathfrak{p}})=\textbf{V}_{\mathfrak{p}} (\tilde{\textbf{t}}_{\mathfrak{p}}^p),
\end{displaymath}
où $\tilde{\textbf{t}}_{\mathfrak{p}}^p= (\tilde{t}_{\mathfrak{p},1}^p, \ldots, \tilde{t}_{\mathfrak{p}, g}^p)$ est l'image de $\tilde{\textbf{t}}_{\mathfrak{p}}$ par le Frobenius.

\bigskip

On a alors la proposition suivante:
\begin{prop}
\label{proposition metrique}
Si $P$ est un point de $p$-torsion se réduisant sur $0$ modulo $\mathfrak{q}$, pour une place $\mathfrak{q}/ \mathfrak{p}$ dans un corps de définition de $P$, on a: 
\begin{displaymath}
\forall \ 1 \leq i \leq \alpha:  |t_{\mathfrak{p}, i}(P)|_{\mathfrak{q}} \leq p^{-n_{\mathfrak{q}}/p}, 
\end{displaymath}
et: 
\begin{displaymath}
\forall \ \alpha < i \leq g:  |t_{\mathfrak{p}, i}(P)|_{\mathfrak{q}} \leq p^{-n_{\mathfrak{q}}/p^{2}}. 
\end{displaymath}
\end{prop}
\begin{dem}
Soit $P$ se réduisant sur $0$ modulo $\mathfrak{q}$, pour $\mathfrak{q}/\mathfrak{p}$ dans un corps de définition de $P$. On sait déjà que: 
\begin{displaymath}
\forall \ 1 \leq i \leq g:  |t_{\mathfrak{p}, i}(P)|_{\mathfrak{q}} < 1.
\end{displaymath} 
De plus, le morphisme $[p]$ en $P$ est donné par le $g$-uplet de séries formelles $\textbf{F}$ appliquées au système de paramètres (voir \cite{Hindry-Silverman00}, page 272). On en déduit: 
\begin{displaymath}
\textbf{F} \circ \textbf{t}_{\mathfrak{p}} (P)= 0.
\end{displaymath}
 
\bigskip

D'après le corollaire \ref{factorisation}, on a une factorisation: $\textbf{V}_{\mathfrak{p}}= \textbf{U}_{\mathfrak{p}} \circ \Phi_{\alpha, \mathfrak{p}}$, et la différentielle de $\textbf{U}_{\mathfrak{p}}$ est inversible. En utilisant la réduction modulo $\mathfrak{p}$ de $\textbf{F}$ et les propriétés de base d'une loi de groupe formel (rappelées dans \cite{Hindry-Silverman00}, page 269), on voit que $\textbf{F}$ est donnée par: 
\begin{displaymath}
\textbf{F}(\textbf{t}_{\mathfrak{p}})= p \textbf{t}_{\mathfrak{p}}+ \textbf{G}(\textbf{t}_{\mathfrak{p}}) + \textbf{H} \circ \Phi_{\alpha} (\textbf{t}_{\mathfrak{p}}^p). 
\end{displaymath} 
Le $g$ uplet de séries formelles $\textbf{G}$ a ses coefficients dans $\mathfrak{p} \mathcal{O}_{\mathfrak{p}}$ et ses premiers termes sont quadratiques; le $g$-uplet $\textbf{H}$ a ses coefficients inversibles modulo $\mathfrak{p}$ et sa différentielle est inversible dans $\mathcal{O}_{\mathfrak{p}}$. Soit $i_0 \in [1, g]$ tel que $|t_{\mathfrak{p}, i_0}(P)|_{\mathfrak{q}}$ soit maximal. En inversant la différentielle,  et par choix de $i_0$, on obtient:
\begin{displaymath}
 t_{\mathfrak{p}, i_0}(P)^{p^{n_{i_0}}} \in \mathfrak{p} \mathcal{O}_{\mathfrak{p}},
\end{displaymath}
où $n_{i_{0}}=1$ si $i_0 \leq \alpha$ et $n_{i_{0}}=2$ sinon. Comme l'indice de ramification $e_{\mathfrak{p}/p}$ vaut $1$ et par définition de $i_0$, on en déduit que pour tout $i \leq g$:
\begin{displaymath}
 |t_{\mathfrak{p}, i}(P)|_{\mathfrak{q}} \leq p^{-n_{\mathfrak{q}}/p²}. 
\end{displaymath}
Ceci étant démontré, on est assuré que les termes non-linéaires dans $\textbf{H}$ donnent une norme $\mathfrak{p}$-adique plus petite que $1$, et n'interfèrent pas avec le terme linéaire. On obtient donc, cette fois, pour tout $i \leq \alpha$:
\begin{displaymath}
t_{\mathfrak{p}, i}(P)^{p} \in \mathfrak{p} \mathcal{O}_{\mathfrak{p}},
\end{displaymath}
et la proposition est entièrement démontrée. 
\end{dem}
\textbf{Remarque} Puisque la ramification initiale: $e_{\mathfrak{p}/p}=1$ et le sous-groupe des points de $p$-torsion se réduisant sur $0$ modulo $\mathfrak{q}$ est galoisien de cardinal inférieur à $p^{2g-\alpha}$, le théorème de Raynaud (corollaire 3.4.4 de \cite{Raynaud74}) donne: 
\begin{displaymath}
|t_{\mathfrak{p}, i}(P)|_{\mathfrak{q}} \leq p^{-n_{\mathfrak{q}}/p^{(2g-\alpha)}},
\end{displaymath}
pour tout $i \leq g$. La théorie galoisienne ne suffit donc pas ici à démontrer la propriété métrique, alors qu'elle donne les mêmes résultats que l'approche des groupes formels dans le cas multi-elliptique. 

\bigskip

\label{coefficient rang}
On pose maintenant: $n_{\alpha}=1$ si $\alpha=g$ et $n_{\alpha}=2$ sinon. On choisit un système de paramètres en $0$: $(t_1, \ldots, t_g)$ associé à la base algébrique de $t_A$ (le tangent de $A$ en l'origine) définie en \ref{densite}. Comme par choix des premiers, cette base de paramètres est encore une base modulo $\mathfrak{p}$, on a immédiatement le corollaire:
\begin{cor}
\label{corollaire metrique}
Si $P$ est un point de $p$-torsion se réduisant sur $0$ modulo $\mathfrak{q}$, pour une place $\mathfrak{q}/ \mathfrak{p}$ dans un corps de définition de $P$, on a: 
\begin{displaymath}
\forall 1 \leq i \leq g:  |t_i(P)|_{\mathfrak{q}} \leq p^{-n_{\mathfrak{q}}/p^{n_{\alpha}}}. 
\end{displaymath}
\end{cor} 

\subsection{Base adaptée pour le tangent}

On montre dans ce paragraphe comment améliorer le corollaire précédent et trouver une base entière du tangent dans laquelle les propriétés $p$-adiques de la proposition \ref{proposition metrique} seront simultanément vraies, pour un nombre fini d'idéaux premiers, de norme bornée par un entier $N$. La hauteur des éléments de la base sera bornée par une fonction explicite de $N$.

\bigskip 

On rappelle que $t_A$ désigne l'espace tangent de $A$ en l'origine et on le munit cette fois de sa base canonique $(f_1, \ldots, f_g)$ pour le produit scalaire hermitien induit par la forme de Riemann. A tout $\mathfrak{p} \in \mathcal{P}_A$, on peut associer une base orthonormale $\textbf{e}_{\mathfrak{p}}=(e_{\mathfrak{p},1}, \ldots, e_{\mathfrak{p},g})$ de $t_A$ dans laquelle on dispose de bonnes propriétés métriques (par la proposition \ref{proposition metrique}). On veut trouver une base dans laquelle toutes les propriétés $\mathfrak{p}$-adiques sont lisibles simultanément, pour $\mathfrak{p} \in \mathcal{P}_A$. Ceci est possible si on borne dès maintenant la norme des premiers avec lesquels on travaille. Soit donc $N$ un entier; on suppose $N \geq g^2$ et on note: 
\begin{displaymath}
\mathcal{P}_{A,N}= \{ \mathfrak{p} \in \mathcal{P}_A, \textbf{N}(\mathfrak{p}) \leq N \}.
\end{displaymath}
On définit enfin la hauteur d'un vecteur $x \in K^g$ comme le maximum des hauteurs de ses coordonnées. 

\bigskip

On commence par démontrer un lemme de géométrie des nombres, qui fournit des éléments de petite hauteur dans une classe modulo $\mathfrak{b}$, pour $\mathfrak{b}$ un idéal de $\mathcal{O}_K$.

\begin{lem}
Soit $\mathfrak{b}$ un idéal de $\mathcal{O}_K$ de norme $\beta \in \mathbb{N}$. Il existe une constante $c_1$ ne dépendant que de $K$ telle que: pour toute classe $c$ modulo $\mathfrak{b}$, il existe un élément $x$ de $\mathcal{O}_K$ dans $c$ de hauteur plus petite que $\mathrm{log} \beta + c_1$.
\end{lem}
\begin{dem}
Au cours de la preuve, on notera par commodité: $n=[K:\mathbb{Q}]$.
On commence par plonger canoniquement $K$ dans $\mathbb{R}^n \simeq \mathbb{R}^{r_1} \times \mathbb{C}^{r_2}$ par un morphisme $\sigma$; on note $r_1$ le nombre de plongements réels et $r_2$ le nombre de plongements complexes à conjugaison près. Le module $\sigma(\mathfrak{b})$ est alors un sous-réseau de $\sigma(\mathcal{O}_K)$ et son covolume est donné par la formule (\cf \cite{Samuel}, 4.2): 
\begin{displaymath}
\mathrm{det} (\sigma(\mathfrak{b}))= 2^{-r_2}|d_K|^{1/2}\beta,
\end{displaymath}
où $d_K$ est le discriminant absolu de $K$. On applique le second théorème de Minkowski au réseau $\sigma(\mathfrak{b})$ et à l'ensemble: 
\begin{displaymath}
B= \{b=(y_1, \ldots, y_{r_1}, z_1, \ldots, z_{r_2}) \in \mathbb{R}^{r_1} \times \mathbb{C}^{r_2}, \left\| b \right\| = \sum_{i=1}^{r_1} |y_i|+2\sum_{j=1}^{r_2} |z_j| \leq 1 \}, 
\end{displaymath}
qui est compact, convexe et symétrique par rapport à l'origine. Son volume est donné par la formule (\cite{Samuel}, 4.2): 
\begin{displaymath}
\mathrm{vol} (B) = \frac{2^{r_1-r_2}\pi^{r_2}}{n!}.
\end{displaymath}
En notant, pour $1 \leq i \leq n$: 
\begin{displaymath}
\lambda_{i} =\mathrm{inf}\{ \lambda,\ \exists \ i \mathrm{\ vecteurs \ de \ } 
\sigma(\mathfrak{b}) \mathrm{\ libres \ dans\ } \lambda B \}, 
\end{displaymath}
on a l'inégalité des minima successifs:
\begin{displaymath}
\lambda_{1} \cdots \lambda_{n} \mathrm{vol} (B) 
\leq 2^{n}\mathrm{det} \sigma(\mathfrak{b}).
\end{displaymath}
On souhaite majorer le dernier de ces minima. Pour y parvenir, on remarque que la norme d'un élément $x$ de $\mathcal{O}_K$ est $\geq 1$ (c'est la valeur absolue du coefficient constant de son polynôme minimal sur $\mathbb{Z}$), et on en déduit, par l'inégalité arithmético-géométrique:
\begin{displaymath}
\big( \sum_{i=1}^n |\sigma_i(x)| \big)^n \geq n^n \textbf{N}(x) \geq n^n.
\end{displaymath}
Ceci permet de minorer le premier minimum: $\lambda_1 \geq n$, et par suite:
\begin{displaymath}
n^{n-1} \lambda_n \leq \lambda_1^{n-1} \lambda_n \leq \prod_{i=1}^n \lambda_i \leq \frac{2^{n}\mathrm{det} \sigma(\mathfrak{b})}{\mathrm{vol} (B)}.
\end{displaymath}
Soit $c$ une classe de $\mathcal{O}_K / \mathfrak{b}$; son image par $\sigma$ est une classe du quotient $\sigma(\mathcal{O}_K) / \sigma(\mathfrak{b})$. Comme la norme $\left\| . \right\|$ vérifie l'inégalité triangulaire, on peut choisir un représentant $x$ de $c$ tel que:
\begin{displaymath} 
\left\| x \right\| \leq n \lambda_n.
\end{displaymath} 
De plus, pour un point entier, seules les contributions archimédiennes interviennent dans la hauteur. Notons $I$ l'ensemble des indices pour lesquels $|\sigma_i(x)| \geq 1$. Par concavité de la fonction log: 
\begin{displaymath}
h(x)=\frac{1}{n} \sum_{i \in I} n_{\sigma_i} \mathrm{log} |x|_{\sigma_i} \leq \mathrm{log} \Big( \frac{1}{|I|} \sum_{i \in I}  n_{\sigma_i} |x|_{\sigma_i} \Big) \leq \mathrm{log} ( \left\| x \right\|) \leq \mathrm{log} (n \lambda_n). 
\end{displaymath}
On peut enfin conclure, à l'aide de l'expression du dernier minimum, et des formules donnant le déterminant du sous-réseau et le volume de $B$:
\begin{displaymath}
h(x) \leq \mathrm{log} \beta + (n+1) \mathrm{log} 2 + \frac{1}{2} \mathrm{log} |d_K|.
\end{displaymath} 
\end{dem} 

\bigskip

Le lemme précédent, allié au lemme des restes chinois, permet alors de trouver une bonne base simultanée, de hauteur correctement contrôlée:
\begin{prop}
Il existe $(e_1, \ldots, e_g)$ une base orthogonale de $t_A(K)$, correspondant par dualité à un système de paramètres $(t_1, \ldots, t_g)$ vérifiant la propriété suivante: pour tout $\mathfrak{p} \in \mathcal{P}_{A,N}$ tel que $\mathfrak{p} \cap \mathbb{Z}=p \mathbb{Z}$, si $P$ est un point de $p$-torsion se réduisant sur $0$ modulo $\mathfrak{q}$, pour une place $\mathfrak{q} / \mathfrak{p}$, on a: 
\begin{displaymath}
\forall 1 \leq i \leq \alpha:  |t_{i}(P)|_{\mathfrak{q}} \leq p^{-n_{\mathfrak{q}}/p}, 
\end{displaymath}
et: 
\begin{displaymath}
\forall \alpha < i \leq g:  |t_{i}(P)|_{\mathfrak{q}} \leq p^{-n_{\mathfrak{q}}/p^2}. 
\end{displaymath}
De plus, la hauteur des $e_i$ est majorée par $c_3 N$, pour une constante $c_3$ ne dépendant que de $g$ et $K$.
\end{prop}
\begin{dem}
Soit $1 \leq i \leq g$. Si on note $e_{i, \mathfrak{p}}^j$ la $j$-ème composante de $e_{i, \mathfrak{p}}$ selon la base canonique, le nombre $e_{i, \mathfrak{p}}^j$ est dans $\mathcal{O}_K$. Par le lemme des restes chinois: 
\begin{displaymath}
\prod_{\mathfrak{p} \in \mathcal{P}_{A,N}} \mathcal{O}_K / \mathfrak{p} \mathcal{O}_K \simeq \mathcal{O}_K / \big( \prod_{\mathfrak{p} \in \mathcal{P}_{A,N}} \mathfrak{p} \big) \mathcal{O}_K.
\end{displaymath}
Il existe donc $d_i^j$, un élément de $\mathcal{O}_K$, tel que: 
\begin{displaymath}
d_i^j \equiv e_{\mathfrak{p}, i}^j \mathrm{\  mod \ } \mathfrak{p},
\end{displaymath}
pour tout $\mathfrak{p} \in \mathcal{P}_{A,N}$ et par le lemme précédent, on peut prendre $d_i^j$ vérifiant:
\begin{displaymath}
h(d_i^j) \leq  \mathrm{log} \textbf{N} \bigg(\prod_{\mathfrak{p} \in \mathcal{P}_{A,N}} \mathfrak{p} \mathcal{O}_K \bigg) + c_1, 
\end{displaymath}
pour une constante $c_1$ ne dépendant que de $K$. La norme étant multiplicative pour les idéaux (voir \cite{Samuel}, 3.5), et en utilisant l'estimation de Tchébychev (\cite{Tenenbaum}, th. 3 page 11), on a:
\begin{displaymath}
h(d_i^j) \leq  \sum_{\mathfrak{p} \in \mathcal{P}_{A, N}} \mathrm{log}(N(\mathfrak{p})) + c_1 \leq  \sum_{p \leq N, p \in \mathcal{P}_{\mathbb{Z}}} \mathrm{log} (p) + c_1 \leq c_2 N,
\end{displaymath}
pour une constante $c_2$ ne dépendant que de $K$, avec $\mathcal{P}_{\mathbb{Z}}$ l'ensemble des premiers rationnels. On note $d_i$ le vecteur dont les coordonnées sont les $d_i^j$ dans la base canonique. La propriété $\mathfrak{p}$-adique ne dépendant que de la classe mod $\mathfrak{p}$ des coordonnées, elle reste vraie pour la base $\textbf{d}= (d_1, \ldots, d_g)$ et pour tout $\mathfrak{p} \in \mathcal{P}_{A,N}$. 
Le procédé d'orthogonalisation de Schmidt appliqué à $\textbf{d}$ fournit une base $\textbf{e}$  pour laquelle la propriété métrique reste inchangée (par l'inégalité ultramétrique). Les $e_i$ sont donnés par la formule: 
\begin{displaymath}
e_i = d_i - \sum_{k=1}^{i-1} \frac{(e_k, d_i)}{(e_k, e_k)} e_k. 
\end{displaymath}
Puisque $N \geq g^2$, on obtient par récurrence: $h(e_i) \leq c_3 N$ avec $c_3=4^gc_2$. 
\end{dem}

\section{Théorie des pentes}
\label{pentes}

Dans cette partie, on commence par définir les objets qui apparaissent dans la théorie des pentes, puis on donne les inégalités de pentes dont on se servira par la suite. On finit par estimer les pentes de fibrés qui apparaîtront dans la suite de ce travail. Un des avantages de la méthode des pentes, en géométrie diophantienne, est de rendre plus facile le calcul des constantes. Si le calcul de la constante ne dépendant que de $A$ reste très délicat dans notre minoration du minimum essentiel, on a tenu, dans cette partie, à donner des estimations précises, sinon en le corps de définition et en le genre de $A$, du moins en la hauteur de Faltings de $A$. 

\bigskip

Enfin, comme dans la partie \ref{partie metrique}, on a donné des résultats un peu plus généraux que ceux dont on se servira par la suite, puisqu'on a défini la multiplicité avec un sous-module (de rang maximal) du module tangent donné par un modèle semi-abélien. Ces calculs pourraient éventuellement permettre, dans un travail ultérieur, de mettre en oeuvre une technique de {\it multiplicités penchées}, qu'on n'a pu appliquer ici.

\subsection{Définitions et inégalité de pentes}

On souhaite mettre en \oe uvre le formalisme des pentes, introduit par Bost dans \cite{Bost95}, et qui s'est développé dans la littérature diophantienne depuis une dizaine d'années. Pour des détails et des exemples d'applications de la théorie des pentes, on renvoie par exemple aux articles de Bost (\cite{Bost95}, \cite{Bost01}) ou à l'article très complet de Gaudron (\cite{Gaudron06}). Le but de cette partie est donc d'écrire une {\it inégalité de pentes}. Sous sa forme basique, celle-ci compare les {\it pentes} de deux $\mathcal{O}_K$-modules hermitiens s'il existe un morphisme $\phi$ injectif entre eux. On va donc définir le degré arithmétique d'un fibré vectoriel hermitien, puis sa pente, sa pente maximale, et la hauteur d'un morphisme de fibrés. 

\bigskip

On note provisoirement (dans ce paragraphe): $S=\mathrm{Spec} (\mathcal{O}_K)$, $S_0$ l'ensemble des points fermés de $S$ et $S_{\infty}$ l'ensemble des places archimédiennes de $\mathcal{O}_K$ (correspondant aux points complexes de $S$); on note enfin $M(K)= S_0 \cup S_{\infty}$ l'ensemble des places de $K$.
Un fibré vectoriel $\mathcal{E}$ sur $S$ est constant, ce qui mène à la définition suivante:

\begin{defn}
Un fibré vectoriel hermitien sur $S$ est un $\mathcal{O}_K$-module $\mathcal{E}$ muni d'une collection $\{ \left\| . \right\|_{v} \}_{v \in S_{\infty}}$ telle que pour tout $v \in S_{\infty}$, $\left\| . \right\|_{v}$ soit une norme hermitienne sur le $K_v$-espace vectoriel $\mathcal{E}_v = \mathcal{E} \otimes K_v$ et qu'on ait la compatibilité suivante:
\begin{center}
$\left\| x \right\|_{v}= \left\| \overline{x} \right\|_{\overline{v}}$  pour tous $v \in S_{\infty}$ et $x \in \mathcal{E}_{v}$,
\end{center}
où $\overline{v}$ désigne la conjuguée de $v$ (via le plongement complexe qui leur est associé).  
\end{defn} 
On note $\overline{\mathcal{E}}$ le fibré $(\mathcal{E},  \{ \left\| . \right\|_{v}$ \}). Si $v$ est une place finie de $K$, on lui associe un anneau de valuation $\mathcal{O}_v$, et étant donnée une base $(e_1, \ldots, e_r)$ de $\mathcal{E}_v= \mathcal{E} \otimes \mathcal{O}_v$ sur $\mathcal{O}_K$, on munit $\mathcal{E}_v$ de la norme $\left\| . \right\|_{v}$ définie comme le maximum des valeurs absolues $v$-adiques des coordonnées dans cette base (voir \ref{densite}).  

\bigskip

On peut maintenant définir le degré arakelovien {\it normalisé} d'un fibré hermitien. On commence par les modules de rang $1$ puis on passe au cas général grâce au déterminant:
\begin{defn}
Soit $\overline{\mathcal{E}}$ un $\mathcal{O}_K$-fibré hermitien de rang $1$ et $s$ un élément non nul de $\mathcal{E}$. On définit le degré arithmétique (ou arakélovien) normalisé de $\mathcal{E}$ de la manière suivante: 
\begin{displaymath}
\widehat{\mathrm{deg}} \ \overline{\mathcal{E}} = \frac{1}{[K: \mathbb{Q}]}\bigg( \mathrm{log} \# (\mathcal{E}/s\mathcal{O}_K) - \sum_{v \in S_{\infty}} \mathrm{log} \| s \|_v  \bigg).
\end{displaymath}
Si $\overline{\mathcal{E}}$ est un $\mathcal{O}_K$-fibré hermitien de rang $r$, on pose: 
\begin{displaymath}
\widehat{\mathrm{deg}} \ \overline{\mathcal{E}} = \widehat{\mathrm{deg}} \ \mathrm{det} \overline{\mathcal{E}},
\end{displaymath}
où les normes sur le déterminant sont celles obtenues par puissance tensorielle et quotient à partir de celles de $\mathcal{E}$.
\end{defn}
\textbf{Remarques} La formule du produit montre que cette définition ne dépend pas du choix de $s$. Si $K= \mathbb{Q}$, le degré d'Arakelov est l'opposé du logarithme du covolume de $\mathcal{E}$ vu comme réseau de $\mathcal{E} \otimes_{\mathbb{Z}} \mathbb{R}$ (voir \cite{Bost-Gillet-Soule94}, formule (2.1.13) pour le cas général).

\begin{defn}
Soit $\overline{\mathcal{E}}$ un fibré hermitien de rang non nul, on définit sa pente de la façon suivante:
\begin{displaymath}
\widehat{\mu}(\overline{\mathcal{E}})= \frac{\widehat{\mathrm{deg}}\overline{\mathcal{E}}}{\mathrm{rg}\mathcal{E}}.
\end{displaymath}
\end{defn}
Les pentes des sous-modules de $\mathcal{E}$ sont bornées (par l'inégalité d'Hadamard), ce qui justifie la définition:
\begin{defn}
La pente maximale de $\overline{\mathcal{E}}$ est définie par:
\begin{displaymath}
\widehat{\mu}_{\mathrm{max}}(\overline{\mathcal{E}})=\mathrm{max} \ \widehat{\mu}(\overline{\mathcal{F}}),
\end{displaymath}
où $\mathcal{F}$ décrit l'ensemble des sous-fibrés non-nuls de $\mathcal{E}$ munis des métriques déduites de celles de $\mathcal{E}$ par restriction.
\end{defn}

\bigskip

Soit $\phi$ un morphisme entre deux $\mathcal{O}_K$-fibrés $\mathcal{E}$ et $\mathcal{F}$. Si ces fibrés sont hermitiens, pour toute place $v \in M(K)$, on note $\left\| \phi \right\| _v$ la norme d'opérateur du morphisme: 
\begin{displaymath}
\phi: \mathcal{E}_v=\mathcal{E} \otimes K_v \rightarrow \mathcal{F}_v=\mathcal{F} \otimes K_v.  
\end{displaymath}
On a donc:
\begin{displaymath}
\left\| \phi \right\| _v= \mathrm{sup}_{x \in \mathcal{E}_v, x \neq 0} \frac{\left\| \phi(x) \right\| _v}{\left\| x \right\| _v}.
\end{displaymath}  
On peut alors définir la hauteur de $\phi$:
\begin{defn}
Si $\phi$ est un morphisme entre deux $\mathcal{O}_K$-fibrés hermitiens $\overline{\mathcal{E}}$ et $\overline{\mathcal{F}}$, on appelle hauteur de $\phi$: 
\begin{displaymath}
h(\phi)=\frac{1}{[K:\mathbb{Q}]} \sum_{v \in M(K)} \mathrm{log} \left\| \phi \right\| _v.
\end{displaymath}
\end{defn}

\bigskip

On est alors en mesure d'écrire une première inégalité de pentes:
\begin{lem}
Si le morphisme $\phi: \mathcal{E}_K \rightarrow \mathcal{F}_K$ est injectif:
\begin{displaymath}
\widehat{\mathrm{deg}} \ \overline{\mathcal{E}} \leq \mathrm{rg} (\mathcal{E}) \Big( \widehat{\mu}_{\mathrm{max}}(\overline{\mathcal{F}}) + h(\phi) \Big).
\end{displaymath}
\end{lem}
\begin{dem}
C'est une conséquence de l'inégalité d'Hadamard; on renvoie à \cite{Chen06}, page 40, pour plus de détails (où la même convention est faite sur les normes ultramétriques).
\end{dem}
On utilise généralement une version filtrée de cette inégalité. 
Soit $\phi: \mathcal{E}_K \rightarrow \mathcal{F}_K$ une application $K$-linéaire injective. 
On suppose qu'il existe une filtration d'espaces vectoriels:
\begin{displaymath}
\{ 0 \} = \mathcal{F}_{K,N+1} \subset \cdots \subset \mathcal{F}_{K,0}=\mathcal{F}_K.
\end{displaymath}
et que les quotients $\mathcal{G}_{K,i}=\mathcal{F}_{K,i}/\mathcal{F}_{K, i+1}$ sont les tensorisation avec $K$ de fibrés hermitiens $\overline{\mathcal{G}}_{i}$ sur $\mathrm{Spec} \mathcal{O}_K$. 
On définit une filtration sur $\mathcal{E}_K$ par image réciproque: 
\begin{displaymath}
\mathcal{E}_{K, i}= \phi_K^{-1}(\mathcal{F}_{K, i}),
\end{displaymath}
et on considère:
\begin{displaymath}
\phi_{i}:\mathcal{E}_{K, i} \rightarrow \mathcal{G}_{K, i},
\end{displaymath}
l'application linéaire naturellement induite sur le quotient. La version filtrée de l'inégalité des pentes est la suivante:
\begin{lem}
\label{lemme pentes}
Si $\phi$ est injective, on a l'inégalité:
\begin{displaymath}
\widehat{\mathrm{deg}} \ \overline{\mathcal{E}} \leq \sum_{i=0}^{N} \mathrm{dim}\big(\mathcal{E}_{K,i}/\mathcal{E}_{K,i+1}\big) \big( \widehat{\mu}_{\mathrm{max}}(\overline{\mathcal{G}}_{i}) + h(\phi_{i}) \big).
\end{displaymath}
\end{lem}
Cette version est particulièrement utile dans une preuve de transcendance: la filtration correspond alors aux différents ensembles et ordres d'annulation d'une fonction auxiliaire.

\subsection{Quelques fibrés hermitiens et leurs pentes}
 
On peut maintenant préciser les fibrés hermitiens auxquels on compte appliquer la méthode des pentes. Rappelons que $A$ est une variété abélienne munie d'un fibré ample et symétrique $L$, et définie sur un corps de nombres $\mathcal{O}_K$.

\subsubsection{Tangent et espace symétrique}
\label{tangent}
Quitte à prendre une extension finie de $K$ ne dépendant que de $A$, on suppose que $A$ admet réduction semi-abélienne sur $K$. Il existe donc un modèle semi-abélien $\pi: \mathcal{A} \rightarrow \mathrm{Spec} \mathcal{O}_{K}$; on note $\epsilon: \mathrm{Spec} \mathcal{O}_{K} \rightarrow \mathcal{A}$ la section nulle et on pose: 
\begin{displaymath}
t_{\mathcal{A}}:=\epsilon^{*}T_{\mathcal{A}/\mathrm{Spec} \mathcal{O}_{K}},
\end{displaymath} 
où $T_{\mathcal{A}/\mathrm{Spec} \mathcal{O}_{K}}$ est le fibré tangent de $\mathcal{A}$ sur $\mathcal{O}_K$.
On a une structure de $\mathcal{O}_K$-module de type fini sur $t_{\mathcal{A}}$, qui est un sous-module de l'espace vectoriel tangent de $A$ à l'origine. Si $\sigma$ est une place complexe de $K$, il existe un isomorphisme: $t_{\mathcal{A}} \otimes_{\sigma} \mathbb{C} \simeq t_{\mathcal{A}_{\sigma}}(\mathbb{C})$ et la forme de Riemann $\omega$ associée au fibré ample $L$ induit une forme hermitienne $\omega_{\sigma} \in \bigwedge^{1,1} t_{\mathcal{A}_{\sigma}}^{\vee}(\mathbb{C})$. Si $\omega_{\sigma}$ s'écrit sous la forme: 
\begin{displaymath}
\omega_{\sigma}=\frac{i}{2} \sum_{1 \leq h,l \leq g}a_{h,l}f_h^* \wedge \overline{f_l^*},
\end{displaymath}
où $(f_1^*, \ldots, f_g^*)$ est la base duale d'une base $(f_1, \ldots, f_g)$ de $t_{\mathcal{A}_{\sigma}}(\mathbb{C})$, pour $(z_1, \ldots, z_g) \in \mathbb{C}^g$, on pose:
\begin{displaymath}
\left\| \sum_{i=1}^g z_if_i \right\|_{\sigma}^2=\sum_{1 \leq h,l \leq g} a_{h,l}z_h \overline{z_l}.
\end{displaymath}
Ces métriques font de $t_{\mathcal{A}}$ un fibré vectoriel hermitien et cette structure se transporte par dualité à $t_{\mathcal{A}}^{\vee}$.

\bigskip

On sait calculer explicitement la pente de $\overline{t_{\mathcal{A}}^{\vee}}$ en fonction de la hauteur de Faltings de $A$ (\cf \cite{Gaudron06}, proposition 4.7). Pour définir celle-ci, on considère le module $\omega_{\mathcal{A}/ \mathcal{O}_K} \simeq \mathrm{det} \ t_{\mathcal{A}}^{\vee}$ muni de la norme associée à chaque place archimédienne $\sigma$ de $K$: 

\begin{displaymath}
\forall s \in \omega_{\mathcal{A}/ \mathcal{O}_K} \otimes_{\sigma} \mathbb{C}, \ \  \left\|s\right\|_{\sigma}^2:=\frac{i^{g^2}}{(2 \pi)^g}\int_{\mathcal{A}_{\sigma}(\mathbb{C})} s \wedge \overline{s}.
\end{displaymath}

\begin{defn}
La hauteur de Faltings de $A$ est le degré d'Arakelov normalisé du fibré $\overline{\omega_{\mathcal{A}/ \mathcal{O}_K}}$. Elle est notée $h_F(A)$ et ne dépend pas du choix de $K$, le corps de définition, et de $\mathcal{A}$, le modèle semi-abélien de $A$ sur $\mathcal{O}_K$.
\end{defn}
Cette fonction vérifie les axiomes d'une hauteur sur les classes d'isomorphismes de variétés abéliennes de dimension $g$ sur $\overline{\mathbb{Q}}$ (voir le théorème (2.1) de \cite{Bost95}). Majorer la pente maximale de $\overline{t_{\mathcal{A}}^{\vee}}$ en fonction de $h_F(A)$ revient donc à comparer la pente maximale de deux modules isomorphes mais non isométriques; on obtient (\cf \cite{Bost95}, (5.41)):
\begin{lem}
Il existe une constante $c_4>0$ ne dépendant que de $g$ telle que:
\begin{displaymath}
\hat{\mu}_{\mathrm{max}}(\overline{t_{\mathcal{A}}^{\vee}}) \leq c_4 \mathrm{max} \{ 1, h_F(A), \mathrm{log} h^0(A,L) \},
\end{displaymath}
où $h^0(A,L)$ désigne la dimension de l'espace des sections du fibré $L$ sur $A$.
\end{lem} 
On a construit une base du tangent (sur $\mathcal{O}_K$) dans laquelle les propriétés $\mathfrak{p}$-adiques de $A$ sont lisibles, pour un grand nombre de premiers $\mathfrak{p}$ de $\mathcal{O}_K$. On souhaiterait pouvoir majorer, en général, la pente maximale associée à des sous-modules du tangent de rang maximal. Soit $\mathcal{W}$ un sous-module de $t_{\mathcal{A}}$ de rang $g$ engendré par des vecteurs $(e_1, \ldots, e_g)$; on munit ce module des métriques hermitiennes (issues de la forme de Riemann) de $t_{\mathcal{A}}$ par restriction. On peut majorer la pente maximale de $\overline{\mathcal{W}^{\vee}}$ en fonction de la hauteur d'une base de $\mathcal{W}$:
\begin{lem}
Soit $c(\mathcal{W})>0$  tel qu'il existe une base de $\mathcal{W}$ formée d'éléments dont toutes les coordonnées sont dans $\mathcal{O}_K$ et de hauteur plus petite que $c(\mathcal{W})$.
La pente maximale de $\overline{\mathcal{W}^{\vee}}$ vérifie alors:
\begin{displaymath}
\hat{\mu}_{\mathrm{max}}(\overline{\mathcal{W}^{\vee}}) \leq c_7 \mathrm{max} \{ 1, h_F(A), \mathrm{deg}_L(A), c(\mathcal{W})\},
\end{displaymath}
pour une constante $c_7 >0$ ne dépendant que de $g$ et $K$. 
\end{lem}
\begin{dem}
La valeur absolue du degré normalisé de $\mathcal{W}$ est majorée, à partir de la formule (2.1.13) de \cite{Bost-Gillet-Soule94}, de l'inégalité d'Hadamard et en prenant une base de $\mathcal{O}_K$ d'éléments dont la hauteur est bornée par $c(\mathcal{W})$: 
\begin{displaymath}
 | \widehat{\mathrm{deg}} \ \overline{\mathcal{W}}| \leq c_5 c(\mathcal{W}),
\end{displaymath}
pour une constante $c_5$ ne dépendant que de $g$ et de $K$.
En effet, pour les vecteurs à coordonnées dans $\mathcal{O}_K$, il n'y a pas de contribution ultramétrique dans l'expression de la hauteur, et les contributions archimédiennes sont comparables à la norme $L^2$ sortant de l'inégalité d'Hadamard.
Un raffinement du lemme précédent (formule (41) de \cite{Gaudron06}) donne: 
\begin{displaymath}
\hat{\mu}_{\mathrm{max}}(\overline{\mathcal{W}^{\vee}}) \leq c_6 \mathrm{max} \{ 1, h_F(A), \mathrm{deg}_L(A), |\widehat{\mathrm{deg}} \ \overline{\mathcal{W}}| \},
\end{displaymath}
où $c_6$ ne dépend que de $g$. La majoration de la valeur absolue du degré de $\overline{\mathcal{W}}$ permet de conclure.
\end{dem}

\bigskip

Pour passer aux dérivées d'ordre supérieur, on doit comprendre comment la pente maximale se comporte avec les puissances symétriques. Soit $\overline{\mathcal{E}}$ un fibré hermitien sur $\mathcal{O}_K$. Pour tout $m$, la $m$-ème puissance symétrique $S^m \mathcal{E}$ est munie d'une structure hermitienne par produit tensoriel puis projection; on note $S^m \overline{\mathcal{E}}$ le fibré hermitien ainsi obtenu. Le lemme suivant est démontré dans l'appendice de \cite{Graftieaux01}:
\begin{lem}
Pour $\overline{\mathcal{E}}$ un fibré hermitien de rang $g$:
\begin{displaymath}
\hat{\mu}_{\mathrm{max}}\big(S^m(\overline{\mathcal{E}})\big) \leq m(\hat{\mu}_{\mathrm{max}}(\overline{\mathcal{E}})+2g \mathrm{log}g).
\end{displaymath}
\end{lem}
La combinaison des deux derniers lemmes donne la proposition suivante, qui nous servira à majorer la pente maximale des fibrés d'arrivée dans l'inégalité de pentes:
\begin{prop}
On a la majoration suivante:
\begin{displaymath}
\hat{\mu}_{\mathrm{max}}\big(S^m(\overline{\mathcal{W}^{\vee}})\big) \leq c_8 m \ \mathrm{max} \{ 1, h_F(A), \mathrm{deg}_L(A), c(\mathcal{W})\},
\end{displaymath}
pour une constante $c_8$ ne dépendant que de $g$ et $K$.
\end{prop}

\subsubsection{L'espace des sections d'un fibré ample sur la variété abélienne}
\label{moret bailly}

La variété abélienne $A$ est munie de $L$ un fibré ample et symétrique, et on note $H^0(A,L)$ l'espace des sections de degré $1$ sur ce fibré. Si $s \in H^0(A,L)$, on peut définir une métrique sur $L$ à l'aide de la fonction $\theta$ associée à $L$. Pour $x=\mathrm{exp}_A(z)$ et $z \in t_A$, on pose en effet:
\begin{displaymath}
\left\|s(x)\right\|:=e^{-\frac{\pi}{2}\left\|z\right\|^2} |\theta(z)|.
\end{displaymath}
Cette métrique est appelée cubiste et sa forme de courbure est invariante par translation. 

\bigskip

Il existe alors (voir \cite{Bost96}, §4.3) un modèle de $(A,L)$, appelé modèle de Moret-Bailly et noté $(\mathcal{A}, \overline{\mathcal{L}}, 0)$, constitué d'un schéma abélien: 
\begin{displaymath}
\pi: \mathcal{A} \rightarrow \mathrm{Spec}\mathcal{O}_K,
\end{displaymath}
et d'un fibré hermitien $\overline{\mathcal{L}}$ sur $\mathcal{A}$ vérifiant notamment les propriétés suivantes:

\bigskip

\begin{itemize}
\item  Il existe un isomorphisme de variétés abéliennes sur $\overline{\mathbb{Q}}$: $i: A \stackrel{~}{\rightarrow} \mathcal{A}_{\overline{\mathbb{Q}}}$.

\item Il existe un isomorphisme de fibrés en droites sur $A$: $i^*\mathcal{L}_{\overline{\mathbb{Q}}} \stackrel{~}{\rightarrow} L$.

\item L'origine de $A$ se relève en une section $\epsilon: \mathrm{Spec} \mathcal{O}_K \rightarrow \mathcal{A}$. 

\item Pour toute place $\sigma$ archimédienne de $K$, la métrique sur $\mathcal{L} \otimes_{\sigma} \mathbb{C}$ est la métrique cubiste définie plus haut. 
\end{itemize}

\bigskip

On note $\mathcal{E}$ le $\mathcal{O}_K$-module $H^0(\mathcal{A}, \mathcal{L})$ et on le munit des métriques hermitiennes suivantes aux places archimédiennes: pour tout plongement $\sigma: K \hookrightarrow \mathbb{C}$ et $s \in \mathcal{E} \otimes_{\sigma} \mathbb{C} \simeq H^0(\mathcal{A}_{\sigma}, \mathcal{L}_{\sigma})$, on pose:
\begin{displaymath}
\left\|s\right\|_{L^2, \sigma}= \bigg( \int_{\mathcal{A}_{\sigma}(\mathbb{C})} \left\|s(x)\right\|²d\mu(x) \bigg)^{1/2},
\end{displaymath}
où $d\mu$ est la mesure de Haar de masse totale égale à $1$ sur $\mathcal{A}_{\sigma}(\mathbb{C})$.
On sait que $\overline{\mathcal{E}}$ est {\it semi-stable}: sa pente est égale à sa pente maximale. De plus, on peut la calculer de façon totalement explicite:
\begin{prop} 
\label{theoreme moret bailly}
La pente de $\overline{\mathcal{E}}$ est donnée par la formule suivante:
\begin{displaymath}
\hat{\mu}(\overline{\mathcal{E}})=-\frac{1}{2}h_F(A)+\frac{1}{4}\mathrm{log}\frac{\chi(A,L)}{(2\pi)^g};
\end{displaymath}
où $\chi(A,L)=\frac{\mathrm{deg}_L(A)}{g!}$ désigne la caractéristique d'Euler-Poincaré de $L$.
\end{prop}
Cette formule a été démontrée par Moret-Bailly (\cite{Moret-Bailly90}) dans le cas $\chi(A, L)=1$ et en général par Bost (\cite{Bost96}). La semi-stabilité de $\overline{\mathcal{E}}$ est démontrée dans \cite{Bost95}, 4.2.

\subsubsection{Pente maximale et sous-variétés}
\label{pentesousvariete}

Le morphisme de fibrés qu'on définira dans la partie suivante associera à une section de $L^{\otimes M}$ sur $A$ (pour $M \geq 1$ un entier) sa restriction à un fermé de Zariski $X$, avec multiplicité (dans un plongement étiré). On aura donc besoin de majorer la pente maximale du fibré des sections sur des sous-variétés de $A$ avec multiplicité.

\bigskip

Soit $X$ un fermé de Zariski lisse et équidimensionnel de $A$, de dimension $d$, défini sur une extension $K'$ de $K$. On en prend un modèle, c'est-à-dire un $\mathcal{O}_{K'}$-schéma propre et plat $\mathcal{X}$ tel que $\mathcal{X} \times_{\mathrm{Spec}\mathcal{O}_{K'}} K' =X$ et on choisit un $\mathcal{O}_K-$fibré $\mathcal{L}$ sur $\mathcal{X}$ dont la restriction à $X$ est $L$. Soit de plus $\mathcal{W}$ un sous-module du tangent $t_{\mathcal{A}}$, muni de la restriction des métriques hermitiennes provenant de la forme de Riemann. Pour chaque place archimédienne $\sigma$ de $K'$, on munit le $\mathcal{O}_K$-fibré: 
\begin{displaymath}
\mathcal{L}^{\otimes M} \otimes S^m (\mathcal{W}^{\vee})
\end{displaymath} 
de la norme obtenue par produit tensoriel de la norme cubiste sur $\mathcal{L} \otimes_{\sigma} \mathbb{C}$ et de la norme symétrique sur $S^m (\mathcal{W}^{\vee}) \otimes_{\sigma} \mathbb{C}$; puis on munit l'espace de sections, {\it vu comme} $\mathcal{O}_K$-{\it module}: 
\begin{displaymath}
\mathcal{H}_{m}^M=H^0 \Big(\mathcal{X}, \mathcal{L}^{\otimes M} \otimes S^m (\mathcal{W}^{\vee})\Big)
\end{displaymath} 
de la métrique de Löwner $\left\|.\right\|_{L, \sigma}$ associée à la norme du sup, pour toute place archimédienne $\sigma$. Il s'agit d'une norme hermitienne proche de la norme du sup dans le sens suivant (\cf \cite{Gaudron07}, 2.2 pour le cas euclidien): 
\begin{displaymath}
\forall x \in \mathcal{H}_m^M \otimes \mathbb{C}: \left\|x\right\|_{L, \sigma} \leq  \left\|x\right\|_{\sup, \sigma} \leq \sqrt{2\mathrm{rg} \mathcal{H}^M_m} \left\|x\right\|_{L, \sigma}.
\end{displaymath}
\textbf{Remarque} En termes géométriques, l'existence de la norme de Löwner est reliée à l'existence, étant donné un corps $C$ convexe et symétrique, d'un ellipsoïde de volume minimal contenant $C$.

\bigskip

On étudie d'abord le cas sans multiplicité: $m=0$. S'il s'agit d'estimer la pente et non la pente maximale, on sait obtenir une formule asymptotique lorsque le degré tend vers l'infini. En effet, en combinant les théorèmes de Hilbert-Samuel arithmétique et géométrique, on obtient: 
\begin{thm}
En munissant $\mathcal{H}_{0}^M$ de la norme du sup, on a l'équivalent asymptotique suivant pour la pente:
\begin{displaymath}
\hat{\mu}(\overline{\mathcal{H}_{0}^M})= M\frac{h_{\overline{L}}(X)}{\mathrm{deg}_L(X)}+ o(M).
\end{displaymath}
\end{thm}
\textbf{Remarque} La définition de la hauteur $h_{\overline{L}}(V)$ d'une variété $V$ est donnée dans \cite{Bost-Gillet-Soule94}, partie 3; elle correspond à la définition de la hauteur projective donnée par Philippon en faisant le choix de la norme $L^2$ pour les places archimédiennes (voir \cite{Hindry-Silverman00}). \\
\\
\begin{dem}
Il suffit de combiner le théorème de Hilbert-Samuel arithmétique (voir \cite{Gillet-Soule92} ou le théorème (1.4) de \cite{Zhang95}) et le théorème de Hilbert-Samuel géométrique (\cf \cite{Hartshorne77}, page 51).  
\end{dem}

On désire obtenir une majoration efficace de la pente maximale par le premier terme de ce développement asymptotique et des termes d'erreur bien contrôlés. Dans sa thèse (\cite{Chen06}), H. Chen montre que le quotient 
\begin{displaymath}
\frac{\hat{\mu}_{\mathrm{max}}(\overline{\mathcal{H}_{0}^M})}{M}
\end{displaymath} 
converge aussi. Plus précisément, il donne une majoration assez fine de la pente maximale en munissant le fibré de la norme $L^2$. Cette inégalité fait apparaître une constante non-explicite provenant de l'inégalité de Gromov, et c'est pour contourner cette difficulté qu'on introduit la norme de Löwner. De plus, on précise le terme principal de la majoration à l'aide du ``théorème de Wirtinger'' et de l'inégalité des minima successifs de Zhang. 
\begin{prop}
\label{wirtinger}
Il existe une constante explicite $c_9$ ne dépendant que de $K$ telle qu'on ait la majoration suivante pour la pente maximale de $\overline{\mathcal{H}_{0}^M}$:
\begin{displaymath}
\hat{\mu}_{\mathrm{max}}(\overline{\mathcal{H}_{0}^M}) \leq M \frac{h_{\overline{L}}(X)}{\mathrm{deg}_{L}(X)}+ d \mathrm{log} M + \mathrm{log} \ \mathrm{deg}_L(X) + \  c_9.
\end{displaymath} 
\end{prop} 
\begin{dem}
On commence par majorer la pente maximale en introduisant la plus petite norme $\epsilon(\overline{\mathcal{H}_{0}^M})$ d'un élément non-nul du réseau $\mathcal{H}_{0}^M$ via le premier théorème de Minkowski (\cf \cite{Bost-Kunnemann06}, inégalité (3.24) et la majoration de la fonction $\psi$ page 35); cette norme est la norme hermitienne sur la somme orthogonale des $\mathcal{H}_{0}^M \otimes_{\sigma} \mathbb{C}$. On a:

\begin{displaymath}
\hat{\mu}_{\mathrm{max}}(\overline{\mathcal{H}_{0}^M}) \leq -\mathrm{log} \ \epsilon(\overline{\mathcal{H}_{0}^M}) + \frac{1}{2}\mathrm{log \ rg}_{\mathbb{Z}}(\mathcal{H}_{0}^M)+ \frac{\mathrm{log}|\Delta_K|}{2[K:\mathbb{Q}]}.
\end{displaymath}
On majore le terme $-\mathrm{log} \ \epsilon(\overline{\mathcal{H}_{0}^M})$ à l'aide de la théorie de l'intersection arithmétique. On peut appliquer la formule reliant la hauteur d'une variété à la hauteur d'un diviseur sur cette variété (\cf \cite{Bost-Gillet-Soule94}, proposition 3.2.1) et comme $L$ est ample avec $c_1(\overline{L})$ définie positive (par la proposition 3.2.4 de \cite{Bost-Gillet-Soule94}), pour $s$ un diviseur effectif de $\overline{L}^M$, on a:
\begin{displaymath}
h_{\overline{L}}(\mathrm{div}(s))=M h_{\overline{L}}(X)+ \int_{X(\mathbb{C})}\mathrm{log}\left\| s \right\|c_1(\overline{L})^{d} \geq 0.
\end{displaymath}
De plus, comme $c_1(\overline{L})$ est définie positive, on a: 
\begin{displaymath}
\int_{X(\mathbb{C})}\mathrm{log}\left\| s \right\|c_1(\overline{L})^{d} \leq \mathrm{max}_{\sigma \in \Sigma_{\infty}} \mathrm{log} \left\| s \right\|_{\mathrm{sup}, \sigma} \int_{X(\mathbb{C})}c_1(\overline{L})^{d}.
\end{displaymath}
Par le théorème de Wirtinger (voir \cite{Griffiths-Harris78}, page 171), cette dernière intégrale n'est autre que 
\begin{displaymath}
\mathrm{deg}_L(X).  
\end{displaymath}
On a donc: 
\begin{displaymath}
-\mathrm{max}_{\sigma \in \Sigma_{\infty}} \mathrm{log} \left\| s \right\|_{\mathrm{sup}, \sigma} \leq M  \frac{h_{\overline{L}}(X)}{\mathrm{deg}_{L}(X)}.
\end{displaymath}
Soit $\sigma$ tel que $\left\| s \right\|_{\mathrm{sup}, \sigma}$ réalise le maximum sur toutes les places archimédiennes. Par choix de la norme de Löwner associée à la norme du sup sur $\mathcal{H}_{0}^M$, on a: 
\begin{displaymath}
\mathrm{max}_{\sigma \in \Sigma_{\infty}} \left\| s \right\|_{\mathrm{sup}, \sigma} \leq \sqrt{2\mathrm{rg} \mathcal{H}^M_0} \left\| s \right\|_{L, \sigma}, 
\end{displaymath}
et par suite: 
\begin{displaymath}
-  \mathrm{log} \left\| s \right\| = -\frac{1}{2}\mathrm{log} \Big(\sum_{\sigma \in \Sigma_{\infty}} \left\| s \right\|_{L, \sigma}^2\Big)\leq - \mathrm{log} \left\| s \right\|_{L, \sigma}  \leq  - \mathrm{max}_{\sigma \in \Sigma_{\infty}} \mathrm{log} \left\| s \right\|_{\mathrm{sup}, \sigma} +  \mathrm{log} \ \mathrm{rg} \mathcal{H}^M_0.
\end{displaymath}
On en déduit que: 
\begin{displaymath}
-\mathrm{log} \ \epsilon(\overline{\mathcal{H}_{0}^M}) \leq  M \frac{h_{\overline{L}}(X)}{\mathrm{deg}_{L}(X)}  + \mathrm{log} \ \mathrm{rg} \mathcal{H}^M_0.
\end{displaymath}
Il reste à majorer le rang de $\mathcal{H}_{0}^M$, ce qui suit du théorème de Chardin (\cf \cite{Chardin89}): 
\begin{displaymath}
\mathrm{rg}_{\mathbb{Z}} \mathcal{H}_{0}^M \leq  [K: \mathbb{Q}] \mathrm{rg} \mathcal{H}_0^M \leq [K: \mathbb{Q}]M^d \mathrm{deg}_L(X). 
\end{displaymath}
La proposition est donc entièrement démontrée.
\end{dem}

\bigskip

On peut maintenant majorer la pente maximale du fibré des sections avec multiplicité. Le choix, classique en géométrie diophantienne, de prendre le tangent de la variété abélienne, qui est un fibré constant (et non le tangent à la sous-variété), simplifie le calcul. Il semble possible d'obtenir une majoration dans le cas du tangent à la sous-variété en calculant les classes de Segré géométriques (voir \cite{Fulton84}) et arithmétiques (voir \cite{Mourougane04}) associées au fibré tangent $t_{\mathcal{X}}$. Dans le cas qui nous intéresse, le fibré $\mathcal{W}^{\vee}$ étant constant, on a l'isomorphisme (isométrique):
\begin{displaymath}
\mathcal{H}_{m}^M \simeq \mathcal{H}_{0}^M \otimes S^m \mathcal{W}^{\vee}.  
\end{displaymath}

Bost a conjecturé que la pente maximale du produit tensoriel est la somme des pentes maximales. Les meilleurs résultats connus dans cette direction sont ceux de Bost-Künnemann (\cf \cite{Bost-Kunnemann06}), et de H. Chen, qui a démontré dans sa thèse (\cf \cite{Chen06}):
\begin{thm}
Soient $\overline{\mathcal{E}_1}, \ldots, \overline{\mathcal{E}_n}$ des fibrés hermitiens non-nuls sur $\mathrm{Spec}\mathcal{O}_K$, alors:
\begin{displaymath}
\hat{\mu}_{\mathrm{max}}(\overline{\mathcal{E}_1} \otimes \cdots \otimes \overline{\mathcal{E}_n}) \leq \sum_{i=1}^n \Big(\hat{\mu}_{\mathrm{max}}(\overline{\mathcal{E}_i})+ \mathrm{log} (\mathrm{rg}\mathcal{E}_i)\Big).
\end{displaymath} 
\end{thm}
\textbf{Remarque} Pour un produit de deux fibrés, le résultat de Bost et Künnemann fait apparaître un facteur $1/2$ pour le terme logarithmique (mais aussi, en contrepartie, une constante ne dépendant que de $K$). 

\bigskip

Rappelons que par $c(\mathcal{W})$, on désigne un majorant de la hauteur d'une base entière de $\mathcal{W}$. On peut maintenant démontrer le corollaire suivant:
\begin{cor}
\label{pente sous variete}
On a la majoration, pour la pente maximale:
\begin{eqnarray*}
\hat{\mu}_{\mathrm{max}}(\overline{\mathcal{H}_{m}^M})  \leq  M \frac{h_{\overline{L}}(X)}{\mathrm{deg}_{L}(X)} + 2d \mathrm{log} (M) + 2\mathrm{log} \ \mathrm{deg}_L(X)  + c_{10} m \ \mathrm{max} \{ 1, h_F(A), \mathrm{deg}_L(A), c(\mathcal{W}) \},
\end{eqnarray*} 
où $c_{10}$ ne dépend que de $K$ et de $g$.
\end{cor}
\begin{dem}
Le théorème précédent nous montre que:
\begin{displaymath}
\hat{\mu}_{\mathrm{max}} (\overline{\mathcal{H}_{m}^M}) \leq \hat{\mu}_{\mathrm{max}} (\overline{\mathcal{H}_{0}^M}) + \hat{\mu}_{\mathrm{max}} (S^m (\overline{\mathcal{W}^{\vee}})) + \mathrm{log}(\mathrm{rg}\mathcal{H}_{0}^M)+\mathrm{log}(\mathrm{rg} S^m (\mathcal{W}^{\vee})).
\end{displaymath} 
Les termes relatifs à $\overline{\mathcal{H}_{0}^M}$ ont été calculés dans la proposition précédente. Le rang de $S^m (\mathcal{W}^{\vee})$ est donné par la formule classique:
\begin{displaymath}
\mathrm{rg} \ S^m (\mathcal{W}^{\vee})= \Big( ^{m+g-1}_{\ \ g-1}\Big) \leq m^{g-1},
\end{displaymath}
et sa pente maximale a été majorée en (3.2.1).
\end{dem}

\section{Choix des fibrés et du morphisme}
\label{choix fibres}

Après ces préliminaires, la preuve du théorème \ref{theoreme} commence véritablement ici. Soit $A$ une variété abélienne définie sur un corps de nombres $K$, munie d'un fibré $L$ ample et symétrique, qu'on pourra supposer très ample quitte à considérer $L^{\otimes 3}$. Par la suite, lorsque cela ne sera pas précisé, le degré et la hauteur seront définis par rapport à ce fibré. Quitte à prendre une extension finie de $K$ (ne dépendant que de $A$), on prend un modèle de Moret-Bailly $(\mathcal{A}, \mathcal{L}, 0)$ de $(A,L)$, suivant la terminologie de Bost (\cite{Bost96}). Ce modèle est en particulier semi-abélien (voir \cite{Gaudron06}, définition-théorème 4.3). On suppose de plus que $\mathcal{A}$ vérifie \textbf{H3}.  

\bigskip

On prend une sous-variété propre (et irréductible) $V$ de $A$, de codimension $r$, qui n'est pas incluse dans le translaté d'une sous-variété abélienne. On pose aussi: $X= V + \Sigma_0$, en vue de la dernière phase de la preuve. L'ensemble $\Sigma_0$ sera un sous-groupe fini de $A$. On va construire un fibré hermitien $\overline{\mathcal{E}_M}$ associé à un espace vectoriel $E$, une suite de fibrés $\overline{\mathcal{G}_k}$, $k \in I$ (pour un certain ensemble fini $I$), correspondant à une fibration d'un espace vectoriel $F$, et un morphisme $\phi$ de restriction entre $E$ et $F$. Les $\overline{\mathcal{G}_k}$ seront définis à partir de l'espace des sections d'une puissance de $\mathcal{L}$ sur des modèles entiers de $X$ et de ses translatés par des points de torsion bien choisis. La partie précédente nous permettra de calculer les termes de pentes associés à ces fibrés. A la fin de cette partie, on fixera les paramètres et on supposera par l'absurde que le minimum essentiel de $V$ est majoré en fonction des paramètres. Ce n'est qu'après avoir obtenu une contradiction qu'on en déduira le théorème \ref{theoreme}, en calculant explicitement les paramètres.

\subsection{Le plongement \'etir\'e}

Pour \'eliminer la constante de comparaison entre hauteur projective et hauteur
de Néron-Tate sur $A$, on considère classiquement un plongement \'etir\'e.

Soit $M$ un entier supérieur ou égal à $1$. La multiplication par $M$ sur $A$ est notée $[M]$ et on d\'efinit $\phi_M$:
\begin{eqnarray*}
  A & \rightarrow & \ \ A \times A\\
  x & \rightarrow & (x, [M] x).
\end{eqnarray*}
Ce plongement a été utilisé pour la première fois par Laurent pour étudier le
problème de Lehmer elliptique (\cf \cite{Laurent83}). Son principe est le suivant: les techniques diophantiennes nous renseignent sur la hauteur projective, et le minimum essentiel fait intervenir la hauteur de Néron-Tate associée au plongement. On sait que la différence entre ces deux hauteurs est bornée mais la hauteur de Néron-Tate peut être très petite; il y a donc une perte d'information sur la hauteur projective. Le plongement étiré multiplie la hauteur par un paramètre assez grand, qui rend négligeable la constante de comparaison.

\bigskip

On notera $L_M$ (resp. $\mathcal{L}_M$) l'image par $\phi_M$ de $L$ (resp. $\mathcal{L}$). On a: 
\begin{displaymath}
L_M \simeq L^{\otimes M^2+1}; 
\end{displaymath}
par abus de langage, on utilise la même notation pour le fibré induit sur $A$ par le plongement. On note deg$_M$ (resp. $\hat{h}_M$) le degré (resp. la hauteur canonique) par rapport au fibré $L_M$. Le lemme suivant indique la variation de la hauteur et du degr\'e par changement de fibré:
\begin{lem}
\label{second etirement}
Si $V$ est une sous-variété de $A$, on a la variation suivante de la hauteur canonique:
\begin{displaymath}
\hat{h}_M(V)= (M^2+1)^{\mathrm{dim}(V)+1}\hat{h}(V);
\end{displaymath}
et la formule suivante pour le degré:
\begin{displaymath}
\mathrm{deg}_M (V)=(M^2+1)^{\mathrm{dim}(V)}\mathrm{deg} (V).
\end{displaymath}
\end{lem}
\begin{dem}
Ces formules sont démontrées, par exemple, dans \cite{Philippon95}, proposition 7.
\end{dem}

Dans cette partie et la suivante, on travaillera donc dans le plongement étiré. On confondra la variété abélienne $A$ et ses sous-variétés avec leurs images par $\phi_M$, mais on précisera toujours avec soin le fibré.   

\subsection{Le fibré de départ}

Commençons par décrire le premier fibré. On pose $\mathcal{E}_M=H^0(\mathcal{A}, \mathcal{L}_M)$, le module des sections sur $\mathcal{A}$ sur le fibré $\mathcal{L}_M$. On note aussi $E$ le tensorisé de $\mathcal{E}_M$ avec $K$. On a vu en \ref{moret bailly} comment munir ce fibré d'une structure hermitienne et on a donné son degré dans le cas des sections de degré $1$; le cas général s'en déduit en remplaçant $\mathcal{L}$ par $\mathcal{L}^{\otimes (M^2+1)}$. On fait une hypothèse qui sera aisément vérifiée par la suite, en supposant que:
\begin{displaymath}
\mathrm{log} (M^2) > 2 \Big(\frac{2}{g}h_F(A) + \mathrm{log} (2 \pi g!)\Big). 
\end{displaymath}

\begin{lem}
\label{degre e}
On a la minoration suivante pour le degré normalisé de $\overline{\mathcal{E}_M}$:
\begin{displaymath}
\widehat{\mathrm{deg}} \ \overline{\mathcal{E}_M} \geq c_{11} (M^2)^g \mathrm{log} (M),
\end{displaymath}
pour une constante $c_{11}>0$ ne dépendant que de $A$.
\end{lem}
\begin{dem}
En multipliant le degré arithmétique par le rang, la proposition \ref{theoreme moret bailly} donne:
\begin{displaymath}
\hat{\mu}(\overline{\mathcal{E}_M})=(M^2+1)^g \frac{\mathrm{deg}_L(A)}{g!}\bigg(-\frac{1}{2}h_F(A)+\frac{1}{4}\mathrm{log}\frac{(M^2+1)^g \mathrm{deg}_L(A)}{(2\pi g!)^g}\bigg).
\end{displaymath}
Notons que la hauteur de Faltings ne dépend que de la classe d'isomorphisme de la variété abélienne, et qu'elle est donc invariante par le plongement étiré. L'hypothèse faite sur les paramètres donne immédiatement le lemme, avec une constante $c_{11}>0$ ne dépendant que de $A$ et facilement explicitable.
\end{dem}

\subsection{Voisinages infinitésimaux}
\label{derivations}

On définit maintenant le fibré d'arrivée pour appliquer l'inégalité des pentes. Dans la littérature existante, ce fibré est très souvent formé à partir d'un nombre fini de points. Mais ici, on veut que les sections s'annulent sur le fermé $X$ avec multiplicité, puis sur ses translatés. Pour former le fibré d'arrivée, on commence par introduire la notion de voisinage infinitésimal (\cf \cite{Bost96}). Ce voisinage sera donné par un fermé $Y$ de $A$ et un ordre de dérivation $l$, par rapport au tangent $t_{A_M}$, image de $t_{A}$ par l'étirement (voir \ref{tangent} pour des détails sur le fibré tangent).

\bigskip

On fixe désormais une base $(f_1, \ldots, f_g)$ de l'espace tangent $t_{A_M}$. On commence par choisir la base du tangent $(e_1, \ldots, e_g)$ sur $A$ correspondant à la base de dérivations algébriques fixée en \ref{densite}. L'action de l'isogénie $[M]$ sur le tangent étant la multiplication par $M$, on en déduit d'abord une base $(e_1, \ldots, e_{2g})$ de $t_{A^2}$, puis on pose: $f_i=e_i+Me_{g+i}$. Il correspond à $(f_1, \ldots, f_g)$ une base de dérivations qu'on note encore (par abus de langage) $(\partial_1, \ldots, \partial_g)$. Si $x \in A$, on en déduit par translation une base du tangent $t_{A_M,x}$ en $x$ (dans le plongement étiré), associée à la base de dérivations $(\partial_{1,x}, \ldots, \partial_{g,x})$. On notera aussi par la suite: 
\begin{displaymath}
\partial_{x}^{\lambda}= \frac{1}{\prod_{1=i}^{g} \lambda_i!}\partial_{1,x}^{\lambda_1} \cdots \partial_{g,x}^{\lambda_g}.
\end{displaymath}
Soit $Y$ un fermé de Zariski de $A$ (dans le plongement étiré). On définit le schéma $V(Y, t_{A_M},l)$ de la façon suivante:

\bigskip

\begin{itemize}
\item si $l=0$, $V(Y, t_{A_M},l)$ est le sous-schéma réduit de $A$ défini par $Y$.
\item si $l=1$ et $Y=0$, $V(0,t_{A_M},1)$ est le voisinage infinitésimal d'ordre $1$ de $0$.
\item si $l \geq 1$, $V(Y,t_{A_M},l)$ est l'image dans $A$ du schéma $Y \times V(0, t_{A_M}, 1)^l$ par le morphisme d'addition $A^{l+1} \rightarrow A$.
\end{itemize}

\bigskip

Le schéma $V(Y, t_{A_M} ,l)$ admet pour support le fermé $Y$ et son faisceau d'idéaux $\mathcal{I}$ est défini par:
\begin{displaymath}
s \in \mathcal{I} \Leftrightarrow \Big( \forall y \in Y, \ \  \forall \lambda \in \mathbb{N}^{g} \ \mathrm{tq:} \sum_{i=1}^{g} \lambda_i \leq l: \  \partial_{y}^{\lambda} s =0 \Big).
\end{displaymath}

\subsection{Fibration de l'espace d'arrivée}

L'espace vectoriel d'arrivée sera formé à partir des espaces de sections sur des modèles de $X$ et de nombreux translatés de $X$ par des points de torsion, avec multiplicité. On doit donc d'abord préciser les points de torsion en lesquels on extrapole, et la multiplicité, puis mettre un ordre sur cet ensemble. 

\bigskip

On se donne $T_0>0$. Pour $1 \leq n \leq r$ ($r$ étant la codimension de $V$ et de $X$), on se donne aussi deux nombres positifs $T_n$ (qui correspond à la multiplicité après $n$ extrapolations), et $N_n$ (qui borne les normes des premiers d'extrapolation). 
Puis on définit les ensembles $\mathcal{P}_n$: 
\begin{displaymath}
\mathcal{P}_n= \{\mathfrak{p} \in \mathcal{P}_{A} \ \mathrm{avec} \ \textbf{N}(\mathfrak{p}) \in  [N_n/2; N_n] \}.
\end{displaymath}
On note $\mathcal{P}$ l'ensemble des premiers de $\mathbb{Z}$. La projection des $\mathcal{P}_n$ sur $\mathbb{Z}$ est donnée par:
\begin{displaymath}
\mathcal{P}_{n, \mathbb{Z}}= \{p \in \mathcal{P}, \exists \ \mathfrak{p} \in \mathcal{P}_n, \mathfrak{p}/p \}.
\end{displaymath}
Le choix des paramètres sera tel que les ensembles $\mathcal{P}_n$ (resp. $\mathcal{P}_{n, \mathbb{Z}}$) soient disjoints. On note: 
\begin{center}
$\mathrm{Tor}_{A,n}= \{ P$, points de $p$-torsion se réduisant sur $0$ mod $\mathfrak{q}$, pour $\mathfrak{q} / \mathfrak{p} / p$ et $\mathfrak{p} \in \mathcal{P}_{n} \}$, 
\end{center}
où $\mathfrak{q}$ est un idéal premier dans un corps de définition du point de torsion $P$. On note pour toute la suite $K'$ le corps engendré par $K$, par un corps de définition de $V$ et par les coordonnées projectives des points de la réunion: $\bigcup_{1 \leq n \leq r} \mathrm{Tor}_{A,n}$. Les fibrés et les morphismes qu'on va considérer seront tous définis sur $K'$; en fait, comme on considère des $\mathcal{O}_K$-fibrés, les calculs de pentes valables sur $\mathcal{O}_K$ se transportent sans changer à $\mathcal{O}_{K'}$, et le corps $K'$ intervient uniquement dans les estimations ultramétriques.

\bigskip

On ordonne les points de $\mathrm{Tor}_{A,n}$ en les classant d'abord selon le plus petit premier $p$ de torsion (pour l'ordre naturel sur $\mathbb{Z}$) puis en choisissant arbitrairement un ordre à $p$ fixé. Pour $i=(i_1, \ldots, i_r)$, où $i_n$, pour $1 \leq n \leq r$,  est un indice dans $[1, |\mathrm{Tor}_{A,n}|] $, on note $P_{i}$ le point $P_{i_1}+ \cdots + P_{i_r}$. On note $I$ l'ensemble de ces multi-indices, qu'on ordonne avec l'ordre lexicographique. On confondra dans la suite l'ensemble $I$ et son image dans $\mathbb{N}$ par l'indexation, et les éléments $i \in I$ pourront être vus comme des entiers via cette identification. 
On construit, pour $i \in I$, une suite de fermés $X_{i}=X+P_{i}$. On pose enfin, pour $i \in I$: $T_{(i)}=T_{n_{i}}$ et $N_{(i)}=N_{n_i}$, où $n_{i}$ est le plus grand des $j$ tels que $P_{i_j} \neq 0$. 

\bigskip

On définit maintenant $S$, un schéma, et $F$ un espace vectoriel (de sections):
\begin{displaymath}
S=\bigcup_{i \in I} V(X_{i}, t_{A_M}, T_{(i)}) \mathrm{ \ \ \ et: \ } F=H^0(S, L_{M}); 
\end{displaymath}
et pour $k \in I$ un entier, on pose:
\begin{displaymath} 
S_k=\bigcup_{i\leq k} V(X_{i}, t_{A_M}, T_{(i)}).
\end{displaymath}
Ceci permet de définir $F_k$, le noyau du morphisme de restriction: 
\begin{displaymath}
p_k: H^0(S, L_{M}) \rightarrow H^0(S, L_{M})_{|S_k}. 
\end{displaymath}
La suite décroissante des $F_k$ est une filtration de $F$ et on en déduit une filtration de $E$ en posant: $E_k:=\phi^{-1}(F_{k-1})$. 
Soit enfin $G_k= F_{k-1}/F_k$. En appliquant le lemme des serpents aux deux suites exactes: 
\begin{displaymath}
0 \rightarrow F_k \rightarrow H^0(S, L_{M}) \rightarrow H^{0}(S, L_{M})_{|S_k} \rightarrow 0
\end{displaymath}
et:
\begin{displaymath}
0 \rightarrow F_{k-1} \rightarrow H^0(S, L_M) \rightarrow H^{0}(S, L_M)_{|S_{k-1}} \rightarrow 0 ,
\end{displaymath}
on voit que $G_k$ s'identifie au noyau de l'application:
\begin{displaymath}
H^{0}(S, L_M)_{|S_k} \rightarrow H^{0}(S, L_M)_{|S_{k-1}}.
\end{displaymath}
Ce noyau est constitué de sections sur $S_k$ nulles sur $S_{k-1}$, donc est un sous-espace vectoriel de: 
\begin{displaymath}
H_k :=H^0 \Big(X_{k}, \mathrm{Sym}^{T_{(k)}}(t_{A_M}^{\vee}) \otimes L_M \Big) \simeq \mathrm{Sym}^{T_{(k)}}(t_{A_M}^{\vee}) \otimes H^0 (X_{k}, L_M), 
\end{displaymath}
puisque le fibré $\mathrm{Sym}^{T_{(k)}}(t_{A_M}^{\vee})$ est constant. On choisit un modèle entier $H^0 (\mathcal{X}_{k}, \mathcal{L}_M)$ de $H^0 (X_{k}, L_M)$ et sa métrique, ainsi qu'un modèle entier $t_{\mathcal{A}}$ de $t_A$ comme en \ref{tangent}; on en déduit un modèle entier $t_{\mathcal{A}_M}$ de $t_{A_M}$ (voir \ref{derivations}); on dispose aussi d'une métrique sur $t_{\mathcal{A}_M}$ et on en déduit une métrique sur $\mathrm{Sym}^{T_{(k)}}(t_{\mathcal{A}_M}^{\vee})$ par passage au dual, tensorisation et quotient. \\
\\
\textbf{Remarque} On a choisi ici de noter $\mathrm{Sym}$ le fibré symétrique, pour éviter toute confusion avec le schéma $S$.

\bigskip

On obtient donc un modèle entier $\mathcal{G}_k$ de $G_k$ et par restriction des métriques un fibré hermitien $\overline{\mathcal{G}_k}$; sa pente maximale est majorée par le corollaire \ref{pente sous variete}: 
\begin{lem}
Il existe une constante $c_{12}$ ne dépendant que de $A$ telle que:
\begin{displaymath}
\hat{\mu}_{\mathrm{max}}(\overline{\mathcal{G}_{k}}) \leq c_{12} \Big(M^2 \hat{\mu}_{\mathrm{ess}}(X) + \mathrm{log} \big(\mathrm{deg}(X)\big) + T_{(k)}\mathrm{log}(M)\Big).
\end{displaymath} 
\end{lem}
\begin{dem}
L'injection: 
\begin{displaymath}
\mathcal{G}_k \hookrightarrow \mathcal{H}_k := \mathrm{Sym}^{T_{(k)}}(t_{\mathcal{A}_M}^{\vee}) \otimes H^0 (\mathcal{X}_k, \mathcal{L}_M)
\end{displaymath}
est isométrique et on a, par définition de la pente maximale:
\begin{displaymath}
\hat{\mu}_{\mathrm{max}}(\overline{\mathcal{G}_{k}}) \leq \hat{\mu}_{\mathrm{max}} ( \overline{\mathcal{H}_k} ).
\end{displaymath}
On peut donc appliquer le corollaire \ref{pente sous variete} avec les fermés $X_{k}$. En notant $h_M$ la hauteur projective associée à $L_M$, on commence par remarquer qu'il existe une constante $c_{12}$ ne dépendant que de $A$ telle que: 
\begin{eqnarray*}
\frac{h_M(X_{k})}{\mathrm{deg}_M(X_{k})} & \leq & \frac{\hat{h}_M(X_{k})}{\mathrm{deg}_M(X_{k})}+c_{12} \\ 
 & \leq & (M^2+1) \frac{\hat{h}(X_{k})}{\mathrm{deg}(X_{k})}+c_{12}, \\
\end{eqnarray*}
en comparant la hauteur projective et la hauteur canonique dans le plongement étiré (voir \cite{David-Philippon02}, propositions 3.9 et 3.14), puis grâce au lemme \ref{second etirement}. L'inégalité des minima successifs (démontrée par Zhang dans \cite{Zhang95}, theorem 5.2) donne donc, comme $M \geq 1$:
\begin{displaymath}
\frac{h_M(X_{k})}{\mathrm{deg}_M(X_{k})}  \leq 2g M^2\hat{\mu}_{\mathrm{ess}}(X) +c_{12}. 
\end{displaymath}
On remarque que le minimum essentiel, tout comme le degré, n'est pas modifié lorsqu'on translate par des points d'ordre fini.
Par choix de la base de dérivation sur $\mathcal{A}_M$, en tenant compte de l'action de $[M]$ sur le tangent et en bornant la hauteur de la base de dérivations algébriques par une constante ne dépendant que de $A$, on peut choisir:
\begin{displaymath}
c(t_{\mathcal{A}_M}) \leq \mathrm{log}(M) +c_{12}. 
\end{displaymath}
On a donc:
\begin{displaymath}
\hat{\mu}_{\mathrm{max}}(\overline{\mathcal{H}_{k}}) \leq c_{12} \Big(M^2 \hat{\mu}_{\mathrm{ess}}(X)+ \mathrm{log} \big(\mathrm{deg}(X)\big) + T_{(k)} \mathrm{log}(M) \Big),
\end{displaymath} 
quitte à prendre $c_{12}$ assez grande en fonction de $A$ pour obtenir ces deux inégalités. Le lemme est donc entièrement démontré.
\end{dem}

\bigskip

On note $\phi$ le morphisme de restriction de $E$ vers $F$ et on pose, pour $k \in I$: $\phi_k: E_k \rightarrow G_k$ l'application linéaire déduite de $\phi$ et de la projection canonique: $F_{k-1} \rightarrow G_k$. Pour être en mesure d'écrire l'inégalité de pentes, il restera à démontrer que $\phi$ est injectif, ce qu'on fera à l'aide d'un lemme de zéros. 

\subsection{Choix des paramètres et hypothèse sur le minimum essentiel}
\label{parametres}

Le choix des paramètres, auquel nous procédons maintenant, doit permettre d'assurer l'injectivité du morphisme $\phi$ et de contredire l'inégalité des pentes. La stratégie qui guide ce choix est la suivante:\\

\begin{itemize}
\item on définit les paramètres $M$, $T_0$, $N_1$ de telle manière que le terme (correspondant au fermé $X$ non-translaté): 
\begin{displaymath}
\mathrm{rg}(\mathcal{G}_0)\big( \widehat{\mu}_{\mathrm{max}}(\overline{\mathcal{G}_{0}}) + h(\phi_{0}) \big),
\end{displaymath}
soit inférieur à $\widehat{\mathrm{deg}}(\overline{\mathcal{E}_M})$.\\
\item les relations entre les paramètres $T_i$ et $N_i$ sont telles que les contributions des termes suivants dans l'inégalité des pentes soient négatives.  \\
\item le paramètre d'étirement $M$ est pris aussi petit que possible, ce qui permet de montrer que le morphisme $\phi$ est injectif; il est aussi choisi de telle sorte que $M^2\hat{\mu}_{\mathrm{ess}}(X)$ soit majoré par une constante, ce qui détermine la minoration obtenue pour le minimum essentiel.
\end{itemize}

\bigskip

On commence par introduire l'indice d'obstruction avec poids $\omega(x,X)$; celui-ci permet classiquement de prendre en compte la hauteur des dérivées dans un lemme de Siegel. Il aura un emploi similaire dans le cadre de la théorie des pentes.
\begin{defn}
Soit $X$ un fermé de Zariski propre de $A$ et $x$ un réel positif.
On pose:
\begin{displaymath}
\omega(x, X)= \mathrm{inf} \{(x
\mathrm{deg}(Z))^{1/\mathrm{codim}(Z)} \},
\end{displaymath}
où l'infimum porte sur l'ensemble des fermés de Zariski (lisses, équidimensionnels) propres de $A$ contenant $X$.
\end{defn}
On utilisera souvent le lemme suivant, qui compare l'indice d'obstruction simple, ne prennant en compte que les hypersurfaces, et l'indice d'obstruction avec poids: 
\begin{lem}
\label{obstruction avec poids}
Soit $X$ un fermé propre (lisse et équidimensionnel) de $A$ de codimension $r$. Il existe une constante $c_{13}$ ne dépendant que de $A$ telle que pour tout réel $x$ positif:
\begin{displaymath}
c_{13} x^{1/r}\omega(X) \leq \omega(x,X) \leq x \omega(X).
\end{displaymath}
\end{lem}
\begin{dem}
L'inégalité de droite est claire et celle de gauche est une conséquence d'un résultat de Chardin (\cf \cite{Chardin90}, corollaire 2, page 8 et exemple 1, page 9). Remarquons qu'un fibré très ample étant fixé, seule la dimension du projectif dans lequel on plonge $A$ intervient dans $c_{13}$.
\end{dem}

\bigskip

On introduit maintenant la constante $C_0$ annoncée à la fin de l'introduction; il s'agit d'une constante ne dépendant que de $A$, grande devant toutes les constantes du problème, dans un sens explicitable. Soit $\Delta$ le paramètre: 
\begin{displaymath}
\Delta= C_0² \mathrm{log} \big(3 \mathrm{deg}(V)\big).
\end{displaymath}
C'est à partir de ce paramètre, qui est {\it grosso modo} de l'ordre de log$\big(\mathrm{deg}(V)\big)$, qu'on va définir tous les autres paramètres. Son avantage -comparé à $\mathrm{log} \big( \mathrm{deg} (V) \big)$- est d'être inconditionnellement grand devant les constantes intervenant au cours de la preuve, par choix de $C_0$.

En vue de la phase de descente, qui clôt la preuve, nous introduisons deux paramètres, qui permettront d'itérer toute la construction. Soit donc $R$ un réel positif et $\rho$ un entier vérifiant les hypothèses suivantes:
\begin{displaymath}
\Delta \geq  \mathrm{log} (R)  \ \  \mathrm{et:} \ \ 1 \leq \rho \leq \big(9(2r)^{r+1}\big)^{r-1}.
\end{displaymath}

Dans le formalisme des pentes, on regarde la restriction d'une section aux translatés de $V$ par $n$ points de torsion, où $0 \leq n \leq r$. On prend un paramètre $T_0$ correspondant à la multiplicité initiale et on associe à tout $n \in [1;r]$ des paramètres spécifiques, à savoir une multiplicité $T_{n}$, et une borne $N_{n}$ pour la norme des premiers de torsion, ce qui détermine un ensemble de premiers $\mathcal{P}_n$. Le lien entre ces paramètres est choisi de telle sorte que la contribution des termes de pente, pour des indices strictement positifs, soient négatifs. On prend d'abord:
\begin{displaymath}
N_n=\Delta^{\rho (2r)^{r+2-n}}.
\end{displaymath}
Chaque $N_n$ est donc négligeable devant le précédent; les raisons de ce choix seront plus claires au cours de la preuve de la proposition \ref{quasi injectivite}, qui montrera {\it quasiment} l'injectivité du morphisme $\phi$. Passons aux paramètres de multiplicités. Rappelons qu'on a posé: $n_{\alpha}=1$ si $\alpha=g$ et $n_{\alpha}=2$ sinon (voir \ref{coefficient rang}). Le choix de cet indice correspond à la propriété métrique obtenue à la fin de la partie \ref{partie metrique}. Si $\alpha=g$, toutes les directions du tangent sont associées à une propriété métrique en $p^{-1/p}$; sinon, il existe des directions associées à des propriétés en $p^{-1/p^2}$. Par récurrence descendante, pour $0 \leq n \leq r-1$, on pose:
\begin{displaymath}
T_n=  T_{n+1} \big(N_{n+1}^{n_{\alpha}} \Delta^2\big).
\end{displaymath}
Puis on pose: 
\begin{displaymath}
T_r= 1.
\end{displaymath}
Ces formules déterminent complètement $T_0$ par itération: 
\begin{displaymath}
T_{0}=\Delta^{2r} \big( N_1 \cdots N_r \big)^{n_{\alpha}}.
\end{displaymath}

On finit par le paramètre $M$. Le but est de prendre $M^2$ assez grand pour que le premier terme dans la somme de l'inégalité de pentes soit plus petit que le degré de $\overline{\mathcal{E}_M}$. Cette condition s'apparente à celle qu'on obtiendrait par un lemme de Siegel dans une preuve classique de transcendance. Pour montrer l'injectivité de $\phi$, on aura au contraire besoin que $M^2$ ne soit pas trop grand. On choisit donc: 
\begin{displaymath}
M=\bigg[ \Big(T_0 \omega(\Delta T_0, X) \Big) ^{1/2}\bigg]+1.
\end{displaymath}

\bigskip

On rappelle que $V$ est une sous-variété propre de $A$ qui n'est pas incluse dans un translaté de sous-variété abélienne propre, et que: $X=V + \Sigma_0$, où $\Sigma_0$ est un sous-groupe fini de $A$. On suppose que le cardinal de $\Sigma_0$ n'est pas trop gros: 
\begin{eqnarray}
\label{sigma zero}
\mathrm{log}(|\Sigma_0|) \leq \Delta.
\end{eqnarray}
Puis on fait l'hypothèse suivante sur le minimum essentiel de $X$:
\begin{eqnarray}
\label{minimum essentiel}
\hat{\mu}_{\mathrm{ess}}(X) < \frac{\Delta}{T_0 \omega(\Delta T_{0}, X)}.
\end{eqnarray}
Notre but désormais sera d'obtenir une contradiction. Ceci fait, un rapide calcul nous donnera le théorème \ref{theoreme}.

\section{Utilisation de l'inégalité des pentes}
\label{calcul fibres}

On a déjà calculé le degré de $\overline{\mathcal{E}_M}$ et les pentes maximales. On veut maintenant majorer le rang des $\mathcal{G}_k$ (via les $\mathcal{H}_k$) et la hauteur des $\phi_k$. Le calcul du rang des $\mathcal{H}_k$ correspond au calcul du rang du système linéaire dans un lemme de Siegel. Les estimations ultramétriques dans la hauteur des $\phi_k$ sont le point crucial de la preuve et ont été préparées par la partie \ref{partie metrique}. Les estimations archimédiennes, enfin, utilisent essentiellement l'inégalité de Cauchy (en plusieurs variables).

\subsection{Majoration du rang}

Pour majorer le rang des $\mathcal{G}_k$, on doit d'abord compter les dérivations puis estimer la fonction de Hilbert d'un fermé lisse et équidimensionnel assez précisément. Remarquons que dans notre cas, c'est une majoration du rang pour $k=0$ qui est cruciale.

On rappelle qu'on a l'injection suivante:
\begin{displaymath}
\mathcal{G}_k \hookrightarrow \mathcal{H}_k = \mathrm{Sym}^{T_{(k)}}(t_{\mathcal{A}_M}^{\vee}) \otimes H^0 (\mathcal{X}_k, \mathcal{L}_M).
\end{displaymath}
On peut donc se contenter de majorer le rang du fibré de droite, qui est plus explicite.

\subsubsection{Nombre de dérivations}

Le calcul du nombre de dérivations est classique dans les preuves de transcendance. C'est un raisonnement de géométrie différentielle reposant sur l'``astuce de Philippon-Waldschmidt'', qu'on trouve par exemple dans \cite{Amoroso-David03}, lemme 2.5, et \cite{David-Hindry00}, 5.3, dans un contexte abélien. 

\bigskip

Si $X$ est un fermé lisse et équidimensionnel, $L$ un fibré ample sur $X$ et $n$ un entier, on note $h^0(X, L^{\otimes n})$ la fonction de Hilbert en degré $n$ associée à $X$ et $L$. La prise en compte de l'indice d'obstruction avec poids nous amène à faire le raisonnement qui suit avec un fermé lisse et équidimensionnel $Z_k$ contenant $X_k$, de codimension $r' \leq r$ ($r$ étant la codimension commune à $V$ et à tous ses translatés).

\begin{prop}
On a l'inégalité suivante pour le rang de $\mathcal{H}_k$:
\begin{displaymath}
\mathrm{rg}(\mathcal{H}_k) \leq T_{(k)}^{r'} h^0\big(Z_{k}, L_M \big).
\end{displaymath}
\end{prop}
\begin{dem}
Commençons par remarquer que l'inclusion: $X_k \subset Z_k$, nous permet de nous ramener à la majoration du rang du fibré  $\mathcal{G}_k$, dans lequel on a remplacé $X_k$ par $Z_k$. Quitte à permuter les indices, on peut ensuite supposer que les dérivations $(\partial_{1,x} \ldots, \partial_{r',x})$ se projettent sur une base du cotangent pour $x \in U$, un ouvert dense de $Z_{k}$. Cet ouvert est défini par la non-annulation du déterminant d'une sous-matrice de taille maximale de la matrice jacobienne sur un ouvert affine.

\bigskip

On note $\mathcal{I}^{(k)}$ le faisceau d'idéaux associé au voisinage infinitésimal $V(Z_{k}, t_{A_M}, T_{(k)})$. Pour $s$ une section sur $Z_{k}$, on va démontrer que:
\begin{eqnarray}
\label{astuce}
  \Big(\forall x \in U, \ \forall \lambda \in \mathbb{N}^{r'} \mathrm{\ tel \ que \ } \sum_{i \leq r'} \lambda_i \leq T_{(k)} \ : \ \partial_{x}^{\lambda} s = 0  \Big)
  & \Longrightarrow & s \in \mathcal{I}^{(k)},  
\end{eqnarray}
où on plonge $\mathbb{N}^{r'}$ dans $\mathbb{N}^g$ en complétant
les $r'$-uplets par des zéros. 

\bigskip

Faisons une premi\`ere r\'ecurrence
sur $T_{(k)}$. Le r\'esultat est vrai pour $T_{(k)} = 0$. On le suppose vrai
au rang $T_{(k)} - 1 \geq 0$. Soit $s \in H^0(Z_{k}, L_M)$ vérifiant
l'assertion de gauche dans (\ref{astuce}). Soit $(\lambda_1, \ldots,
\lambda_g)$ tels que:
\begin{displaymath}
\sum_{i \leq g} \lambda_i \leq T_{(k)}.
\end{displaymath}
On fait une nouvelle r\'ecurrence sur la longueur résiduelle:
\begin{displaymath}
|\lambda|'= \lambda_{r'+1}+ \ldots + \lambda_{g},
\end{displaymath}
pour montrer que $\partial^{\lambda}_{x} s = 0$ si $x \in U$. Pour $|\lambda|' = 0$, l'assertion de gauche
dans (\ref{astuce}) implique que: $\partial_x^{\lambda} s = 0$
pour tout $x \in U$. Supposons le r\'esultat d\'emontr\'e pour
$|\lambda|' - 1 \geq 0$ et soit $x \in U$. Quitte à renuméroter les
indices, comme $|\lambda|' \geq 1$, on peut écrire:
\begin{displaymath}
\partial^{\lambda}_{x} s = \gamma \partial_{g, x}
  \circ \partial^{\mu}_{x} s ;
\end{displaymath}
pour une certaine constante $\gamma$ (correspondant aux factorielles dans la d\'efinition
des $\partial_i)$. Comme $(\partial_{1, x}, \ldots, \partial_{r',x})$ se projettent sur une base du cotangent de $Z_{k}$ en $x$, on d\'ecompose:
\begin{displaymath}
\partial_{g,x} = \delta_x +
  \delta'_x,
\end{displaymath}
o\`u $\delta_x$ est une combinaison linéaire des $\partial_{i,x}$ ($1 \leq i \leq r'$) et $\delta'_x$ est un vecteur du tangent de $Z_{k}$,
d\'ependant de $x$. On a:
\begin{displaymath}
\delta_x \circ \partial^{\mu}_{x} s =0,
\end{displaymath}
par hypoth\`ese de deuxi\`eme r\'ecurrence et: $\delta'_x \circ
\partial^{\mu}_{x} s = 0$. En effet, $v \rightarrow
\partial^{\mu}_{v} s$ est une section de $\mathcal{I}^{(k)}$ par hypoth\`ese de premi\`ere r\'ecurrence et l'égalité suit de
la d\'efinition du tangent. Ainsi, $\partial^{\lambda}_{x} s = 0$ pour tout $x \in U$ donc sur $Z_{k}$ tout entier car $U$
est Zariski-dense dans $Z_k$. L'implication (\ref{astuce}) est donc vérifiée.

\bigskip

On obtient ainsi une majoration du nombre de conditions
linéaires à écrire pour être dans $\mathcal{I}^{(k)}$. Puisque le nombre de $r'$-uplets
d'entiers de somme $i$ est égal au terme binomial
$\Big(^{i+r'-1}_{\ \ r'-1}\Big)$ et puisque la dérivation est de degré $1$ sur les idéaux homogènes:
\begin{displaymath}
\mathrm{rg}(\mathcal{H}_k) \leq \sum_{i=0}^{T_{(k)}-1} \Big(^{i+r'-1}_{\ \
r'-1} \Big) h^0(Z_{k},L_M) \leq \Big(^{T_{(k)}+r'-1}_{\ \ \
r'-1}\Big) h^0(Z_{k}, L_M).
\end{displaymath}
On a enfin la majoration: 
\begin{displaymath}
\Big(^{T_{(k)}+r'-1}_{\ \ \
r'-1}\Big) \leq  T_{(k)}^{r'},
\end{displaymath} 
ce qui finit de prouver la proposition. 
\end{dem}
\textbf{Remarque}
En suivant une approche légèrement différente, on peut chercher des résultats métriques aussi précis que possible, dans les cas intermédiaires où le $p$-rang est différent de $0$ ou $g$, et mettre en oeuvre la méthode des pentes avec un système de {\it multiplicités penchées}. Le calcul du rang ci-dessus s'adapte bien à ce contexte, et on peut trouver, par exemple sous l'hypothèse que $A$ est simple, des dérivations $(\partial_{1,x} \ldots, \partial_{r',x})$ se projetant sur une base du cotangent de $Z_{k}$ sur un ouvert Zariski-dense, correspondant aux directions du tangent les moins ``poussées'' (associées aux directions de meilleures propriétés $p$-adiques). On utilise pour cela le fait que $V$, donc $X$ et les $X_k$, ne sont pas incluses dans un translaté de sous-variété abélienne propre. Ceci permet d'avoir une majoration satisfaisante pour le rang. Il semble cependant difficile d'avoir une bonne majoration de la pente maximale de $\overline{\mathcal{H}_k}$ dans ce cas.

\subsubsection{Majoration de la fonction de Hilbert géométrique}

Pour majorer la fonction de Hilbert de la variété $Z_{k}$, on ne peut pas se contenter du théorème de Hilbert-Samuel géométrique. Celui-ci donne une estimation asymptotique, ce qui oblige à introduire une constante indéterminée dépendant de $Z_{k}$ (et par suite de l'indice d'obstruction), ce qu'on souhaite éviter. On dispose d'une estimation inconditionnelle démontrée par Chardin, qui donne immédiatement: 

\begin{prop}
On a:
\begin{displaymath}
h^0(Z_{k}, L_M) \ \leq \ \big( ^{1+g-r'}_{\ \ g-r'} \big)
\ \mathrm{deg}_M(Z_{k}).
\end{displaymath}
\end{prop}
\begin{dem}
C'est une conséquence du r\'esultat principal de \cite{Chardin89}, valable pour un fermé lisse et équidimensionnel.
\end{dem}
On regroupe les deux derniers résultats dans le corollaire suivant:
\begin{cor}
\label{rang}
\begin{displaymath}
\mathrm{rg} (\mathcal{G}_k) \leq \frac{g(2M^2)^g}{ \Delta T_{0}}.
\end{displaymath}
\end{cor}
\begin{dem}
On commence par combiner les deux dernières propositions en remarquant que le coefficient binômial apparaissant dans la dernière est inférieur ou égal à $g$. Et on obtient, grâce au lemme \ref{second etirement}:
\begin{displaymath}
\mathrm{rg} (\mathcal{H}_k) \leq g T_{(k)}^{r'} \mathrm{deg}_M(Z_{k}) \leq g T_{(k)}^{r'}(2M^2)^{g-r'} \mathrm{deg}(Z_{k}).
\end{displaymath}
Ceci est vrai pour toute variété $Z_{k}$ contenant $X_{k}$. Soit donc $Z_{k}$ une sous-variété, de codimension $r'$, telle que (on se ramène à $X$ car le degré est invariant par translation): 
\begin{displaymath}
\omega(\Delta T_{0}, X)= \Big( \Delta T_{0} \mathrm{deg} (Z_{k}) \Big)^{1/r'}.
\end{displaymath} 
On a donc, par choix du paramètre $M$:
\begin{displaymath}
\mathrm{rg} (\mathcal{H}_k) \leq \frac{g(2M^2)^g}{\Delta T_{0}} \bigg( \frac{T_{(k)}} {T_0} \frac{\omega(\Delta T_{0}, X)} {\omega(\Delta T_{0} , X)} \bigg)^{r'},
\end{displaymath}
et le résultat suit puisque $T_{(k)} \leq T_0$.
\end{dem}

\subsection{Normes ultramétriques des morphismes}

Soit $k \in I$ un entier, et $\mathfrak{p}$ un premier de $\bigcup_{1 \leq n \leq r}\mathcal{P}_{n}$ au-dessus d'un nombre premier $p$. Pour $\mathfrak{q} / \mathfrak{p}$, un idéal premier de $K'$, la norme $\mathfrak{q}$-adique de $\phi_k$ n'est correctement majorée que si la $k$-ème étape dans la filtration correspond à la translation par un point de $p$-torsion. Autrement, on se contente d'une majoration simple. On introduit donc la définition suivante:
\begin{defn}
On dit que $k$ est $n$-lié à $\mathfrak{q}$ si le point $P_{k}=P_{k_1}+ \cdots + P_{k_{n}}$ associé à $k$ est tel que $P_{k_n}$ soit un point de $p$-torsion non-nul se réduisant sur $0$ modulo $\mathfrak{q}$, et $\mathfrak{q}/p$. 
\end{defn}  
La proposition qui suit est une conséquence du lemme $p$-adique de la partie 3.
\begin{prop}
\label{ultrametrique}
Pour tout $k \in I$, et tout idéal premier $\mathfrak{q} / \mathfrak{p}$ de $K'$, on a:
\begin{displaymath}
\mathrm{log} \left\| \phi_k \right\|_{\mathfrak{q}} \leq 0.
\end{displaymath}
De plus, si $k$ est $n$-lié à $\mathfrak{q}$, on a: 
\begin{displaymath}
\mathrm{log} \left\| \phi_k \right\|_{\mathfrak{q}} \leq - n_{\mathfrak{q}}\big(T_{n-1}-T_{n}\big) \frac{\mathrm{log}p}{p^{n_{\alpha}}}.
\end{displaymath}
\end{prop}
\textbf{Remarque} On rappelle que $n_{\alpha}$ vaut $1$ si $\alpha=g$ (densité positive de premiers de réduction ordinaire) et vaut $2$ sinon (voir \ref{coefficient rang}).\\
\\
\begin{dem}
Soit $k \in I$ et $\mathfrak{q}/ \mathfrak{p}$ un idéal premier de $K'$. Soit $s$ un élément de $E_k$ tel que: $\left\|s\right\|_{\mathfrak{q}} \leq 1$; soit $\mathcal{X}_{k}$ le modèle entier de $X_{k}$ choisi dans la définition de $\mathcal{G}_k$, et $x$ un élément de $\mathcal{X}_{k}$. Soit aussi $\lambda \in \mathbb{N}^g$ tel que $\sum_{1 \leq i \leq g} \lambda_i \leq T_{(k)}$. On a: 
\begin{displaymath}
\partial^{\lambda}_x s = \partial^{\lambda} (s \circ \tau_x).
\end{displaymath} 
Par définition de  $\mathcal{P}_{A}$, les coefficients de $\tau_x$ sont $\mathfrak{p}$-entiers (voir (\ref{condition 1})). De plus, le choix de $\mathcal{P}_{A}$ permet aussi que l'opérateur différentiel $\partial^{\lambda}$ appliqué en chaque fonction abélienne s'exprime comme un polynôme en les fonctions abéliennes à coefficients $\mathfrak{p}$-entiers (par (\ref{condition 2})). En remarquant que la base de dérivations sur $A$ est algébrique (ce qui fait disparaître les factorielles au dénominateur quand on dérive des polynômes) et que la base de dérivations sur $A_M$ est obtenue par combinaisons linéaires à coefficients entiers des dérivations sur $A$, on a: 
\begin{displaymath}
|\partial^{\lambda}_x s |_{\mathfrak{q}} \leq 1,
\end{displaymath}
ceci étant vrai pour tout $x \in \mathcal{X}_{k}$ et tout $\lambda \in \mathbb{N}^g$ tel que $\sum \lambda_i \leq T_{(k)}$, on en déduit: 
\begin{displaymath}
\left\|\phi_k(s)\right\|_{\mathfrak{q}} \leq 1. 
\end{displaymath}
La première partie de la proposition est donc démontrée en passant au sup sur $s \in E_k$. 

\bigskip

Supposons que $k$ soit $n$-lié à $\mathfrak{q}$. On désigne par $\textbf{t}=(t_1, \ldots, t_g)$ une base de paramètres correspondant à une base orthonormale pour la forme de Riemann sur $t_A$. On reprend les notations convenues avant la proposition \ref{proposition metrique}. L'isomorphisme: 
\begin{displaymath}
\widehat{\mathcal{A}}_{\mathcal{O}_\mathfrak{p}} \simeq \mathcal{O}_{\mathfrak{p}}[[t_{1}, \ldots, t_{g}]],
\end{displaymath}
induit une application:
\begin{displaymath}
H^0(\mathcal{A}_{\mathcal{O}_{\mathfrak{p}}}, \mathcal{O}_{\mathcal{A}_{\mathcal{O}_{\mathfrak{p}}}}) \rightarrow \mathcal{O}_{\mathfrak{p}} [[t_1, \ldots, t_g]].
\end{displaymath}
Celle-ci associe à l'image d'une section $f$ de $H^0(\mathcal{A}_{\mathcal{O}_{\mathfrak{p}}}, \mathcal{O}_{\mathcal{A}_{\mathcal{O}_{\mathfrak{p}}}})$ son développement de Taylor en $0$:
\begin{displaymath}
\sum_{\mu \in \mathbb{N}^g} (\partial^{\mu} f ) \ \textbf{t}^{\mu}.
\end{displaymath}

\bigskip

On écrit $P_k= P_{k'}+ P_{k_n}$, et pour $x \in \mathcal{X}_{k}$: $x=x'+P_{k_n}$. Le point de torsion $P_{k_n}$ se réduit sur $0$ modulo $\mathfrak{q}$. On en déduit, par la proposition \ref{proposition metrique}, les majorations suivantes:
\begin{displaymath}
\forall 1 \leq i \leq g:  |t_{i}(P_{k_n})|_{\mathfrak{q}} \leq p^{-n_{\mathfrak{q}}/p^{n_{\alpha}}}. 
\end{displaymath}
De plus, on a encore, pour tout $\mu \in \mathbb{N}^g$:
\begin{displaymath}
|\partial^{\mu} (s \circ \tau_{x'})|_{\mathfrak{q}} \leq 1.
\end{displaymath}
La norme étant ultramétrique, ceci suffit à voir que la série de Taylor de $\partial^{\lambda} (s \circ \tau_{x'})$ converge; on en déduit l'égalité:
\begin{displaymath}
\partial^{\lambda}_x s = \sum_{\mu \geq \lambda \in \mathbb{N}^g} \partial^{\mu} (s \circ \tau_{x'}) \textbf{t}(P_{k_n})^{\mu-\lambda},
\end{displaymath}
où: 
\begin{displaymath}
\mu \geq \lambda \Leftrightarrow  \forall i \leq g: \ \mu_i \geq \lambda_i.
\end{displaymath}

Par définition de $E_k$, $s$ est nulle à un ordre $T_{n-1}$ en $x'$, donc si $\partial^{\mu} (s \circ \tau_{x'}) \neq 0$, on a: 
\begin{displaymath}
\sum_{i \leq g} \mu_i > T_{n-1}. 
\end{displaymath}
On en déduit: 
\begin{displaymath}
\left|\prod_{i=1}^g t_{i} (P_{k_n})^{\mu_i-\lambda_i}\right|_{\mathfrak{q}} \leq p^{-n_{\mathfrak{q}}(T_{n-1}-T_n)/p^{n_{\alpha}}}. 
\end{displaymath}
La norme étant ultramétrique, l'inégalité fine de la proposition en résulte.  
\end{dem}

\subsection{Normes archimédiennes des morphismes d'évaluation}

Soit $k$ un entier, et $\sigma$ une place archimédienne de $K'$. Tous les $\phi_k$ sont définis sur $K'$. Le but de ce paragraphe est de majorer $\left\| \phi_k \right\|_{\sigma}$, par des méthodes d'analyse complexe, en particulier l'inégalité de Cauchy. On suit pour cela le paragraphe 5.9 de \cite{Gaudron06}.  Cette majoration sera la même pour toutes les places archimédiennes; plus précisément, on obtient la proposition suivante:
\begin{prop}
\label{archimedien}
Il existe une constante $c_{16}$ ne dépendant que de $A$ telle que:
\begin{displaymath}
\frac{1}{[K': \mathbb{Q}]} \Big(\sum_{\sigma: K' \hookrightarrow \mathbb{C}} \mathrm{log} \left\|\phi_k \right\|_{\sigma} \Big) \leq c_{16} T_{(k)} \mathrm{log}(M).
\end{displaymath}
\end{prop}
\begin{dem}
Soit $s \in E_k \otimes \mathbb{C}$ et $x \in \mathcal{X}_{k}$.  On a alors:
\begin{displaymath}
\phi_k (s)(x) \in \big(S^{T_{(k)}} t_{\mathcal{A}_M}^{\vee} \otimes x^* \mathcal{L}_M \big)\otimes_{\sigma} \mathbb{C},
\end{displaymath}
qui est isomorphe (non isométriquement) à:
\begin{displaymath}
\mathrm{Hom}_{\mathbb{C}}(S^{T_{(k)}} t_{\mathcal{A}_{M, \sigma}}, x^* \mathcal{L}_{M, \sigma}).
\end{displaymath} 
Si on fixe une base d'ouverts affines, l'image de $\phi_k (s)(x)$ par cet isomorphisme est le morphisme qui associe à une dérivation $D$ d'ordre $T_{(k)}$ la valeur de $Ds$ en $x$ (la dérivation $D$ donnant par translation une dérivation en $x$). On note $\Theta_{T_{(k)}}$ cet isomorphisme qui est défini dans \cite{Gaudron06}, 4.1, où sa norme d'opérateur, ainsi que celle de son inverse, sont majorées: 
\begin{displaymath}
\left\| \Theta_{T_{(k)}} \right\|_{\sigma} \leq \mathrm{max}_{\mathbf{i} \in \mathbb{N}^g, |\textbf{i}|=T_{(k)}} \Big\{ \frac{T_{(k)}!}{\mathbf{i}!} \Big\} \mathrm{ \ et \ } \left\| \Theta_{T_{(k)}} ^{-1} \right\|_{\sigma} \leq 1.
\end{displaymath}
On en déduit que: 
\begin{displaymath}
\left\| \phi_k (s)(x) \right\|_{\sigma} \leq \left\| \Theta_{T_{(k)}}\big(\phi_k (s)(x)\big) \right\|_{\sigma}.
\end{displaymath}

\bigskip

Soit $(f_{1, \sigma}, \ldots, f_{g, \sigma})$ une base de $t_{\mathcal{A}_{M, \sigma}}$ (composée à partir d'une base orthonormée de $t_{\mathcal{A}}$ pour la forme de Riemann induite par $\sigma$), correspondant à des dérivations $(\partial_{1, \sigma}, \ldots, \partial_{g, \sigma})$. Soit $D$ une dérivation d'ordre $T_{(k)}$ le long de $t_{\mathcal{A}_{M, \sigma}}$. On écrit $D= \sum_{\textbf{i} \in \mathbb{N}^g, |\textbf{i}|=T_{(k)}} d_\textbf{i} \partial_{1, \sigma}^{i_1} \cdots \partial_{g, \sigma}^{i_g}$. La norme sur $S^{T_{(k)}} t_{\mathcal{A}_{M, \sigma}}$ est la norme quotient déduite de la projection: 
\begin{displaymath}
t_{\mathcal{A}_{M, \sigma}}^{\otimes T_{(k)}} \rightarrow S^{T_{(k)}} t_{\mathcal{A}_{M, \sigma}}. 
\end{displaymath}
On a donc: 
\begin{displaymath}
\left\| D \right\|^2= \sum_{|\textbf{i}|=T_{(k)}} |d_\textbf{i}|^2 \frac{\textbf{i}!}{T_{(k)}!} \geq \Big( \sum_{|\textbf{i}|=T_{(k)}} |d_\textbf{i}| \Big)^2 \times g^{-T_{(k)}}.
\end{displaymath}
De plus, on a: 
\begin{displaymath}
\left\| Ds(x) \right\| \leq \sum_{|\textbf{i}|=T_{(k)}} |d_\textbf{i}| \left\| \big(\partial_{1, \sigma}^{i_1} \cdots \partial_{g, \sigma}^{i_g}\big)s(x)\right\|,
\end{displaymath}
et on en déduit: 
\begin{displaymath}
\left\| \phi_k (s)(x) \right\|_{\sigma} \leq g^{T_{(k)}/2} \mathrm{max}_{|\textbf{i}|=T_{(k)}} \Big\{ \left\| \big(\partial_{1, \sigma}^{i_1} \cdots \partial_{g, \sigma}^{i_g}\big)s(x) \right\| \Big\}.
\end{displaymath}
En reprenant la définition de la métrique cubiste, on peut se ramener à la fonction $\theta$ correspondant à $s$; notons que via le plongement étiré, la section $s$ est une section sur $A$ de degré $\leq 2M^2$. On note $\textbf{u}$ un logarithme de $\sigma(x)$ ayant une norme hermitienne minimale sur $\mathbb{C}^g$. En tenant compte de l'action de $\phi_M$ sur les dérivations, il vient:
\begin{displaymath}
\left\| \phi_k (s)(x) \right\|_{\sigma} \leq (Q\sqrt{g})^{T_{(k)}} e^{-3 \pi M^2 \left\| \textbf{u} \right\|_{\sigma}^2} \times \mathrm{max}_{|\textbf{i}|=T_{(k)}} \Big\{ \Big| \frac{1}{\textbf{i}!} \Big(\frac{\partial}{\partial z}\Big)^\textbf{i} \theta (\textbf{u} + \textbf{z})_{|\textbf{z}=0} \Big| \Big\}.
\end{displaymath}

\bigskip

De plus, par l'inégalité de Cauchy appliquée à $\theta$, pour tout réel $r_{C} > 0$ (à ne pas confondre avec $r$, qui désigne la codimension de $V$), on a:
\begin{displaymath}
\mathrm{max}_{|\textbf{i}|=T_{(k)}} \Big\{ \Big| \frac{1}{\textbf{i}!} \Big(\frac{\partial}{\partial z}\Big)^\textbf{i} \theta (\textbf{u} + \textbf{z})_{|\textbf{z}=0} \Big| \Big\} \leq \bigg(\frac{1}{r_C^{T_{(k)}}}\bigg) \mathrm{sup}_{\left\| \textbf{z} \right\|_{\sigma} \leq r_C} |\theta(\textbf{u}+\textbf{z})|.
\end{displaymath}
En revenant aux métriques cubistes, on trouve: 
\begin{displaymath}
\left\| \phi_k (s)(x) \right\|_{\sigma} \leq \bigg(\frac{M \sqrt{g}}{r_C}\bigg)^{T_{(k)}}e^{3 \pi M^2(r_C^2+2r_C \left\| \textbf{u} \right\|_{\sigma})} \left\| s \right\|_{\mathrm{sup}, \sigma}.
\end{displaymath}
On choisit alors $r_C$ de façon à optimiser la majoration:
\begin{displaymath}
r_C=\frac{\sqrt{g}}{M^2 \mathrm{max} \{ 1, \left\| \textbf{u} \right\|_{\sigma} \}}.
\end{displaymath}

\bigskip

Comme $M^2 \geq \sqrt{g}$ (par le choix des paramètres), cela donne: 
\begin{displaymath}
\left\| \phi_k (s)(x) \right\|_{\sigma} \leq \big(M^3 \mathrm{max} \{1, \left\| \textbf{u} \right\|_{\sigma} \} \big)^{T_{(k)}} e^{9 \pi g}\left\| s \right\|_{\mathrm{sup}, \sigma}.
\end{displaymath}
La norme sur $\mathcal{G}_k$ étant la norme de Löwner $\left\|. \right\|_L$ associée à la norme du sup, elle est plus petite que celle-ci et on a finalement: 
\begin{displaymath}
\left\| \phi_k (s) \right\|_{L, \sigma} \leq \big(M^3 \mathrm{max} \{1, \left\| \textbf{u} \right\|_{\sigma} \} \big)^{T_{(k)}} e^{9 \pi g}\left\| s \right\|_{\mathrm{sup}, \sigma},
\end{displaymath} 
Et, vu le choix de la norme $L^2$ sur $\mathcal{E}_M$: 
\begin{displaymath}
\left\| \phi_k  \right\|_{\sigma} \leq \big(M^3 \mathrm{max} \{1, \left\| \textbf{u} \right\|_{\sigma} \} \big)^{T_{(k)}} e^{9 \pi g} \mathrm{sup}_{f \in \mathcal{E}_M, s \neq 0} \bigg\{ \frac{ \left\|  s \right\|_{\mathrm{sup}, \sigma}}{\left\|  s \right\|_{\mathrm{L^2}, \sigma}} \bigg\}.
\end{displaymath} 

\bigskip

La remarque 4.18 de \cite{Gaudron06} montre que pour une certaine constante $c_{14}$ ne dépendant que de $g$: 
\begin{displaymath}
\frac{1}{[K':\mathbb{Q}]} \sum_{\sigma: K' \rightarrow \mathbb{C}} \mathrm{log} \mathrm{max} \{1, \left\| \textbf{u} \right\|_{\sigma} \} \leq c_{14} \mathrm{max} \{1, \mathrm{log}(h_F(A)), \mathrm{log} (h^0(A, L)) \}.
\end{displaymath}
De plus, d'après la même référence, lemme 4.16, il existe une autre constante $c_{15}$ telle que: 
\begin{displaymath}
\mathrm{log} \ \mathrm{sup}_{s \in \mathcal{E_M}_{\sigma}, s \neq 0} \bigg\{ \frac{ \left\|  s \right\|_{\mathrm{sup}, \sigma}}{\left\|  s \right\|_{\mathrm{L^2}, \sigma}} \bigg\} \leq c_{15} \mathrm{max}\{1, \mathrm{log}(h_F(A)), \mathrm{log}(h^0(A,L)) \} \mathrm{log} (M^2).
\end{displaymath}
La proposition suit en sommant sur les places archimédiennes, pour une constante $c_{16}$ dépendant de $A$.
\end{dem}
\textbf{Remarques} On aurait pu choisir différemment le paramètre $r_C$ issu de l'inégalité de Cauchy, en prenant à la place: $r_C'=r_C \times M^2$. On aurait alors obtenu une majoration par $c_{16}(T_{(k)}+M^2)$, qui aurait été plus mauvaise dans notre contexte puisque tous les paramètres sauf $M$ sont logarithmiques en $\mathrm{deg}(X)$.

Jusqu'à la phase de descente (où cela ne semble plus possible), on pourrait travailler avec des termes en  log$\Big(\omega(V)\Big)$ à la place de log$\Big(\mathrm{deg}(V)\Big)$, quitte à améliorer la proposition \ref{wirtinger}.

\subsection{Inégalité de pentes et conséquences}
\label{pentes consequences}

Pour pouvoir appliquer le théorème des pentes, on doit s'assurer que le morphisme $\phi$ est injectif. On y travaillera dans la partie suivante, et on suppose par avance ici cette injectivité.\\
\\
\textbf{Contradiction sous l'injectivité de $\phi$} \\
\\
On suppose par l'absurde que (\ref{minimum essentiel}) est vérifié. Puisqu'on a aussi supposé que $\phi$ est injectif, on peut donc utiliser l'inégalité des pentes du lemme \ref{lemme pentes}. On a donc: 
\begin{displaymath}
\widehat{\mathrm{deg}} \ \overline{\mathcal{E}_M} \leq \sum_{k \in I} \mathrm{dim}(E_k/E_{k+1})\big( \widehat{\mu}_{\mathrm{max}}(\overline{\mathcal{G}}_{k}) + h(\phi_{k}) \big).
\end{displaymath}

\bigskip

Soit $k>0$ un entier; on veut d'abord montrer que la contribution du $k$-ème terme dans la somme précédente est négative. On regroupe pour commencer les estimations faites sur le morphisme $\phi_k$, à savoir les majorations archimédiennes et ultramétriques. L'entier $k$ est par définition $n$-lié à tout idéal premier $\mathfrak{q}$ de $\mathcal{O}_{K'}$ au-dessus d'un seul premier $\mathfrak{p}$ de $\mathcal{P}_A$. On note $p$ son image dans $\mathbb{Z}$; on a alors (par les propositions \ref{ultrametrique} et \ref{archimedien}):
\begin{eqnarray*}
h(\phi_k) = \frac{1}{[K':\mathbb{Q}]} \sum_{v \in M(K')} \mathrm{log} \left\| \phi_k \right\|_v  & \leq  & c_{16} T_{n} \mathrm{log}(M) - \frac{1}{[K':\mathbb{Q}]} \sum_{\mathfrak{q} / \mathfrak{p}}  n_{\mathfrak{q}}\big(T_{n-1}-T_{n}\big)\frac{\mathrm{log}p}{p^{n_{\alpha}}} \\
& \leq & c_{16} T_{n} \mathrm{log}(M) - \frac{1}{[K:\mathbb{Q}]} \big(T_{n-1}-T_{n}\big)\frac{\mathrm{log}p}{p^{n_{\alpha}}}, \\
\end{eqnarray*}
où on a utilisé l'égalité classique (voir\cite{Lang86}, II, Corollary 1 to Theorem 2):
\begin{displaymath}
\sum_{\mathfrak{q} / \mathfrak{p}} n_{\mathfrak{q}}=[K':K]. 
\end{displaymath}
On majore ensuite:
\begin{eqnarray*}
\hat{\mu}_{\mathrm{max}}(\overline{\mathcal{G}_k}) + h(\phi_k) & \leq &  c_{17} \Big( M^2 \hat{\mu}_{\mathrm{ess}}(X) + T_n\mathrm{log} (M) + \mathrm{log} \big(\mathrm{deg}(X)\big) \Big) -  \frac{1}{[K:\mathbb{Q}]}T_{n-1}\frac{\mathrm{log}p}{p^{n_{\alpha}}} \\
 & \leq & c_{17} \Big( M^2 \hat{\mu}_{\mathrm{ess}}(X) + T_n \mathrm{log} (M) + \mathrm{log} \big(\mathrm{deg}(X)\big) \Big) -  \frac{1}{[K:\mathbb{Q}]}T_{n-1}\frac{\mathrm{log}N_{n}}{N_{n}^{n_{\alpha}}}, \\
\end{eqnarray*}
pour une constante $c_{17}$ ne dépendant que de $A$, car $p \leq N_n$ et la fonction log$(x)/x^{n_{\alpha}}$ est décroissante pour $x \geq 3$. Le choix du paramètre $M$ (\cf \ref{parametres})  et l'hypothèse (\ref{minimum essentiel}) sur le minimum essentiel donnent:
\begin{displaymath}
M^2 \hat{\mu}_{\mathrm{ess}}(X) \leq 2 \Delta \frac{ T_0 \omega( \Delta T_0, X)}{ T_0 \omega( \Delta T_0, X)} \leq 2 \Delta. 
\end{displaymath}
Puis:
\begin{eqnarray*}
\hat{\mu}_{\mathrm{max}}(\overline{\mathcal{G}_k}) + h(\phi_k) & \leq & 5c_{17} \Delta T_n - \frac{T_{n-1}}{N_n^{n_{\alpha}}} \leq 0. \\
\end{eqnarray*}
On a d'abord utilisé l'hypothèse (\ref{sigma zero}) sur $\Sigma_0$, puis on a majoré log$(M)$ par $\Delta$ grâce au lemme \ref{obstruction avec poids}, en comparant log$\Delta$ à $\Delta$ (quitte à prendre $C_0$ assez grand). On a ensuite éliminé $[K:\mathbb{Q}]$ avec log$N_n$ (pour $C_0$ assez grand là encore). La dernière inégalité résulte de la définition des $T_n$ par récurrence descendante (voir encore le paragraphe \ref{parametres}).

\bigskip

Il suit, grâce à l'estimation du terme de hauteur qu'on vient de faire (qui est encore valable pour $k=0$, sans le raffinement ultramétrique): 
\begin{displaymath}
\widehat{\mathrm{deg}} \ \overline{\mathcal{E}_M} \leq \mathrm{dim}(E_0/E_1)\big( \hat{\mu}_{\mathrm{max}}(\overline{\mathcal{G}_{0}}) + h(\phi_{0}) \big) \leq \mathrm{rg}(\mathcal{G}_0)(5 c_{17}\Delta T_0).
\end{displaymath}
On a pu remplacer dim$(E_0/E_1)$ par rg$(\mathcal{G}_0)$ en utilisant l'injectivité de $\phi$, et parce que la hauteur est majorée par un terme positif. 
En combinant le lemme \ref{degre e} (il est clair par le choix des paramètres que log$(M)$ est plus grand que n'importe quelle constante)et la proposition  \ref{rang} avec $k=0$, on obtient: 
\begin{displaymath}
c_{11}(M)^{2g} \mathrm{log}(M) \leq \frac{g(2M^2)^g}{\Delta T_{0}}(5 c_{17}\Delta T_0),
\end{displaymath}
et on rappelle que $c_{11} > 0$. Puis:
\begin{displaymath}
\mathrm{log} (M) \leq \frac{5g2^gc_{17}}{c_{11}}.
\end{displaymath}
On en déduit une contradiction, puisque log$(M)$ est plus grand que n'importe quelle constante du problème, par définition de $\Delta$. 

{\begin{flushright}$\Box$\end{flushright}}

\section[Injectivité du morphisme]{Lemme de zéros et injectivité du morphisme}

Il nous reste donc à nous assurer que le morphisme de restriction est injectif. On procède par l'absurde, en commençant par écrire un lemme de zéros, et le choix des paramètres fait dans \ref{parametres} doit mener à une contradiction. Ceci s'avère assez délicat, et on y parvient en deux étapes. La première est de nature combinatoire (paragraphe \ref{combinatoire}), et mène à une {\it quasi-contradiction} en \ref{quasi contradiction}; la seconde est un argument de descente sur des variétés (paragraphe \ref{descente 1}), qui imposera de travailler en petite codimension: $r \leq 2$.

\subsection{Lemme de zéros}

Le lemme de zéros dont on a besoin ici s'inscrit dans la tradition des théorèmes démontrés par P. Philippon, dont on reprend le formalisme (\cf \cite{Philippon95}). Ce lemme est l'analogue abélien du théorème utilisé dans \cite{Amoroso-David03}, à une différence près: on prend en compte les multiplicités, et ce à l'aide de la notion de dessous d'escalier, qui permet d'envisager des multiplicités variables et différentes selon les directions. Dans le cas qui nous intéresse, il est utile d'envisager la multiplicité finale dans le lemme de zéros, qui permet une légère amélioration par rapport au cas torique. On ne fait pas usage, cependant, de multiplicités différentes selon les directions.\\
\\
\textbf{Remarque} On a pris la multiplicité finale $T_r$ dans l'inégalité de pentes égale à $1$ mais la multiplicité finale dans le lemme de zéros est un certain $T_{r_0}$, pour $r_0 \leq r$, et n'est pas forcément nulle. 

\bigskip

Dans la suite, si $l \in \mathbb{N}$ et $Z$ est une sous-variété de $A$, on notera $\{l\}Z$ le cycle $Z$ avec la multiplicité $l$. Cette notation peu conventionnelle a pour but d'éviter toute confusion avec l'image de $Z$ sous la multiplication par $l$, notée $[l]Z$.

\bigskip

Faisons d'abord quelques rappels nécessaires pour écrire le lemme de zéros. Soit $A$ une variété abélienne munie d'une fibré ample $L$; on considère la base de dérivations sur $A$ définie en \ref{derivations}. Un ensemble $E \subset \mathbb{N}^g$ est un escalier si pour tout $\beta \in E$, on a $\beta + \mathbb{N}^g \subset E$. Un sous-ensemble de $\mathbb{N}^g$ est un dessous d'escalier s'il est le complémentaire d'un escalier. Si $W$ est le dessous d'un escalier $E$ de $\mathbb{N}^g$, et si on a des indices: $1 \leq i_1 < \cdots < i_d \leq g$, on note $\mathcal{C}_{i_1, \ldots, i_d}(W)$ l'enveloppe convexe dans $\mathbb{R}_+^d$ de la trace de $E$ sur la $d$-face de $\mathbb{N}^g$ définie par $(i_1, \ldots, i_d)$. 

\bigskip

On appelle aussi ensemble pondéré un sous-ensemble $\Sigma$ de $\mathbb{N}^g \times A$ tel que pour tout $x \in A$, l'ensemble $W_{x, \Sigma}=(\mathbb{N}^g \times \{ x \}) \cap \Sigma$ soit un dessous d'escalier (éventuellement vide). On appelle support de $\Sigma$, noté Supp$(\Sigma)$, sa projection sur $A$. Si $\Sigma$ et $\Sigma'$ sont deux ensembles pondérés, on définit $\Sigma+ \Sigma'$ comme l'ensemble des couples $(x+x', \lambda+ \lambda')$, pour $(x, \lambda) \in \Sigma$, et $(x', \lambda') \in \Sigma'$; c'est aussi un ensemble pondéré. On a $E+ \emptyset= \emptyset$ et si $E$ est un sous-ensemble de $A$, on l'identifie à l'ensemble pondéré $\{0\} \times E$.

\bigskip

On dit que $f \in H^0(A, L)$ s'annule sur un ensemble pondéré $\Sigma$ si pour tout $(x, \lambda) \in \Sigma$, on a: $\partial_x^{\lambda} f=0$.
Si $V$ est une sous-variété de codimension $r$ de $A$ et $W$ un dessous d'escalier, on pose:
\begin{displaymath}
m_W(V)= r! \mathrm{max}_{ x \in V, 1 \leq i_1 < \cdots < i_d \leq g} \{ \mathrm{vol}(\mathbb{R}_+^{r} / \mathcal{C}_{i_1, \ldots, i_d}(W))\},
\end{displaymath}
où le maximum porte sur $x \in X$ et les $d$-faces de $\mathbb{N}^g$ telles que $(\partial_{i_1, x}, \ldots, \partial_{i_d, x})$ forment une base du quotient $t_{A, x}/t_{V,x}$. 

\bigskip

On peut maintenant énoncer le théorème dont on aura besoin:
\begin{thm}
Soit $V$ une sous-variété irréductible de $A$, de codimension $r$, $\tilde{M} \geq 1$ un entier et $\Sigma_0, \ldots, \Sigma_r$ des ensembles pondérés finis à support dans $A(\overline{K})$ tels que pour tout $1 \leq n \leq r$:
\begin{displaymath}
\mathrm{Supp} (\Sigma_n) = \bigcup_{l=1 \ldots s_n} H_{n,l},
\end{displaymath} 
où les $H_{n,l}$ sont des sous-groupes de $A(\overline{K})$; on suppose aussi que les dessous d'escaliers associés aux $\Sigma_n$ ne dépendent que de $n$. Soit de plus $f \in H^0(A, L^{\otimes \tilde{M}})$, non nulle, qui s'annule sur $V+ \Sigma_0 + \cdots + \Sigma_r$. 
Alors il existe une constante $c_{18}$ ne dépendant que de $A$, deux entiers $1 \leq r_0 \leq r_1 \leq r$, des indices $j_0, \ldots, j_{r_0-1}$ avec $1 \leq j_l \leq s_l$ pour $l=0 \ldots r_0-1$, et des sous-variétés algébriques  $Z_j \ (j=1, \ldots, s_{r_0})$ de $A$, propres, $\overline{K}$-irréductibles et de codimension $r_1$, contenant au moins une composante isolée de
\begin{displaymath}
 H_{0, j_0} + \cdots + H_{r_0-1, j_{r_0-1}} + \Sigma_{r_0} + \cdots + \Sigma_r+V, 
\end{displaymath}
telles que:
\begin{displaymath}
\mathrm{deg} \bigg( \bigcup_{x \in H_{0, j_0} + \cdots + H_{r_0-1, j_{r_0-1}}} \bigcup_{j=1 \ldots s_{r_0}, \  y \in H_{r_0,j}}  \{m_{W_{y}}(x+y+Z_{j})\} \big(x + y + Z_j\big) \bigg) \leq c_{18} \tilde{M}^{r_1}.
\end{displaymath}
\end{thm}
\begin{dem}
Le même théorème dans le cadre plus général des groupes algébriques est l'objet d'un article en préparation de David et Amoroso	(\cite{Amoroso-David07}).
\end{dem}

\subsection{Degré d'une sous-variété obstructrice}
\label{combinatoire}

Il s'agit d'adapter le lemme de zéros au cas qui nous intéresse; le degré $\tilde{M}$ du théorème qu'on vient d'énoncer sera égal à $M^2+1$ dans notre cas. 

On reprend les hypothèses et notations des parties \ref{choix fibres} et \ref{calcul fibres} et on rappelle que:
\begin{displaymath}
X=V + \Sigma_0,
\end{displaymath}
où $V$ est une variété irréductible et $\Sigma_0$ est un sous-groupe fini de $A$. On suppose enfin que le cardinal $|\Sigma_0|$ est premier à tous les premiers des $\mathcal{P}_{n, \mathbb{Z}}$, pour $1 \leq n \leq r$. Si $p$ est un nombre premier de $\bigcup_{1 \leq n \leq r} \mathcal{P}_{n, \mathbb{Z}}$ et $\mathfrak{q} / p$ dans $K'$, on désigne par: Ker$[p]_{\mathfrak{q}}$ le groupe des points de $p$-torsion sur $K'$ se réduisant sur $0$ modulo $\mathfrak{q}$.

\bigskip

Si $l= \prod_{n=1}^r p_n$ avec, pour tout $1 \leq n \leq r \ : \ p_n \in \mathcal{P}_{n, \mathbb{Z}}$ ou $p_n=1$, on note: 
\begin{displaymath}
\mathrm{Ker} [l]^*=  \bigoplus_ {n} \mathrm{Ker} [p_n]_{\mathfrak{q}_n}. 
\end{displaymath}
Cette somme est bien directe car le choix des paramètres implique que les $\mathcal{P}_{n, \mathbb{Z}}$ sont deux-à-deux disjoints.
Notre but, jusqu'à la fin de cette partie, sera de démontrer la proposition suivante, pour un bon choix du fermé $X$ (dont dépend la construction du morphisme):
\begin{prop}
\label{injectivite}

Le morphisme $\phi: E \rightarrow F$ est injectif. 
\end{prop}
On va supposer que ce n'est pas le cas et obtenir une contradiction en appliquant le lemme de zéros du paragraphe précédent. Celui-ci permet de majorer le degré d'une réunion de sous-variétés. On souhaite se ramener à une seule sous-variété obstructrice, et  utiliser le fait que la réunion est largement distincte. On y arrive par un travail sur le stabilisateur. Commençons donc par une définition:
\begin{defn}
\label{stabili}
Si $Z$ est une variété propre et irréductible de $A$, on appelle stabilisateur de $Z$, noté Stab$(Z)$, l'ensemble:
\begin{displaymath}
\{ x \in A, x+ Z= Z \} = \bigcap_{x \in Z} (Z - x).
\end{displaymath}
\end{defn}
On a les propriétés suivantes pour le stabilisateur:
\begin{displaymath}
\mathrm{dim}( \mathrm{Stab}(Z)) \leq \mathrm{dim} (Z) \ \ \mathrm{et:} \ \ \mathrm{deg}( \mathrm{Stab}(Z)) \leq \mathrm{deg}(Z)^{\mathrm{dim}(Z)+1}.
\end{displaymath}
La première suit de la définition, et la preuve torique de la seconde (dans \cite{Amoroso-David99}, 2) se transpose sans changement aux variétés abéliennes.   

\bigskip

On peut maintenant démontrer la proposition suivante: 
\begin{prop}
Il existe une constante $c_{20}$, des entiers $r_0 \leq r_1 \leq r$ strictement positifs, un entier $l \in \mathcal{P}_{1, \mathbb{Z}} \cdots \mathcal{P}_{r_0, \mathbb{Z}}$, et une sous-variété $Z$ de $A$, propre et $\overline{\mathbb{Q}}$-irréductible, de codimension $r_1$ contenant un translaté de $V$ par un point de torsion, tels que:
\begin{displaymath}
T_{r_0}^{r_1} |\mathcal{P}_{r, \mathbb{Z}}|\frac{\Sigma_0}{|\Sigma_0 \cap \mathrm{Stab}(Z)|}\frac{l^{2g- \alpha}}{|\mathrm{Ker}[l]^* \cap \mathrm{Stab} (Z)|} \mathrm{deg}(Z) \leq  c_{20} M^{2r_1} \Delta.
\end{displaymath}
\end{prop}
\textbf{Remarque} Pour simplifier les calculs qui viennent, on pose: 
\begin{displaymath}
f(\Sigma_0, Z)= \frac{\Sigma_0}{|\mathrm{Stab}(Z) \cap \Sigma_0|}.
\end{displaymath} 
\begin{dem}
Si le morphisme $\phi$ n'est pas injectif, il existe une section $f \in H^0(A, L^{\otimes M^2+1})$ qui s'annule sur:
\begin{displaymath}
\bigcup_{i \in I} V(X_i, t_A, T_{(i)}).
\end{displaymath}
Par définition des voisinages infinitésimaux, ceci implique que $f$ s'annule sur: 
\begin{displaymath}
X+\Sigma_1+ \cdots+ \Sigma_r, 
\end{displaymath}
où l'ensemble $\Sigma_n$, pour $1\leq n \leq r$, est pondéré de support:
\begin{displaymath}
\mathrm{Supp} (\Sigma_n)= \mathrm{Tor}_{A,n}= \bigcup_{\mathfrak{q}/ \mathfrak{p}  \in \mathcal{P}_{n}} \mathrm{Ker}[p]_{\mathfrak{q}};
\end{displaymath} 
cet ensemble est associé au dessous d'escalier simple et ne dépendant que de $n$ défini par: 
\begin{displaymath}
\sum_{k=1}^{g} \lambda_k \leq T_n.
\end{displaymath}
Il existe donc, par le théorème précédent, deux entiers $r_0$ et $r_1$ tels que: $r_0 \leq r_1 \leq r$, des couples d'idéaux premiers $(p_1, \mathfrak{q}_1), \ldots, (p_{r_0-1}, \mathfrak{q}_{r_0-1})$ avec $\mathfrak{q}_n \in \mathcal{P}_n$ (pour $1 \leq n \leq r_{0}-1$) et des sous-variétés algébriques $Z_{\mathfrak{q}}$ de $A$ (pour tout $\mathfrak{q} \in \mathcal{P}_{r_0}$), propres et $\overline{\mathbb{Q}}$-irréductibles, de codimension $r_1$, tels que:
\begin{displaymath}
\mathrm{deg} \bigg( \bigcup_{\mathfrak{q} \in \mathcal{P}_{r_0}} \bigcup_{\zeta \in \Sigma_0 \oplus \bigoplus_n \mathrm{Ker}[p_n]_{\mathfrak{q}_n}}   \{m_{W_{\zeta_{r_0}}}(\zeta + Z_{\mathfrak{q}})\} \big(\zeta+ Z_{\mathfrak{q}}\big)  \bigg) \leq c_{18}M^{2r_1}, 
\end{displaymath}
où on a écrit, pour unifier l'écriture dans la somme directe: $\mathfrak{q}_{r_0}= \mathfrak{q}$; et où $\zeta_{r_0}$ est la composante selon $r_0$ de $\zeta$ dans la somme directe.
De plus, pour tout $\mathfrak{q} \in \mathcal{P}_{r_0}$, la variété $Z_{\mathfrak{q}}$ contient un translaté de $V$ par un point de torsion. 

\bigskip

Le terme de multiplicité se calcule immédiatement. Soit $\mathfrak{q} \in \mathcal{P}_{r_0}$ et $\zeta \in \Sigma_0 \oplus \bigoplus_n \mathrm{Ker}[p_n]_{\mathfrak{q}_n}$; on a:
\begin{displaymath}
m_{W_{\zeta_{r_0}}}(\zeta + Z_{\mathfrak{q}}) = r_1!T_{r_0}^{r_1} \mathrm{vol}(\{ u_1+ \cdots + u_{r_1}<1 \}) = T_{r_0}^{r_1}.
\end{displaymath}
Cette multiplicité ne dépend pas de la variété dans la réunion donc on peut la mettre en facteur.

\bigskip

Les entiers $r_0$ et $r_1$ sont déjà déterminés; posons $l_0= p_1 \cdots p_{r_0-1}$. Choisissons, pour tout premier $p \in \mathcal{P}_{r_0, \mathbb{Z}}$, un idéal premier $\mathfrak{q}$ de $\mathcal{P}_{r_0}$ divisant $p$ tel que la quantité:
\begin{displaymath}
f(\Sigma_0, Z_{\mathfrak{q}}) \frac{(l_0 p)^{2g- \alpha}}{|\mathrm{Ker}[l_0p]^* \cap \mathrm{Stab} (Z_{\mathfrak{q}})|} \mathrm{deg} ( Z_{\mathfrak{q}})
\end{displaymath}
soit minimale parmi les premiers de $\mathcal{P}_{r_0}$divisant $p$.  Prenons aussi $\mathfrak{q}_{r_0} \in \mathcal{P}_{r_0}$ (et $p_{r_0}$) tels que cette même quantité soit minimale parmi tous les premiers de $\mathcal{P}_{r_0}$. On pose $l=l_0 p_{r_0}$ et $Z=Z_{\mathfrak{q}_{r_0}}$. Il suffit donc de majorer cette quantité pour obtenir la proposition. Rappelons que la somme:
\begin{displaymath}
\mathrm{Ker} [l_0]^* = \bigoplus_{n=1}^{r_0-1} \mathrm{Ker}[p_n]_{\mathfrak{q}_n}
\end{displaymath}
est bien directe car, les premiers $p_n$ étant deux-à-deux distincts, on peut écrire une relation de Bézout entre un des $p_n$ et tous les autres. 
On partitionne $\mathcal{P}_{r_0, \mathbb{Z}}$ en introduisant la relation d'équivalence suivante:
\begin{displaymath}
p \sim p' \Longleftrightarrow \bigg( \exists \ \gamma \in  \Sigma_0 \oplus \mathrm{Ker} [l_0]^* \oplus \bigoplus_{p_i \in \mathcal{P}_{r_0, \mathbb{Z}}} \Big( \mathrm{Ker}[p_i]_{\mathfrak{q}_i} \Big), \mathrm{ \  tel \ que \ } \gamma + Z_{\mathfrak{q}} = Z_{\mathfrak{q}'} \bigg),
\end{displaymath}
et on note $(\mathcal{C}_1, \ldots, \mathcal{C}_s)$ les différentes classes d'équivalence associées.
Les variétés $Z_{\mathfrak{q}}$ sont irréductibles et il en va de même pour chacune de leurs translatées.
Si $p$ et $p'$ appartiennent à des classes différentes, les réunions
\begin{displaymath}
\bigcup_{\zeta \in \Sigma_0 \oplus \mathrm{Ker}[l_0]^* \oplus \mathrm{Ker}[p]_{\mathfrak{q}}} \zeta + Z_{\mathfrak{q}} 
\end{displaymath}
n'ont aucune composante en commun et on peut additionner les degrés. Le choix d'un seul idéal de $\mathcal{P}_{r_0}$ au-dessus d'un nombre premier restreint la réunion. On a donc:
\begin{eqnarray}
T_{(r_0)}^{r_1} \ \sum_{j=1}^{s} \mathrm{deg} \bigg( \bigcup_{p \in \mathcal{C}_j} \bigcup_{\zeta \in \Sigma_0 \oplus \mathrm{Ker}[l_0]^{*} \oplus \mathrm{Ker}[p]_{\mathfrak{q}}} \zeta + Z_{\mathfrak{q}}  \bigg) & \leq & c_{18}M^{2r_1}. 
\end{eqnarray}

\bigskip

Soit $p \in \mathcal{P}_{r_0, \mathbb{Z}}$ et $\mathfrak{q}$ l'idéal qui lui est associé; le stabilisateur de $Z_{\mathfrak{q}}$ ne dépend que de la classe d'équivalence de $p$ puisque si $p'$ (associé à $\mathfrak{q}'$) est dans la même classe que $p$, $Z_{\mathfrak{q}'}$ est un translaté de $Z_{\mathfrak{q}}$. On appelle $\mathcal{S}_j$ le stabilisateur commun aux $Z_{\mathfrak{q}}$, pour $\mathfrak{q}$ associé à $p \in \mathcal{C}_j$. Dans chaque classe $\mathcal{C}_j$, on fixe un premier $\rho_j \in \mathcal{C}_j$ et on note $Z_{\rho_j}$ la variété qui lui est associée. Pour tout autre premier $p \in \mathcal{C}_j$, il existe donc un élément
$\alpha_{p} \in \bigoplus_{p_i \in \mathcal{P}_{r_0, \mathbb{Z}}} \mathrm{Ker}[p_i]_{\mathfrak{q}_i}$ et $\eta_{p} \in \Sigma_0 \oplus \mathrm{Ker}[l_0]^*$ tels que:
\begin{displaymath}
\alpha_{p} + \eta_{p} + Z_{\mathfrak{q}}= Z_{\rho_j}.
\end{displaymath}
Remarquons que la somme est directe car tous les $p_i$ sont distincts. Soient $p \neq p'$ dans la même classe $\mathcal{C}_j$; soient $\omega_{p} \in \mathrm{Ker}[p]_{\mathfrak{q}}$ et $\omega_{p'} \in \mathrm{Ker}[p']_{\mathfrak{q}'}$. Si les réunions
\begin{eqnarray}
\bigcup_{\zeta \in \Sigma_0 \oplus \mathrm{Ker}[l_0]^{*}}  \zeta + \omega_{\xi} + Z_{\xi}  
\end{eqnarray}
pour $\xi=p$ et $\xi = p'$, ont au moins une composante commune, c'est qu'il existe un élément $\eta_{p,p'} \in \Sigma_0 \oplus \mathrm{Ker} [l_0]^*$ tel que:
\begin{displaymath}
\omega_{p} + Z_{p}= \eta_{p,p'} + \omega_{p'} + Z_{p'}.
\end{displaymath}
On en déduit, grâce aux deux dernières égalités, que:
\begin{displaymath}
x= \alpha_{p} - \omega_{p} - \alpha_{p'}  + \omega_{p'} + \big( \eta_{p} - \eta_{p'}  + \eta_{p,p'} \big) \in \mathcal{S}_j.
\end{displaymath}
On note $\alpha_{p}^{p_i}$ la composante selon $p_i$ de $\alpha_{p}$. 
On remarque que: $\alpha_{p}^{p} - \alpha_{p'}^{p} - \omega_{p} \in \mathrm{Ker}[p]_{\mathfrak{q}}$. De même: $\alpha_{p'}^{p'} - \alpha_{p}^{p'} - \omega_{p'} \in \mathrm{Ker}[p']_{\mathfrak{q}'}$, et: $\eta_{p'} - \eta_{p}  + \eta_{p,p'} \in \Sigma_0 \oplus \mathrm{Ker}[l_0]^*$. Le nombre $p$ est premier à $p'$, à $l_0$, à $\Sigma_0$ et à tous les autres premiers de $\mathcal{P}_{r_0, \mathbb{Z}}$. Il existe donc une relation de Bézout: 
\begin{displaymath}
up + v l_0\prod_{p_i \neq p \in \mathcal{P}_{r_0, \mathbb{Z}}} p_i=1.
\end{displaymath}
On en déduit que:
\begin{displaymath}
[v l_0 \prod_{p_i \neq p \in \mathcal{P}_{r_0, \mathbb{Z}}} p_i]x=[v l_0 \prod_{p_i \neq p \in \mathcal{P}_{r_0, \mathbb{Z}}} p_i](\alpha_{p}^{p} - \alpha_{p'}^{p} - \omega_{p})= \alpha_{p}^{p} - \alpha_{p'}^{p} - \omega_{p} \in \mathcal{S}_j;
\end{displaymath}
puisqu'il suit de sa définition que le stabilisateur est stable sous la multiplication par $n$, quel que soit $n \in \mathbb{N}$.
Par contraposition, si:
\begin{center}
$\omega_{p} \in \mathrm{Ker}[p]_{\mathfrak{q}} \setminus \big( \alpha_{p}^{p} - \alpha_{p'}^{p} - \mathcal{S}_j \big)$ et $\omega_{p'} \in \mathrm{Ker}[p']_{\mathfrak{q}'} \setminus \big( \alpha_{p'}^{p'} - \alpha_{p}^{p'} + \mathcal{S}_j \big)$, 
\end{center}
les réunions (5) n'ont pas de composantes communes. Il suit:
\begin{displaymath}
\mathrm{deg} \bigg( \bigcup_{p \in \mathcal{C}_j} \bigcup_{\zeta \in \Sigma_0 \oplus \mathrm{Ker}[l_0]^*} \bigcup_{\xi \in \mathrm{Ker}[p]_{\mathfrak{q}}} \zeta + \xi + Z_{\mathfrak{q}} \bigg) \geq \sum_{p \in \mathcal{C}_j} \mathrm{deg} \bigg( \bigcup_{\zeta \in \Sigma_0 \oplus \mathrm{Ker}[l_0]^*} \bigcup_{\xi \in \mathrm{Ker}[p]_{\mathfrak{q}} \setminus \bigcup_i \big(\alpha_{p}^{p} - \alpha_{p_i}^{p} + \mathcal{S}_j\big)} \zeta + \xi + Z_{\mathfrak{q}} \bigg).
\end{displaymath}

\bigskip

Fixons $j$ et $p \in \mathcal{C}_j$. On va calculer le degré de la réunion totale en fonction de deg$(Z_{\mathfrak{q}})$. Par choix de l'ensemble $\mathcal{P}_{A}$, il y a $p^{2g- \alpha}$ points se réduisant sur $0$ mod $\mathfrak{q}$, et il y a $l_0^{2g-\alpha}$ points dans $\mathrm{Ker}[l_0]^*$. Il en résulte:
\begin{eqnarray}
\mathrm{deg} \bigg( \bigcup_{ \zeta \in \Sigma_0 \oplus \mathrm{Ker}[l_0]^*} \bigcup_{ \xi \in \mathrm{Ker}[p]_{\mathfrak{q}}} \zeta + \xi + Z_{\mathfrak{q}} \bigg) =  \frac{f(\Sigma_0, Z_{\mathfrak{q}}) (l_0p)^{2g-\alpha}}{|\mathcal{S}_j \cap \big(\mathrm{Ker} [l_0p]^* \big)|} \mathrm{deg} (Z_{\mathfrak{q}}). 
\end{eqnarray}
A cette réunion, il faut retrancher:
\begin{displaymath}
\mathrm{deg} \bigg( \bigcup_{ \zeta \in \Sigma_0 \oplus \mathrm{Ker}[l_0]^*} \bigcup_{\xi \in \bigcup_i \big( \alpha_{p}^{p} - \alpha_{p_i}^{p} + \mathcal{S}_j \big)} \zeta + \xi + Z_{\mathfrak{q}} \bigg) \leq \frac{|\mathcal{C}_j|l_{0}^{2g-\alpha}}{|\mathcal{S}_j \cap \mathrm{Ker} [l_0]^*|} \mathrm{deg} (Z_{\mathfrak{q}}).
\end{displaymath}
En effet, il y a au plus $|\mathcal{C}_j|$ points de la forme $\alpha_{p_i}^{p}$. Notons $\tilde{\mathcal{C}}_j$ le sous-ensemble de $\mathcal{C}_j$ formé des $p$ divisant $[\mathcal{S}_j  :\mathcal{S}_j^0]$, où $\mathcal{S}_j^0$ désigne la composante connexe de l'identité dans $\mathcal{S}_j$. Si $p \notin \tilde{\mathcal{C}}_j$, les dénominateurs des deux dernières formules sont égaux et comme $\alpha \leq g$, on a:
\begin{eqnarray*}
\mathrm{deg} \bigg( \bigcup_{ \zeta \in \Sigma_0 \oplus \mathrm{Ker}[l_0]^*} \bigcup_{\xi \in \mathrm{Ker} [p]_{\mathfrak{q}} \setminus \bigcup_i \big(\alpha_{p}^{p} - \alpha_{p_i}^{p} + \mathcal{S}_j\big)} \zeta + \xi + Z_{\mathfrak{q}} \bigg) & \geq & \bigg(1- \frac{|\mathcal{C}_j|}{p^2} \bigg) \frac{f(\Sigma_0, Z_{\mathfrak{q}}) (l_0p)^{2g-\alpha}}{|\mathcal{S}_j \cap \big(\mathrm{Ker} [l_0p]^* \big)|} \mathrm{deg} (Z_{\mathfrak{q}}).\\
\end{eqnarray*}
Le quotient $\frac{1}{p^2}$ provient du fait que la ``partie discrète'' du stabilisateur de $Z_{\mathfrak{q}}$ est triviale et que sa composante connexe en $0$ est un groupe algébrique de codimension $\geq r+1 \geq 2$.
En fixant $j$, on somme sur l'ensemble des $p$; en tenant compte de la définition de $Z$, on obtient:
\begin{eqnarray*}
\mathrm{deg} \bigg( \bigcup_{p \in \mathcal{C}_j} \bigcup_{ \zeta \in \Sigma_0 \oplus \mathrm{Ker}[l_0]^{*}} \bigcup_{\xi \in \mathrm{Ker} [p]_{\mathfrak{q}}} \zeta + Z_{\mathfrak{q}} \bigg) & \geq & \ \bigg( |\mathcal{C}_j \setminus \tilde{\mathcal{C}}_j| - |\mathcal{C}_j| \sum_{p \in \mathcal{C}_j} \frac{1}{p^2} \bigg) \  \frac{f(\Sigma_0, Z) (l)^{2g- \alpha}}{|\mathrm{Ker}[l]^* \cap \mathrm{Stab} (Z)|} \mathrm{deg} (Z) \\
 & \geq & \Big( \frac{2}{3}|\mathcal{C}_j|-|\tilde{\mathcal{C}_j}| \Big)\frac{f(\Sigma_0, Z) (l)^{2g- \alpha}}{|\mathrm{Ker}[l]^* \cap \mathrm{Stab} (Z)|} \mathrm{deg} (Z) .
\end{eqnarray*}
On avait en effet choisi l'ensemble $\mathcal{P}_{A}$ tel que la somme $\sum_{\mathfrak{p}/p \in \mathcal{P}_{A}}\frac{1}{p^2} < \frac{1}{3}$.
On a plus directement, par un calcul direct à partir de (6):
\begin{displaymath}
\mathrm{deg} \bigg( \bigcup_{p \in \mathcal{C}_j} \bigcup_{ \zeta \in \mathrm{Ker}[l_0]^{*}} \bigcup_{\xi \in \mathrm{Ker} [p]_{\mathfrak{q}}} \zeta + Z_{\mathfrak{q}} \bigg) \geq    \frac{f(\Sigma_0, Z) (l)^{2g- \alpha}}{|\mathrm{Ker}[l]^* \cap \mathrm{Stab} (Z)|} \mathrm{deg} (Z).
\end{displaymath}
On doit donc estimer le nombre de premiers divisant $[\mathcal{S}_j  :\mathcal{S}_j^0]$. Or:
\begin{displaymath}
|\tilde{\mathcal{C}}_j| \leq \frac{\mathrm{log}[\mathcal{S}_j:\mathcal{S}_j^0]}{\mathrm{log}3} \leq \mathrm{log \ deg} (\mathcal{S}_j) \leq c_{19} \Delta,
\end{displaymath}
pour une constante $c_{19}$. On a ici majoré le degré du stabilisateur en fonction de celui de la variété (\cf \cite{Hindry88}, lemme 6), puis on a utilisé le lemme de zéros pour majorer deg$(Z_{\rho_j})$, et on a majoré log$(M)$ à l'aide du choix des paramètres. Par l'inégalité: max$\{x-y; 1\} \geq \frac{x}{2y}$ pour $x \geq 0$ et $y \geq 1$, on obtient: 
\begin{displaymath}
\mathrm{max} \{ \frac{2}{3}|\mathcal{C}_j|-|\tilde{\mathcal{C}_j}|, 1 \} \geq \frac{|\mathcal{C}_j|}{3 c_{19} \Delta}.
\end{displaymath}
La proposition suit en sommant sur les classes d'équivalence. 
\end{dem}

\subsection{Un premier pas vers l'injectivité}
\label{quasi contradiction}

On suppose maintenant que $A$ vérifie l'hypothèse \textbf{H2}. On vérifie aisément que celle-ci implique \textbf{H3} (voir \ref{densite}).

\bigskip

Le choix des paramètres va nous donner une inégalité ``presque absurde''; on ne pourra cependant pas conclure, car il manquera une hypothèse de coprimalité sur des objets construits simultanément. On devra donc itérer une fois le résultat obtenu, en exigeant cette hypothèse entre la première et la seconde étape, puis on conclura par un argument de descente.  

\bigskip

Rappelons que la dernière proposition nous a donné l'existence d'un couple $(l,Z)$ où $l$ est un entier assez grand (on sait que $l \in \mathcal{P}_1 \cdots \mathcal{P}_{r_{0}}$) et $Z$ est une sous-variété irréductible contenant un translaté de $V$. La proposition suivante résume et précise ce qu'on a obtenu.
\begin{prop}
\label{quasi injectivite}
On suppose que $X$ n'est pas incluse dans le translaté d'une sous-variété abélienne et que son minimum essentiel est majoré de la façon suivante:
\begin{displaymath}
\hat{\mu}_{\mathrm{ess}}(X) \ \omega(X) < \frac{1}{\Delta^{8 \rho (2r)^{r+1}}}.
\end{displaymath}
Alors il existe une sous-variété propre $Z$ de codimension $r_1 \leq r$ contenant un translaté de $V$ par un point de torsion et un entier $l>0$ tels que: 

-L'entier $l$ est premier avec $R$ et:
\begin{displaymath}
l \leq \Delta^{2 \rho(2r)^{r+1}}.
\end{displaymath}

-De plus, on a l'inégalité: 
\begin{displaymath}
\Big(\frac{f(\Sigma_0, Z)l^{gn_{\alpha}}}{|\mathrm{Ker}[l]^* \cap \mathrm{Stab} (Z)|}\mathrm{deg}(Z)\Big)^{1/r_1} < \Delta^{-\rho} l^{n_{\alpha}} \omega(l^{n_{\alpha}},X). 
\end{displaymath}
\end{prop}
\begin{dem}
On commence par montrer que la contrainte portant ici sur le minimum essentiel est plus forte que l'hypothèse (\ref{minimum essentiel}). Si (\ref{minimum essentiel}) n'est pas vérifiée: 
\begin{displaymath}
\hat{\mu}_{\mathrm{ess}}(X) \geq \frac{\Delta}{T_0 \omega(\Delta T_0, X)} \geq \frac{1}{T_0^2 \omega(X)},
\end{displaymath}
par application du lemme \ref{obstruction avec poids}. De plus, par les choix de paramètres faits en \ref{parametres}:
\begin{displaymath}
T_0^2 \leq \Delta^{2r}\prod_{n=1}^r N_n^4 \leq \Delta^{8 \rho (2r)^{r+1}},
\end{displaymath}
ce qui contredit l'hypothèse de la proposition. Pour démontrer cette inégalité, on a d'abord utilisé:
\begin{eqnarray}
\label{grosse majoration}
N_1 \cdots N_{r} \leq \Delta^{ \rho[(2r)^2+ \cdots + (2r)^{r+1}]};
\end{eqnarray}
puis on a majoré l'exposant comme suit:
\begin{displaymath}
\sum_{j=2}^{r+1}(2r)^j \leq -r + \sum_{j=1}^{r+1}(2r)^j \leq 2(2r)^{r+1} -r,
\end{displaymath}
par l'in\'egalit\'e:
\begin{displaymath}
1+x+ \ldots +x^h \leq 2x^h \mathrm{ \ pour \ } h \in \mathbb{N} \mathrm{ \ et \ } x \geq 2.
\end{displaymath}

\bigskip

La proposition précédente nous donne donc l'existence de trois entiers strictement positifs $r_0$, $r_1$ et $l$ avec $r_0 \leq r_1 \leq r$, $l \in \mathcal{P}_{1, \mathbb{Z}} \cdots \mathcal{P}_{r_0, \mathbb{Z}}$, et une sous-variété algébrique $Z$  propre et irréductible de $A$, de codimension $r_1$, contenant un translaté de $V$ par un point de torsion, telle que:

\begin{displaymath}
T_{r_0}^{r_1} |\mathcal{P}_{r_0, \mathbb{Z}}|\frac{f(\Sigma_0, Z)l^{2g- \alpha}}{|\mathrm{Ker}[l]^* \cap \mathrm{Stab} (Z)|} \mathrm{deg}(Z) \leq c_{20} M^{2r_1} \Delta.
\end{displaymath}
Par construction des $\mathcal{P}_{n, \mathbb{Z}}$, l'entier $l$ est premier avec $R$ et on a les inégalités:
\begin{displaymath}
2^{-r_0}N_1 \cdots N_{r_0} \leq l \leq N_1 \cdots N_{r_0}.
\end{displaymath}
Et le premier point suit, par la même majoration que (\ref{grosse majoration}).

\bigskip

Reste \`a prouver la seconde in\'egalit\'e.
Le th\'eor\`eme des nombres premiers et le choix de l'ensemble $\mathcal{P}_{A}$ font que, pour une certaine constante $c_{21}>0$ ne dépendant que de $A$:
\begin{displaymath}
|\mathcal{P}_{r_0, \mathbb{Z}}| \geq c_{21} \frac{N_{r_0}}{\mathrm{log} N_{r_0}} - \frac{\mathrm{log} R}{\mathrm{log} 2}.
\end{displaymath}
Par définition des $N_n$ et de $\Delta$, on a: $\mathrm{log} N_{r_0} \leq \Delta^{1/2}$ pour $C_0$ assez grand dans la définition de $\Delta$. On a aussi:
\begin{displaymath}
N_{r_0} \geq \Delta^9 \geq \mathrm{log}(R)^2.  
\end{displaymath}
On a encore, pour $C_0$ assez grand: $\frac{1}{2} c_{21} \Delta^{1/2} \geq 1$, le facteur $\frac{1}{2}$ correspondant au terme en log$(R)$. On en déduit: 
\begin{displaymath}
|\mathcal{P}_{r_0, \mathbb{Z}}| \geq \Delta^{\rho (2r)^{r+2-r_0}-1}.
\end{displaymath}

Par le lemme 2.4 de \cite{Amoroso-David03}, comme $l^{n_{\alpha}} \leq (N_1 \cdots N_r)^{n_{\alpha}} \leq \Delta T_0$, on a:
\begin{displaymath}
\omega(\Delta T_0,X) \leq \frac{\Delta T_0}{l^{n_{\alpha}}} \omega(l^{n_{\alpha}},X).
\end{displaymath}
Sous l'hypothèse \textbf{H2}, on peut supposer que $\alpha=0$ ou $\alpha=g$. Dans les deux cas, on a:
\begin{eqnarray*}
\Big(\frac{f(\Sigma_0, Z)l^{n_{\alpha}g}}{|\mathrm{Ker}[l]^* \cap \mathrm{Stab} (Z)|} \mathrm{deg}(Z)\Big)^{1/r_1} & \leq & \frac{c_{20} M^2  \Delta^{1/r_1}}{ |\mathcal{P}_{r_0, \mathbb{Z}}|^{1/r_1}T_{r_0}} \\
 & \leq & c_{20} \frac{ T_0 \omega(\Delta T_0,X) \Delta^{1/r_1}}{ |\mathcal{P}_{r_0, \mathbb{Z}}|^{1/r_1}T_{r_0}} \\
 & \leq & c_{20} \frac{ \Delta T_0^2  \Delta^{1/r_1}}{ l^{n_{\alpha}}|\mathcal{P}_{r_0, \mathbb{Z}}|^{1/r_1}T_{r_0}} \omega(l^{n_{\alpha}},X). \\
\end{eqnarray*} 
Or on a:
\begin{displaymath}
\frac{T_0}{T_{r_0}} \leq \Delta^{2r_0}\big(N_1 \cdots N_{r_0}\big)^{n_{\alpha}} \leq (2\Delta)^{2r_0} l^{n_{\alpha}},
\end{displaymath}
et on en déduit:
\begin{displaymath}
\Big(\frac{f(\Sigma_0, Z)l^{n_{\alpha}g}}{|\mathrm{Ker}[l]^* \cap \mathrm{Stab} (Z)|} \mathrm{deg}(Z)\Big)^{1/r_1} \leq c_{22} \frac{ l^{n_{\alpha}} \Delta^{2 r + 1 + 1/r_1} N_{r_0+1}^{n_{\alpha}} \cdots N_r^{n_{\alpha}}}{ |\mathcal{P}_{r_0, \mathbb{Z}}|^{1/r_1}} \omega(l^{n_{\alpha}},X). 
\end{displaymath}
L'exposant $h$ de $\Delta$ dans cette dernière majoration est borné par:
\begin{eqnarray*}
h & := & 4r + 2\rho \big((2r)^2+ \cdots + (2r)^{r+1-r_0} \big)- \big(\rho(2r)^{r+2-r_0}-2 \big)/r_1. \\
& \leq & 2\rho\big( (2r) + \cdots + (2r)^{r+1-r_0}\big)-2\rho(2r)^{r+1-r_0} \\
& \leq & 2\rho(r-r_0)(2r)^{r-r_0}+\rho(2r)^{r+1-r_0}-2\rho(2r)^{r+1-r_0} \\
& \leq & - 2\rho r_0(2r)^{r-r_0} \leq -2 \rho. \\
\end{eqnarray*}
On a donc finalement: 
\begin{displaymath}
\Big(\frac{f(\Sigma_0, Z)l^{n_{\alpha}g}}{|\mathrm{Ker}[l]^* \cap \mathrm{Stab} (Z)|} \mathrm{deg}(Z)\Big)^{1/r_1} < \Delta^{-\rho} l^{n_{\alpha}} \omega(l^{n_{\alpha}}, X), 
\end{displaymath}
en faisant disparaître les constantes avec $\Delta^{\rho}$, et le résultat suit. 
\end{dem}
\textbf{Remarques} On notera dorénavant:
\begin{displaymath}
|\mathrm{Ker}[l]^* \cap \mathrm{Stab}(Z)|= \lambda(Z).
\end{displaymath}
Posons $X=V$ et $\Sigma_0=\{0\}$. Si on savait assurer la coprimalité entre $l$ et $[\mathrm{Stab}(Z): \mathrm{Stab}(Z)^0]$ (deux objets construits simultanément), on pourrait déjà clore la preuve, car on aurait: 
\begin{displaymath}
\lambda(Z) \leq l^{n_{\alpha}\mathrm{dim} \big( \mathrm{Stab} (Z)^0 \big)} \leq l^{n_{\alpha} (g-r_1-1)},
\end{displaymath}
la deuxième inégalité provenant du fait que $V$ n'est pas inclus dans un translaté de sous-variété abélienne. La variété $Z$ contenant un translaté de $V$ par un point de torsion, on a de plus:
\begin{displaymath}
\omega(l^{n_{\alpha}}, V) \leq \Big( l^{n_{\alpha}} \mathrm{deg} (Z) \Big)^{1/r_1},
\end{displaymath}
et une contradiction suivrait immédiatement.

\subsection{It\'eration et descente}
\label{descente 1}

Cette construction ne  permet pas de conclure, et on est amen\'e \`a it\'erer la derni\`ere proposition. C'est à ce stade de la preuve qu'on utilise crucialement l'hypothèse sur la codimension de $V$. En effet, on n'a pas pu mettre en place la stratégie de descente, devenue classique dans les travaux diophantiens sur la minoration de hauteurs, en codimension quelconque.

\bigskip

Cette démarche échoue en grande partie pour la raison suivante: la procédure diophantienne donne l'existence d'une sous-variété obstructrice vérifiant une propriété ``pathologique'' mais on perd trop d'information en extrayant une hypersurface contenant cette variété obstructrice. Le passage par l'hypersurface est pourtant indispensable pour garantir l'emboîtement, et par suite l'égalité des dimensions après $r$ itérations.\\
\\
\textbf{Preuve} (du théorème \ref{theoreme}) \\

On peut maintenant démontrer le théorème \ref{theoreme}, qui découlera de la non-injectivité d'un morphisme $\phi$ par le résultat du paragraphe \ref{pentes consequences}. Soit $V$ une sous-variété propre de $A$ qui n'est pas incluse dans un translaté de sous-variété abélienne de $A$, de codimension $r \leq 2$. On rappelle que:
\begin{displaymath}
\Delta= C_0^2 \mathrm{log}(3 \mathrm{deg}(V)),
\end{displaymath}
et on suppose: 
\begin{eqnarray}
\label{1102}
\omega(V) \hat{\mu}_{\mathrm{ess}}(V) < \Delta^{-\big( 16(2r)^{r+1}\big)^r}.
\end{eqnarray}
\\
{\it Première étape.} Pour utiliser la proposition \ref{quasi injectivite}, on doit définir: 
\begin{center}
$\rho_1=  \big(9(2r)^{r+1}\big)^{r-1}$ et $R_1=[\mathrm{Stab}(V): \mathrm{Stab}(V)^0]$.
\end{center}
On a, en tenant compte des propriétés du stabilisateur suivant la définition \ref{stabili}:
\begin{displaymath}
\mathrm{log}(R_1) \leq \mathrm{log}\Big(\mathrm{deg} (\mathrm{Stab} (V))\Big) \leq g \mathrm{log}(3\mathrm{deg}(V)) \leq \Delta.
\end{displaymath}
Si le morphisme $\phi$ n'est pas injectif, on applique la proposition \ref{quasi injectivite} avec $X=V$, ce qui donne l'existence d'un entier $l_1$ et d'une sous-variété $Z_1$ de $A$, propre et de codimension $k_1$, contenant un translaté de $V$ par un point $x_1$, et telle que:
\begin{displaymath}
\left(\frac{l_1^{n_{\alpha_1}g}\mathrm{deg}(Z_1)}{\lambda_1(Z_1)}\right)^{1/k_1} < \Delta^{-\rho_1} l_1^{n_{\alpha_1}} \omega(l_1^{n_{\alpha_1}}, V). 
\end{displaymath}
De plus, on peut supposer que $V$ est de codimension $2$ et que $Z_1$ est une hypersurface. Sinon, on aurait: $Z_1=x_1+V$, car ces variétés seraient de même codimension et: $x_1+V \subset Z_1$. Dans ce cas, l'entier $l_1$ serait premier à
\begin{displaymath}
[\mathrm{Stab}(V): \mathrm{Stab}(V)^0]=[\mathrm{Stab}(Z_1): \mathrm{Stab}(Z_1)^0],
\end{displaymath} 
et la remarque suivant la preuve de la proposition \ref{quasi injectivite} montre qu'on aurait une contradiction. \\
\\
{\it Deuxième étape.} On itère maintenant la proposition \ref{quasi injectivite} en posant:
\begin{displaymath}
V_1= \bigcup_{x \in \mathrm{Stab} (Z_1) \cap \mathrm{Ker}[l_1]^*} x+V.
\end{displaymath}
Puis:
\begin{center}
$\rho_2= \big(9(2r)^{r+1}\big)^{r-2}$ et $R_2=  [\mathrm{Stab}(V): \mathrm{Stab}(V)^0] \times [\mathrm{Stab}(Z_1): \mathrm{Stab}(Z_1)^0] \times l_1$.
\end{center}
La dernière condition permet que le cardinal de $\Sigma_0$ soit premier à tous les premiers des $\mathcal{P}_{i, \mathbb{Z}}$ dans la phase combinatoire. On vérifie une nouvelle fois (par les majorations du degré de $Z_1$ et de $l_1$ données par la proposition \ref{quasi injectivite}) que:
\begin{displaymath}
\mathrm{log}(R_2) \leq g \mathrm{log}(3 \omega(V)) + g^2 \mathrm{log} (\omega(V)) + 3 \mathrm{log}(l_1) \leq \Delta. 
\end{displaymath}
L'hypothèse (\ref{sigma zero}) est satisfaite pour les mêmes raisons (on a la majoration: $|\Sigma_0| \leq l_1^{2g}$).
On doit aussi majorer:
\begin{displaymath}
\omega(V_1) \hat{\mu}_{\mathrm{ess}}(V_1).
\end{displaymath}
Le minimum essentiel de $V_1$ est celui de $V$ puisqu'on translate par des points de torsion. Quant à l'indice d'obstruction, comme $x_1+ V_1 \subset Z_1$ (par définition de ces deux variétés), l'inégalité sur le degré de $Z_1$ donne:
\begin{eqnarray}
\label{10000}
\omega(V_1) \leq l_1^2 \omega(l_1^2, V) \leq l_1^4 \omega(V).
\end{eqnarray}
On obtient donc:
\begin{displaymath}
\omega(V_1) \hat{\mu}_{\mathrm{ess}}(V_1) \leq l_1^4 \omega(V) \hat{\mu}_{\mathrm{ess}}(V) \leq \Delta^{-\big( 16(2r)^{r+1}\big)^r+8 \rho_1(2r)^{r+1}} \leq \Delta^{\rho_2(2r)^{(r+1)}(-16+8)} \leq \Delta^{-8 \rho_2 (2r)^{r+1}}.
\end{displaymath}
Par (\ref{10000}), on a enfin:
\begin{displaymath}
C_0 \mathrm{log}\big(3 \omega(V_1)\big) \leq \Delta. 
\end{displaymath}
La proposition \ref{quasi injectivite} avec $\Sigma_0= \mathrm{Stab} (Z_1) \cap \mathrm{Ker}[l_1]^*$ donne l'existence d'une variété $Z_2$ de codimension $k_2$ contenant un translaté $x_2+V$, telle que:
\begin{displaymath}
\left(\frac{f(\Sigma_0, Z_2)l_2^{n_{\alpha_2}g}\mathrm{deg}(Z_2)}{\lambda_2(Z_2)}\right)^{1/k_2} < \Delta^{-\rho_2} l_2^{n_{\alpha_2}} \omega(l_2^{n_{\alpha_2}}, V_1). 
\end{displaymath}
On remarque que $Z_2$ contient les translatés de $x_2+V$ par les points de $H_0 \cap \mathrm{Stab} (Z_2)$. On a:
\begin{displaymath}
\mathrm{deg}\Big( \bigcup_{x \in \Sigma_0/(\Sigma_0 \cap \mathrm{Stab} Z_2)} x + Z_2 \Big) \leq f(\Sigma_0, Z_2) \mathrm{deg} (Z_2);
\end{displaymath}
et cette réunion, notée $Z_2'$, contient un translaté de $V_1$. Si $Z_2$ est de codimension $2$, on a encore une égalité: $Z_2=x_2+V$ et on en déduit que $l_2$ est premier à $[\mathrm{Stab}(Z_2): \mathrm{Stab}(Z_2)^0]$. Il suit:
\begin{displaymath}
\Big(l_2^{n_{\alpha_2}} \mathrm{deg}(Z_2')\Big)^{1/2} \leq \Delta^{-\rho_2} \omega(l_2^{n_{\alpha_2}},V_1), 
\end{displaymath}
ce qui est absurde, puisque $Z_2'$ contient un translaté de $V_1$.

\bigskip

Les deux variétés $Z_1$ et $Z_2$ sont donc des hypersurfaces, qui contiennent toutes deux un translaté de $V$, de codimension $2$. Quitte à translater ces deux variétés (ce qui est sans conséquences sur le degré et le stabilisateur), on suppose que $V \subset Z_1 \cap Z_2$. Il reste à comparer $Z_1$ et $Z_2$ pour finir la preuve. \\
\\
{\it Cas 1.} L'intersection $Z_1 \cap Z_2$ est de codimension $1$. Les deux hypersurfaces (irréductibles) sont donc égales. Par construction, $Z_1$ contient $V_1$ et on a:
\begin{displaymath}
\omega(V_1) \leq \mathrm{deg}(Z_1) \leq \mathrm{deg}(Z_2).
\end{displaymath}
En outre, l'égalité des variétés nous montre que $l_2$ est premier à la partie discrète du stabilisateur de $Z_2$, et comme cette hypersurface n'est pas incluse dans un translaté de variété abélienne (puisque cette propriété est vraie pour $V \subset Z_2$):
\begin{displaymath}
\omega(V_1) \leq \Delta^{-\rho_2}l_2^{(2-g)n_{\alpha_2}} \lambda_2(Z_2) \omega(V_1) \leq \Delta^{-\rho_2} \omega(V_1).
\end{displaymath}
On obtient donc une contradiction.\\
\\
{\it Cas 2.} L'intersection $Z_1 \cap Z_2$ est de codimension $2$. Dans ce cas, cette intersection contient $V$, mais elle contient aussi les translatés de $V$ par les points de $\Sigma_0 \cap \mathrm{Stab}(Z_2)$. Comme ce groupe est de cardinal une puissance de $l_1$, la partie discrète du stabilisateur de $V$ n'intervient pas et on a:
\begin{displaymath}
\mathrm{deg} \Big(\bigcup_{x \in \Sigma_0 \cap \mathrm{Stab}(Z_2)} x+V \Big) \geq \frac{|\Sigma_0 \cap \mathrm{Stab}(Z_2)|}{l_1^{n_{\alpha_1}(\mathrm{dim}(V)-1)}}\mathrm{deg}(V).
\end{displaymath}
On a utilisé au passage le fait que $V$ n'était pas un translaté de variété abélienne. Par le théorème de Bézout, il vient:
\begin{eqnarray*}
\frac{\lambda_1(Z_1)l_1^{(3-g)n_{\alpha_1}}}{f(\Sigma_0, Z_2)}\mathrm{deg}(V) & \leq & \mathrm{deg}(Z_1) \mathrm{deg}(Z_2) \\
& \leq & \Delta^{-\rho_1}\frac{\lambda_1(Z_1)l_1^{(1-g)n_{\alpha_1}}}{f(\Sigma_0, Z_2)}l_2^4 \omega(l_1^{n_{\alpha_1}},V)\omega(V_1).\\
\end{eqnarray*}
La majoration des termes en $l_2$ a été grossière car ceux-ci sont négligeables devant $\Delta^{\rho_1}$ par la majoration de $l_2$ suivant la proposition  \ref{quasi injectivite}. On en déduit:
\begin{displaymath}
\mathrm{deg}(V) \leq \Delta^{-\rho_1} l_1^{-2n_{\alpha_1}} l_2^{4}\omega(l_1^{n_{\alpha_1}},V) \omega(V_1).
\end{displaymath}
En raffinant (\ref{10000}) avec $n_{\alpha_1}$, on trouve:
\begin{displaymath}
l_1^{n_{\alpha_1}} \mathrm{deg}(V) \leq \Delta^{- \rho_1} l_2^4 \omega(l_1^{n_{\alpha_1}},V)^2 \leq \Delta^{u} \omega(l_1^{n_{\alpha_1}}, V)^2,
\end{displaymath}
et le réel $u$ vérifie:
\begin{displaymath}
u \leq -\rho_1+8 \rho_2(2r)^{r+1}<0.
\end{displaymath}
C'est à nouveau une contradiction. \\
\\
{\it Conclusion.} On a donc démontré par l'absurde que la proposition \ref{injectivite} était vraie, soit avec $X=V$, soit avec $X=V_1$. Il en résulte dans les deux cas que $V$ contredit la majoration (\ref{1102}). On en déduit: 
\begin{displaymath}
\hat{\mu}_{\mathrm{ess}}(V) \geq \frac{C(A)}{\omega(V)} \times
\Big(\mathrm{log}(3 \mathrm{deg}(V))\Big)^{-\lambda(r)},
\end{displaymath}
où $\lambda(r)=(16(2r)^{(r+1)})^{r}$ et $C(A)= \frac{1}{C_0^{2 \lambda(r)}}$, qui ne dépend que de $A$.
{\begin{flushright}$\Box$\end{flushright}} 

\bigskip

\noindent
\textbf{Remarques} Dans le cas où $A$ est à multiplication complexe, l'existence d'un relèvement du morphisme de Frobenius pour presque toute place de $K$ permet de démontrer la même minoration pour une variété $V$ de codimension $r$ quelconque. En effet, la fin de la preuve, à partir du lemme combinatoire, est alors (à quelques détails près) la même que dans le cas torique, en remplaçant l'isogénie $[l]$ par le relèvement du Frobenius associé (on est dans ce cas assuré d'avoir une densité positive de premiers ordinaires).

\bigskip

Il est envisageable qu'un raffinement de l'argument de descente, couplé aux nouvelles idées introduites par Amoroso dans \cite{Amoroso07}, permette de traiter la codimension quelconque.   

\bibliographystyle{fralpha}
\bibliography{biblio}

\end{document}